# FREELY INDECOMPOSABLE GROUPS ACTING ON HYPERBOLIC SPACES

ILYA KAPOVICH AND RICHARD WEIDMANN

ABSTRACT. We obtain a number of finiteness results for groups acting on Gromov-hyperbolic spaces. In particular we show that a torsion-free locally quasiconvex hyperbolic group has only finitely many conjugacy classes of $n$-generated one-ended subgroups.

## 1. INTRODUCTION

It is a well-known intuitive fact that many and in a sense most finitely generated subgroups of a given word-hyperbolic group are free. A similar informal statement can be made about groups acting on $\delta$-hyperbolic spaces. However, there are few precise finiteness results of this kind. A classic example is the following theorem which was stated by M.Gromov in [24] and proved by T.Delzant in [20]

**Theorem 1.1.** *Let $G$ be a torsion-free word-hyperbolic group. Then $G$ contains only finitely many conjugacy classes of one-ended two-generated subgroups.*

Unfortunately, the above statement is no longer true for three-generated subgroups. Consider the group
$$G = \langle a, b, t | t^{-1}at = ab^3a, t^{-1}bt = ba^3b \rangle.$$
Thus $G$ is an ascending HNN-extension of the free group $F(a,b)$ via an injective endomorphism $a \mapsto ab^3a, b \mapsto ba^3b$. It can be shown using the Combination Theorem of M.Bestvina and M.Feighn [8, 9] that $G$ is word-hyperbolic. For every integer $n \geq 1$ consider the subgroup $H_n = \langle a, b, t^n \rangle \leq G$. It is not hard to see that each $H_n$ is a normal subgroup of index $n$ in $G$ and thus $H_n$ is non-elementary and freely indecomposable. Each $H_n$ is obviously generated by three elements. Since all the subgroups $H_n$ are normal in $G$ (for $n \geq 1$), these subgroups are pair-wise non-conjugate and provide infinitely many conjugacy classes of three-generated one-ended subgroups in $G$. Similar examples can be easily obtained using the Rips construction [37].

In this paper we show that the reason for the failure of Theorem 1.1 for three-generated subgroups in the above example is the fact that $G$ possesses two-generated non-quasiconvex subgroups. Namely, the subgroup $F(a,b) \leq G$ is obviously exponentially distorted and thus non-quasiconvex.

Before formulating our main result let us recall the notion of Nielsen equivalence:

2000 *Mathematics Subject Classification.* 20F67.
The first author was supported by the U.S.-Israel Binational Science Foundation grant BSF-1999298.

**Definition 1.2** (Nielsen equivalence). *Let $G$ be a group and let $M = (g_1, \ldots, g_n) \in G^n$ be an $n$-tuple of elements of $G$. The following moves are called* elementary Nielsen moves *on $M$:*

(N1) *For some $i, 1 \le i \le n$ replace $g_i$ by $g_i^{-1}$ in $M$.*
(N2) *For some $i \ne j$, $1 \le i, j \le n$ replace $g_i$ by $g_i g_j$ in $M$.*
(N3) *For some $i \ne j$, $1 \le i, j \le n$ interchange $g_i$ and $g_j$ in $M$.*

*We say that $M = (g_1, \ldots, g_n) \in G^n$ and $M' = (f_1, \ldots, f_n) \in G^n$ are* Nielsen-equivalent, *denoted $M \sim_N M'$, if there is a chain of elementary Nielsen moves which transforms $M$ to $M'$.*

It is easy to see that if $M \sim_N M'$ then $M$ and $M'$ generate the same subgroup of $G$. For this reason Nielsen equivalence is a very useful tool for studying the subgroup structure of various groups.

We will say that a group $G$ is *almost torsion-free* if any non-trivial element of finite order in $G$ has finite centralizer. It is obvious that torsion-free groups are almost torsion-free. Moreover, almost torsion-free word-hyperbolic groups behave very similarly to torsion-free ones.

In the present paper we are able to obtain the following generalization of Theorem 1.1.

**Theorem 1.3.** *Let $G$ be an almost torsion-free word-hyperbolic group with finite generating set $S$ and denote the word metric with respect to $S$ by $d_S$. Let $l \ge 1$ be such that all $l$-generated subgroups of $G$ are quasiconvex. Let $k \ge l + 1$. Then there exists a constant $C = C(G, S, k, l)$ with the following property.*

*Suppose that $g_1, \ldots, g_k \in G$ such that $U = \langle g_1, \ldots, g_k \rangle \le G$ is freely indecomposable and not infinite cyclic. Then $(g_1, \ldots, g_k)$ is Nielsen-equivalent to a tuple $(f_1, \ldots, f_k)$ such that*

$$d_S(1, g^{-1} f_i g) \le C \text{ for } i = 1, \ldots, l+1.$$

*for some $g \in G$.*

Thus if all $l$-generated subgroups of an almost torsion-free hyperbolic group are quasiconvex, for a one-ended subgroup one can make the length of the first $l+1$ generators short by Nielsen transformations and a conjugation. This statement immediately implies the following finiteness results.

**Corollary 1.4.** *Let $G$ be an almost torsion-free word-hyperbolic group such that all $l$-generated subgroups of $G$ are quasiconvex (where $l \ge 1$). Then $G$ has only finitely many conjugacy classes of freely indecomposable non-elementary $(l+1)$-generated subgroups.*

**Corollary 1.5.** *Let $G$ be an almost torsion-free word-hyperbolic group which is* locally quasiconvex *(that is all finitely generated subgroups of $G$ are quasiconvex). Then for any $l \ge 1$ the group $G$ has only finitely many conjugacy classes of non-elementary freely indecomposable $l$-generated subgroups.*

It is worth noting that by the results of Z.Sela [39] and T.Delzant [22] for a fixed one-ended finitely presented group $H$ and a torsion-free word-hyperbolic group $G$ there are only finitely many conjugacy classes of subgroups of $G$ isomorphic to $H$. In [21] T.Delzant also shows that for a given finitely presented group $H$ and word hyperbolic group $G$ there are only finitely many conjugacy classes of subgroups of $G$ which are homomorphic images of $H$ such that the homomorphism does not factor through a group with more than one end.

The class of locally quasiconvex hyperbolic groups is very rich and includes many interesting combinatorial and geometric examples. Free groups and hyperbolic surface groups are locally quasiconvex. If $M^3$ is a closed hyperbolic 3-manifold with nonempty convex boundary, then the fundamental group of $M^3$ is word-hyperbolic and locally quasiconvex as observed by G.Swarup [40]. Moreover, J.McCammond and D.Wise [35] recently showed that "most" small cancellation groups and one-relator groups are locally quasiconvex as well. In particular, J.McCammond and D.Wise proved [35] (see also C.Hruska-D.Wise [26]) that for any non-trivial freely reduced word $w \in F(X)$ (where $F(X)$ is a free group of finite rank) there exists an integer $r_0$ such that for any $r \geq r_0$ the one-relator group $G = \langle X | w^r = 1 \rangle$ is word-hyperbolic and locally quasiconvex. It is well-known that centralizers in one-relator groups with torsion are cyclic and so Theorem 1.4 applies to these groups. P.Schupp [38] and I.Kapovich-P.Schupp [30] recently showed that many Coxeter groups of extra large type are word-hyperbolic locally quasiconvex or $k$-quasiconvex.

G.Arzhantseva [2] constructed large "generic" classes of hyperbolic groups where all $n$-generated subgroups are free and quasiconvex (see also the work of A.Yu.Olshanskii and G.Arzhantseva [3] and of I.Bumagina [15]). Since every 1-generated (i.e. cyclic) subgroup in an arbitrary hyperbolic group is quasiconvex, Theorem 1.4 implies Theorem 1.1.

Our proof of Theorem 1.3 shows that the constant $C(G, S, k, l)$ is effectively computable. This implies that the rank problem for torsion-free locally quasiconvex groups is decidable.

**Theorem 1.6.** *There is an algorithm $\mathcal{A}$ with the following property.*

*Suppose $G = \langle x_1, \ldots, x_l | r_1, \ldots, r_m \rangle$ is a finite group presentation which is known to define an one-ended almost torsion-free word-hyperbolic group where all $(l-1)$-generated subgroups are quasiconvex (e.g. $G$ is locally quasiconvex). Then $\mathcal{A}$ produces the smallest number $k$ such that $G$ can be generated by $k$ elements.*

*Proof.* First we set $S = \{x_1, \ldots, x_l\}$ and apply the algorithm of P.Papasoglu [36] to compute the hyperbolicity constant $\delta$ of the Cayley graph $\Gamma(G, S)$ with respect to the word metric $d_S$. Once $\delta$ is known, we can find the automatic structure with unique representatives on $G$ with the automatic language $L$ consisting of $S$-geodesic words [23]. Obviously $G$ can be generated by $l$ elements. Since by assumption $G$ is one-ended, it cannot be 1-generated. We now have to decide for each $k = 2, \ldots, l-1$ if $G$ can be generated by $k$ elements. For a given $k, 2 \leq k \leq l-1$ we compute $C_1 = C(G, S, k, k-1)$, the constant provided by Theorem 1.3. By Theorem 1.3 any $k$-generated subgroup of $G$ is conjugate to a subgroup generated by a $k$-element subset of the ball $B$ of radius $C_1$ in $\Gamma(G, S)$ around $1 \in G$. The number of such $k$-element subsets of $B$ is at most $(l^k)^k = l^{k^2} \leq l^{l^2}$. Since $k \leq l-1$, by our assumption on $G$ each $k$-element subset $Y$ of $B$ generates a quasiconvex subgroup of $G$ which is therefore rational with respect to the automatic language $L$. Therefore by the result of I.Kapovich [27] for each $Y$ we can recover both the quasiconvexity constant of the subgroup $H = \langle Y \rangle$ generated by $Y$ and the regular language $L_H$ which is the preimage of $H$ in $L$. Now to see whether $H = L$ it suffices to check if the regular languages $L$ and $L_H$ coincide. □

It is worth noting that the rank problem is undecidable for torsion-free word hyperbolic groups in general, as shown by G.Baumslag, C.F.Miller and H.Short [6]. Not surprisingly, their proof utilizes the Rips construction and thereby implicitly the existence of two-generated non-quasiconvex subgroups.

Our main technical tools involve generalizing the Nielsen method for groups acting on simplicial trees as developed by R.Weidmann in [43] to groups acting on $\delta$-hyperbolic spaces. The present paper substantially relies on the authors' previous paper [32], where the basic ideas are developed and fleshed out in detail. In order to obtain the above results, we first prove a technical statement related to group actions on δ-hyperbolic spaces. This statement is a close analogue of the main result of [43]. This statement turns out to also have applications to 3-manifold groups and group actions on real and simplicial trees [33].

**Theorem 1.7.**
*(1) Let $M$ be a compact hyperbolic 3-manifold with nonempty convex boundary and let $G = \pi_1(M)$. Then for any $k \geq 2$ the group $G$ has only finitely many conjugacy classes of non-elementary freely indecomposable k-generated subgroups.*
*(2) Let $M$ be a closed hyperbolic 3-manifold which fibers over a circle and let $G = \pi_1(M)$. Suppose all finitely generated subgroups of $G$ are topologically tame.*
  *Then for any $k \geq 2$ the group $G$ has only finitely many conjugacy classes of non-elementary freely indecomposable k-generated subgroups of infinite index in $G$.*

Part (1) immediately follows from Theorem 1.4 since $G$ is torsion-free and locally quasiconvex [40]. In part (2) $G$ has many non-quasiconvex subgroups (e.g the group of the fiber surface) but we are still able to use our main technical result to obtain the desired conclusion. Note that the "infinite index" assumption is essential and cannot be dropped. Indeed, if $H$ is a fiber group in $M$ then $G$ has the HNN-presentation

$$G = \langle H, t | t^{-1}ht = \phi(h), \text{ for all } h \in H \rangle,$$

where $\phi$ is some automorphism of $H$. Let $m$ be the rank of $H$. For each $n > 1$ the subgroups $H_n = \langle H, t^n \rangle \leq G$ is normal and has index $n$ in $G$. The subgroups $H_n$ are pair-wise non-conjugate and $(m+1)$-generated. On the other hand the "tameness" assumption may well be superfluous. Indeed, according to a long-standing conjecture of W.Thurston (see for example [29]) for a closed hyperbolic 3-manifold fibering over a circle all finitely generated subgroups of the fundamental group are always topologically tame.

## 2. THE MAIN TECHNICAL RESULT

Our main tool is a technical result motivated by the Kurosh subgroup theorem (see [34, 5]) for free products, which states that a subgroup is the free product of a free group and subgroups that are conjugate to subgroups of the factors. The main result of [43] generalizes this theorem to groups acting on simplicial trees. In the present paper we provide the hyperbolic "quasification" of this fact.

First, we equip every non-trivial subgroup $U \leq Isom(X)$ with a $U$-invariant quasiconvex subset $X(U)$ (see Section 4 for the precise definition). The set $X(U)$ generalizes the definition of the subtree $T_U$ in the simplicial tree action case [43] (in that situation $T_U$ contains the minimal $U$-invariant subtree as well the edges of the ambient tree which are moved a "small distance" by some nontrivial element of $U$). For an infinite subgroup $U$ of an almost torsion-free word-hyperbolic group $G$ the set $X(U)$ turns out to be Hausdorff-close to the convex hull of the limit set of $U$ in $\partial G$.

We generalize and push further the Nielsen methods for groups acting on hyperbolic spaces introduced in [32]. The objects which correspond to the tuples of elements of $G$ are the *G-tuples*:

**Definition 2.1** (*G*-tuple). *Let $G$ be a group.*
  *Let $n \geq 0$, $m \geq 0$ be integers such that $m + n > 0$. We will say that a tuple $M = (U_1, \ldots, U_n; H)$ is a G-tuple if $U_i$ is a non-trivial subgroup of $G$ for $i \in \{1, \ldots, n\}$ and $H =$* $(h_1, \ldots, h_m) \in G^m$ is a m-tuple of elements of $G$. We will denote $\overline{M} = U_1 \cup \cdots \cup U_n \cup \{h_1, \ldots, h_m\}$ and call $\overline{M}$ the underlying set of $M$. Note that $\overline{M}$ is nonempty since $m + n > 0$.

By analogy with the Kurosh subgroup theorem we will sometimes refer to the subgroups $U_i$ as *elliptic subgroups* of $M$. This is justified since in most applications the subgroups $U_i$ are generated by sets of elements with short translation length. We will also refer to $H$ as the *hyperbolic component* of $M$. We have the following notion of equivalence for $G$-tuples which generalizes the classical Nielsen equivalence.

**Definition 2.2** (Equivalence of $G$-tuples). *We will say that two $G$-tuples $M = (U_1, \ldots, U_n; H)$ and $M' = (U'_1, \ldots, U'_n; H')$ are equivalent if $H = (h_1, \ldots, h_m)$ and $H' = (h'_1, \ldots, h'_m)$ and $M'$ can be obtained from $M$ by a chain of moves of the following type:*

(1) *For some $1 \leq j \leq n$ replace $U_j$ by $gU_jg^{-1}$ where*
$$g \in \langle \{h_1, \ldots, h_m\} \cup U_1 \cup \ldots \cup U_{j-1} \cup U_{j+1} \cup \ldots \cup U_n \rangle.$$

(2) *For some $1 \leq i \leq n$ replace $h_i$ by $h'_i = g_1 h_i g_2$ where*
$$g_1, g_2 \in \langle \{h_1, \ldots, h_{i-1}, h_{i+1}, \ldots, h_m\} \cup U_1 \cup \ldots \cup U_n \rangle.$$

**Definition 2.3** ($G$-space). *Let $G$ be a group acting on a metric space $(X, d)$ by isometries. In this case we will term $(X, d)$ together with this action a $G$-space.*

We can now formulate the main technical result of this paper:

**Theorem 2.4.** *For any integer $k \geq 1$ there is a constant $K = K(k)$ with the following property.*

*Suppose $(X, d)$ is a $\delta$-hyperbolic strongly geodesic $G$-space. Let $M = (U_1, \ldots, U_n; H)$ be a $G$-tuple where $H = (h_1, \ldots, h_m)$, $n + m \leq k$ and*

$$U = \langle \overline{M} \rangle = \langle \{h_1, \ldots, h_m\} \cup \bigcup_{i=1}^n U_i \rangle \leq G.$$

*Then either $U = U_1 * \cdots * U_n * F(H)$ or after replacing $M$ by an equivalent $G$-tuple $M' = (U'_1, \ldots, U'_n; H')$ with $H' = (h'_1, \ldots, h'_m)$ one of the following holds:*

(1) $d(X(U'_i), X(U'_j)) < \delta K$ *for some* $1 \leq i < j \leq n$.
(2) $d(X(U'_i), h_j X(U'_i)) < \delta K$ *for some* $i \in \{1, \ldots, n\}$ *and* $j \in \{1, \ldots, m\}$.
(3) *There exists a* $x \in X$ *such that* $d(x, h_j x) < \delta K$ *for some* $j \in \{1, \ldots, m\}$.

The notion of "strongly geodesic" hyperbolic space (see Definition 3.4) in the above theorem is a technical condition which insures that all possible definitions of the boundary of $X$ coincide and that any two points in $X \cup \partial X$ can be connected by a geodesic. In particular, proper hyperbolic spaces, $\mathbb{R}$-trees and complete $CAT(-1)$-spaces are strongly geodesic. See also Remark 9.2.

Theorem 2.4 is proved in three steps. We first prove it for the case when $m = 0$, i.e. that $H = (-)$ is the empty tuple. Here the proof of Theorem 2.4 reduces to an elaborate version of the ping-pong argument. Then we prove it for the case $M = (U_1; H)$. Here we substantially use techniques and ideas from our previous paper [32] where the "elliptic" subgroup $U_1$ was trivial. Finally we combine these two results to handle the general case.

Theorem 2.4 is used to obtain all of our applications. Basically we use it as a "generator transfer" trick to analyze a freely indecomposable subgroup generated by a finite set $Y = \{y_1, \ldots, y_m\}$ with $m$ elements. First we start with a $G$-tuple $M_1 = (; H_Y)$ where $H_Y = (y_1, \ldots, y_m)$. We then construct a sequence of $G$-tuples $M_1, M_2, \ldots$ by repeatedly applying Theorem 2.4 in order to either "drag" elements of the "hyperbolic" component of $M_i$ into the "elliptic" components of our $G$-tuples or to join to elliptic subgroups to one new elliptic subgroup. A simple observation shows that the length of the sequence $M_1, M_2, \ldots$ is bounded by $2m - 1$. The desired results are then obtained by analyzing the terminal member of this sequence.

## 3. Hyperbolic metric spaces

We will give only a quick overview of the main definitions related to Gromov-hyperbolic spaces. For the detailed background information the reader is referred to [24], [19], [25], [1], [17], [7].

A *geodesic segment* in a metric space $(X, d)$ is an isometric embedding $\gamma : [0, s] \longrightarrow X$, where $[0, s] \subseteq \mathbb{R}$. Similarly, a *geodesic ray* and a *biinfinite geodesic* are defined as isometric embeddings $\gamma : [0, \infty) \longrightarrow X$ and $\gamma : (-\infty, \infty) \longrightarrow X$ accordingly. A metric space $(X, d)$ is said to be *geodesic* if any two points in $X$ can be joined by a geodesic segment. We will often denote a geodesic segment connecting $a \in X$ to $b \in X$ by $[a, b]$ and identify this geodesic segment with its image in $X$.

Two sets $A, B \subseteq X$ are said to be *$K$-Hausdorff close* if $A$ is contained in the $K$-neighborhood of $B$ and $B$ is contained in the $K$-neighborhood of $A$. If $A, B \subseteq X$ are $K$-Hausdorff close for some $K \geq 0$, they are said to be *Hausdorff close*. Two paths $\gamma : I \longrightarrow X$ and $\gamma' : J \longrightarrow X$ are said to be $K$-Hausdorff close (where $I, J$ are sub-intervals of the real line) if their images $\gamma(I)$ and $\gamma'(J)$ are $K$-Hausdorff close. Similarly, $\gamma, \gamma'$ are said to be Hausdorff-close if they are $K$-Hausdorff close for some $K \geq 0$.

Two geodesic rays $\gamma : [0, \infty) \longrightarrow X$ and $\gamma' : [0, \infty) \longrightarrow X$ are called *asymptotic* if there is $C \geq 0$ such that for any $t \geq 0$ we have $d(\gamma(t), \gamma'(t)) \leq C$. It is easy to see that $\gamma$ and $\gamma'$ are asymptotic if and only if they are Hausdorff close.

We recall the definition of a $\delta$-hyperbolic metric space.

**Definition 3.1** (Hyperbolic space). *A geodesic metric space $(X, d)$ is said to be $\delta$-hyperbolic if for any geodesic triangle in $X$ with geodesic sides $\alpha = [x, y], \beta = [x, z], \gamma = [y, z]$ and vertices $x, y, z \in X$ each side of the triangle is contained in the $\delta$-neighborhood of the union of two other sides. That is for any $p \in \alpha$ there is $q \in \beta \cup \gamma$ such that $d(p, q) \leq \delta$. A geodesic metric space is called* hyperbolic *if it is $\delta$-hyperbolic for some $\delta \geq 0$.*

If $(X, d)$ is a metric space and $x, y, z \in X$, one defines the *Gromov product* $(y, z)_x$ as

$$(y, z)_x := \frac{1}{2}(d(y, x) + d(z, x) - d(y, z)).$$

If $(X, d)$ is a $\delta$-hyperbolic space then the Gromov product measures for how long two geodesics stay close together. Namely, the initial segments of length $(y, z)_x$ of any two geodesics $[x, y]$ and $[x, z]$ in $X$ are $2\delta$-Hausdorff close.

One can attach to a hyperbolic space $X$ a "space at infinity", called the boundary of $X$. The boundary is usually defined as a set of equivalence classes of geodesic rays originating at a base-point $x \in X$. Two geodesic rays $\gamma, \gamma'$ (not necessarily starting at the same point) are said to be *equivalent*, denoted $\gamma \sim \gamma'$, if they are asymptotic. An equivalence class of $\gamma$ is denoted $[\gamma]$.

**Definition 3.2** (Geodesic Boundary). *Let $(X, d)$ be a $\delta$-hyperbolic geodesic metric space and let $x \in X$ be a base-point of $X$. One then defines the* boundary of $X$ relative $x$, *denoted $\partial_x X$, as the set of equivalence classes of geodesic rays originating at $x$. Similarly the* boundary $\partial X$ *of $X$ is defined as the set of equivalence classes of geodesic rays in $X$. For any $x \in X$ there is an obvious map $i_x : \partial_x X \to \partial X$ which sends $[\gamma]$ to $[\gamma]$ for any geodesic ray $\gamma$ starting at $x$.*

*The relative boundary $\partial_x X$ is topologized by saying that two points $[\gamma], [\gamma'] \in \partial_x X$ (where $\gamma, \gamma'$ are geodesic rays starting at $x$) are close if for a large $T > 0$ the segments $\gamma([0,T])$ and $\gamma'([0,T])$ are $2\delta$-Hausdorff close. More precisely, for each $p \in \partial_x X$ we set the basis of neighborhoods of $p$ in $\partial_x X$ to be the collection $\{V_x(p,r) | r > 0\}$ where*

$$V_x(p,r) = \{q \in \partial_x X | \text{ for some geodesic rays } \gamma_p, \gamma_q \text{ with } \gamma_p(0) = \gamma_q(0) = x,$$
$$[\gamma_p] = p, [\gamma_q] = q, \text{ such that the geodesic segments } \gamma_p|_{[0,r]}, \gamma_q|_{[0,r]}$$
$$\text{are } 2\delta - \text{Hausdorff close}\}.$$

An alternative and more general notion of boundary can be given in terms of sequences of points of $X$. If $(X,d)$ is a $\delta$-hyperbolic space and $x \in X$ is a base-point, we say that a sequence $(x_n)_{n=1}^\infty$ of points of $X$ *converges at infinity* if

$$\lim_{n,m \to \infty} \inf (x_n, x_m)_x = \infty.$$

This definition is easily seen to be independent on the choice of $x \in X$. Two sequences $(x_n), (y_n)$ converging at infinity are said to be *equivalent*, denoted $(x_n) \sim (y_n)$, if

$$\lim_{n,m \to \infty} \inf (x_n, y_m)_x = \infty.$$

Once again, the notion of equivalence does not depend on the choice of a base-point $x \in X$.

**Definition 3.3** (Sequential Boundary). *Let $(X,d)$ be a $\delta$-hyperbolic space. We define the sequential boundary $\partial^s X$ as the set of equivalence classes of sequences converging at infinity.*

Let $x \in X$ be a base-point. Then there is a map $j_x : \partial_x X \to \partial^s X$ defined as $j([\gamma]) := [(\gamma(n))_{n=1}^\infty]$. In good situations this map is a bijection for any $x \in X$.

In the remainder of this article we will require the hyperbolic spaces not only to be geodesic but to be *strongly geodesic*. This ensures that not only points in $X$ but also points "at infinity" can be connected by geodesics with points of $X$ and other points "at infinity" and guarantees that sequential and geodesic boundaries of $X$ coincide.

**Definition 3.4** (Strongly geodesic). *We say that a geodesic metric space $(X,d)$ is strongly geodesic if the following two conditions are satisfied:*

(1) *Let $\gamma : [0, \infty) \longrightarrow X$ be a geodesic ray in $X$ and let $y \in X$. Then there exists a geodesic ray $\gamma' : [0, \infty) \longrightarrow X$ with $\gamma'(0) = y$ such that $\gamma$ and $\gamma'$ are asymptotic.*
(2) *Let $\gamma_1 : [0, \infty) \longrightarrow X$ and $\gamma_2 : [0, \infty) \longrightarrow X$ be geodesic rays such that $\gamma_1(0) = \gamma_2(0)$ and such that $\gamma_1$ and $\gamma_2$ are not asymptotic. Then there exists a biinfinite geodesic $\gamma : (-\infty, \infty) \longrightarrow X$ such that $\gamma$ is Hausdorff close to $\gamma_1 \cup \gamma_2$.*
(3) *For any $x \in X$ and any sequence $(x_n)$ converging at infinity there exists a geodesic ray $\gamma : [0, \infty) \to X$ with $\gamma(0) = x$ and such that the sequences $(\gamma(n))_{n \geq 1}$ and $(x_n)$ are equivalent.*

The following straightforward proposition summarizes some good properties of strongly geodesic hyperbolic spaces.

**Proposition 3.5.** *Let $(X, d)$ be a strongly geodesic $\delta$-hyperbolic metric space. Then:*

*(1) For any $x \in X$ the maps $j_x : \partial_x X \to \partial^s X$ and $i_x : \partial_x X \to \partial^s X$ are bijections.*

*(2) For any $x, y \in X$ the maps $j_x : \partial_x X \to \partial^s X$, $j_y : \partial_y X \to \partial^s X$ induce the same topology on $\partial^s X$.*

*(3) For any $x, y \in X$ the maps $j_x : \partial_x X \to \partial X$, $j_y : \partial_y X \to \partial X$ induce the same topology on $\partial X$.*

Because of Proposition 3.5 for a strongly geodesic hyperbolic space $(X, d)$ we can identify both the relative geodesic boundaries $\partial_x X$ and the sequential boundary $\partial^s X$ with the full geodesic boundary $\partial X$, which inherits a canonical topology. Moreover, $X \cup \partial X$ also has a natural topological structure in this case. For points of $X$ we use the original topology of $(X, d)$. After choosing a base-point $x \in X$ for an arbitrary point $p \in \partial X$ we define the basis of neighborhoods of $p$ in $X \cup \partial X$ to be the collection $\{U_x(p, r) | r > 0\}$ where

$$U_x(p, r) = V_x(p, r) \cup \{y \in X | \text{ for some geodesic ray } \gamma, \text{ with } \gamma(0) = x, [\gamma] = p,$$
$$\text{we have } \liminf_{t \to \infty} (y, \gamma(t))_x \geq r\}.$$

It is easy to see that for a strongly geodesic hyperbolic space $(X, d)$ the topology on $X \cup \partial X$ does not depend on the choice of a base-point $x \in X$. Moreover, if $p = [\gamma] \in X$ and $x_n \in X$ is a sequence of points in $X$, then $\lim_{n \to \infty} x_n = p$ in $X \cup \partial X$ if and only if the sequences $(x_n)$ and $(\gamma(n))_{n \geq 1}$ are equivalent.

Recall that a metric space is called *proper* if all closed metric balls are compact. It is well known that proper hyperbolic metric spaces [19, 25, 1] and complete $CAT(-1)$-spaces [4, 14] as well as $\mathbb{R}$-trees are strongly geodesic. Moreover, for a proper hyperbolic space $X$ both $\partial X$ and $X \cup \partial X$ are compact.

Let $(X, d)$ be a strongly geodesic $\delta$-hyperbolic space. Let $p, q \in \partial X$ be two distinct points. We will say that a biinfinite geodesic $\gamma : (-\infty, \infty) \longrightarrow X$ *joins $p$ to $q$* if $\lim_{n \to \infty} \gamma(n) = q$ and $\lim_{n \to \infty} \gamma(-n) = p$. In this case we will often denote $\gamma$ by $[p, q]$. Similarly, a geodesic ray $\gamma$ is said to *join* the point $x = \gamma(0)$ with $p = [\gamma] \in \partial X$. In this case $\gamma$ is often denoted by $[x, p]$. It follows from the definitions that in a strongly geodesic hyperbolic space $X$ any two distinct points in $X \cup \partial X$ can be joined by a geodesic in $X$.

The following notion plays an important role in the theory of hyperbolic spaces.

**Definition 3.6** (Quasigeodesics). *Let $(X, d)$ be a metric space. A path $\sigma : I \longrightarrow X$ (where $I$ is an interval in the real line) is called a $(\lambda, \epsilon)$-quasigeodesic if $\sigma$ is parameterized by the arc-length and for any $s_1, s_2 \in I$ we have*

$$|s_1 - s_2| \leq \lambda d(\sigma(s_1), \sigma(s_2)) + \epsilon$$

*That is to say, for any points $x, y$ on $\sigma$ we have $d_\sigma(x, y) \leq \lambda d(x, y) + \epsilon$, where $d_\sigma(x, y)$ is the length of the $\sigma$-segment between $x$ and $y$.*

*A naturally parameterized path $\sigma$ in $X$ is called an $N$-local $(\lambda, \epsilon)$-quasigeodesic if every subsegment of $\sigma$ of length $N$ is a $(\lambda, \epsilon)$-quasigeodesic.*

It is well-known that in hyperbolic spaces local quasigeodesics are global quasigeodesics, provided the local parameter $N$ is big enough.

**Lemma 3.7.** [19, 25, 1] *For any $C > 0$ and $\delta \geq 0$ there is $K = K(C, \delta) > 0$ with the following property. Let $\alpha$ be a $K$-local $(C, C)$-quasigeodesic in a $\delta$-hyperbolic metric space $X$. Then $\alpha$ is a $(K, K)$-quasigeodesic.*

**Definition 3.8** (Quasiconvexity). *A subset $A$ of a geodesic metric space $(X, d)$ is said to be $\epsilon$-quasiconvex in $X$ (where $\epsilon \geq 0$) if any geodesic segment joining two points of $A$ is contained in the $\epsilon$-neighborhood of $A$. A subset $A \subseteq X$ is called* quasiconvex *if it is $\epsilon$-quasiconvex for some $\epsilon \geq 0$.*

Let $X$ be a strongly geodesic hyperbolic metric space and let $A \subseteq X \cup \partial X$ be a nonempty subset. If $A$ has just one element, we put $Conv(A) = A$. If $A$ contains at least two elements we define $Conv(A)$ to be the union of all geodesics joining distinct points of $A$. This set $Conv(A)$ is termed the *convex hull* of $A$ in $X$.

We shall need the following simple lemma.

**Lemma 3.9.** *Let $(X, d)$ be a strongly geodesic $\delta$-hyperbolic metric space. Then:*
  (1) *Any geodesic quadrilateral with vertices in $X$ is $2\delta$-thin (meaning that each side is contained in the $2\delta$-neighborhood of the union of the other sides) and $2\delta$-quasiconvex.*
  (2) *Any geodesic triangle in $X$ with vertices in $X \cup \partial X$ is $4\delta$-thin and $4\delta$-quasiconvex.*
  (3) *Any geodesic quadrilateral in $X$ with vertices in $X \cup \partial X$ is $8\delta$-thin and $8\delta$-quasiconvex.*
  (4) *For any subset $A \subseteq X$ the convex hull $Conv(A)$ is connected and $2\delta$-quasiconvex.*
  (5) *For any subset $A \subseteq X \cup \partial X$ containing at least one point of $X$ the convex hull $Conv(A)$ is $8\delta$-quasiconvex.*

*Proof.* This lemma is a standard hyperbolic exercise of the type that will be often left to the reader later on in this paper. We will present a complete proof in this instance for demonstration purposes.

(1) Let $\Sigma$ be a geodesic quadrilateral in $X$ with sides $[a, b]$, $[b, c]$, $[c, d]$ and $[a, d]$. Consider also a geodesic $[b, d]$. Since the triangles $[a, d] \cup [a, b] \cup [b, d]$ and $[b, c] \cup [c, d] \cup [b, d]$ are $\delta$-thin, the segment $[a, d]$ is contained in the $2\delta$-neighborhood of $[a, b] \cup [b, c] \cup [c, d]$. Thus $\Sigma$ is $2\delta$-thin, as required.

Now let $x, y$ be two points in $\Sigma$. It is clear that any segment $[p, q]$ is the side of a triangle or a quadrilateral such that all other sides lie in $\Sigma$. It follows that $[p, q]$ lies in the $2\delta$-neighborhood of $\Sigma$, in the case of a triangle from the definition of hyperbolicity and in the case of a quadrilateral from the first assertion.

(2) Let $\Delta = [a, b] \cup [b, c] \cup [a, c]$ be a geodesic triangle in $X$ where $a, b, c \in X \cup \partial X$.

If $a, b, c \in X$, the statement is obvious. Thus we only need to consider ideal triangles, that is when at least one vertex of $\Delta$ is in $\partial X$. We will consider the most complicated case, when all three vertices $a, b, c$ belong to $\partial X$ (the other cases are treated similarly).

Thus suppose $a, b, c \in \partial X$ and let $x \in [a, b]$. Since $[c, a]$ and $[b, a]$ are asymptotic, there exist a point $p \in [a, b]$ such that the ray $[p, a] \subseteq [b, a]$ is contained in the $2\delta$-neighborhood of $[a, c]$. Similar statements hold by symmetry for other vertices of $\Delta$. Choose a point $p' \in [x, a] \subseteq [b, a]$ such that $d(p', x) > 10\delta$ and $d(p', p'') \leq 2\delta$ for some $p'' \in [a, c]$. Similarly, choose a point $q' \in [x, b] \subseteq [a, b]$ such that $d(x, q') > 10\delta$ and $d(q', q'') \leq 2\delta$ for some $q'' \in [c, a]$. Finally, choose a point $z' \in [p'', c] \subseteq [a, c]$ such that there is a point $z'' \in [q'', c] \subseteq [b, c]$ with $d(z', z'') \leq 2\delta$ and such that $d(z', x) > 10\delta$.

Since the triangle $[p', q'] \cup [p', p''] \cup [p'', q']$ is $\delta$-thin and $d(p', x) > 10\delta$ and $d(p', p'') \leq 2\delta$, there is a point $x' \in [p'', q']$ such that $d(x, x') \leq \delta$. Since the triangle $[p'', q'] \cup [p'', z'] \cup [z', q']$ is $\delta$-thin, there is a point $x'' \in [a, c] \cup [z', q']$ such that $d(x', x'') \leq \delta$. If $x'' \in [a, c]$ then $d(x, x'') \leq 2\delta$, as required. Suppose now $x'' \in [z', q']$. Recall that by the choice of $z', q'$ we have $d(x'', q') > 8\delta$, $d(x'', z') > 8\delta$ and $d(z', z'') \leq 2\delta$, $d(q', q'') \leq 2\delta$. Since the quadrilateral $[z', q'] \cup [z', z''] \cup [q', q''] \cup [z'', q'']$ is $2\delta$-thin, this implies that there is a point $x''' \in [z'', q''] \subseteq [c, b]$ such that $d(x'', x''') \leq 2\delta$. Hence $x''' \in [c, b]$ and $d(x, x''') \leq 4\delta$, as required. We have established that the ideal triangle $\Delta$ is $4\delta$-thin, and so all triangles with vertices in $X \cup \partial X$ are $4\delta$-thin. This immediately implies that such triangles are $4\delta$-quasiconvex and part (2) is verified.

Part (3) follows from (2) similar to how (1) follows from the definition of hyperbolicity. Part (4) is a direct consequence of (1). Similarly, part (5) is a direct consequence of (3). □

**Definition 3.10** (Projection). *Let $A$ be a subset of a $\delta$-hyperbolic metric space $X$ and let $p \in X$. We say that a point $a_p \in A$ is* a projection *of $p$ on $A$ if for any $a' \in A$ we have $d(p, a_p) \leq d(p, a') + \delta$.*

If $\delta > 0$ then a projection always exists although it may not be unique. However, it is easy to see that in a hyperbolic space a projection on a quasiconvex set is "almost unique".

**Lemma 3.11.** *Let $A$ be an $\epsilon$-quasiconvex set in a $\delta$-hyperbolic space $X$. Let $p \in X$ and $a_p \in A$ and $a'_p \in A$ be projections of $p$ on $A$. Then*

  (1) $d(a_p, a'_p) \leq 2\epsilon + 8\delta$.
  (2) *If $a \in A$, $x \in [a, a_p]$ such that $x$ lies in the $R$-neighborhood of $[a_p, p]$
      then $d(a_p, x) \leq 2\delta + \epsilon + 2R$.*

*Proof.* (1) easily follows from (2) by considering the $\delta$-thin triangle $[p, a_p]$, $[a_p, a'_p]$, $[p, a'_p]$ and using $R = \delta$.

We will now prove (2). Suppose (2) fails. Choose a point $x$ on $[a_p, p]$ such that $d(a_p, x) > 2\delta + \epsilon + 2R$ and $d(x, y) \leq R$ for some $y \in [a_p, p]$. Then $d(a_p, y) > 2\delta + \epsilon + R$ and $d(p, y) = d(p, a_p) - d(a_p, y) < d(p, a_p) - \epsilon - 2\delta - R$. Since $A$ is $\epsilon$-quasiconvex, there is a point $a' \in A$ such that $d(x, a') \leq \epsilon$. Hence

$$d(p, a') \leq d(p, y) + d(y, x) + d(x, a') \leq d(p, y) + R + \epsilon <$$
$$< d(p, a_p) - (\epsilon + 2\delta + R) + R + \epsilon \leq d(p, a_p) - 2\delta$$

which contradicts the assumption that $a_p$ is a projection of $p$ on $A$. □

The following lemma is a simple but important consequence of Lemma 3.11.

**Lemma 3.12.** *Let $(X, d)$ be a $\delta$-hyperbolic geodesic metric space (with $\delta > 0$) and let $[x_p, x_q]$ be a geodesic segment in $X$. Let $p, q \in X$ be such that $x_p$ is a projection of $p$ on $[x_p, x_q]$ and that $x_q$ is a projection of $q$ on $[x_p, x_q]$.*
  *Then:*

  (1) *If $d(x_p, x_q) \geq 100\delta$ then $[p, x_p] \cup [x_p, x_q] \cup [x_q, q]$ is a $(1, 30\delta)$-quasigeodesic
      in $X$.*

(2) *Either the path $[p, x_p] \cup [x_p, x_q] \cup [x_q, q]$ is a $(1, 30\delta)$-quasigeodesic or there are points $x \in [x_p, p], y \in [x_q, q]$ such that $d(x_p, x) = d(x_q, y)$ and $d(x, y) \leq 100\delta$.*

*Proof.* Since the geodesic $[x_p, x_q]$ is $\delta$-quasiconvex, part (1) follows easily from Lemma 3.11. Part (1) immediately implies part (2). □

## 4. Limit sets and convex hulls

Till the end of this section, unless specified otherwise, we assume that $(X, d)$ is a strongly geodesic $\delta$-hyperbolic $G$-space. In this section we associate to any non-trivial subgroup $U \leq G$ a set $X_U$ which roughly corresponds to the minimal invariant subtree if $X$ is a tree and a set $X(U)$ which corresponds to the tree $T_U$ of [43].

Let $x \in X$. We define the *limit set* of $U$, denoted $\Lambda(U)$, to be the collection of all $p \in \partial X$ such that $p = \lim_{n \to \infty} u_n x$ for some sequence $u_n \in U$. If $y \in X$ is a different point then the orbits $Ux$ and $Uy$ are $d(x, y)$-Hausdorff close. Therefore the definition of $\Lambda(U)$ does not depend on the choice of $x \in X$. We can now define the analogue of the minimal $U$-invariant subtree:

**Definition 4.1.** *Let $U \leq G$ be a non-trivial subgroup such that $\Lambda(U)$ has at least two distinct points. Then we define the* weak convex hull $X_U$ *of $U$ as*

$$X_U := Conv(\Lambda(U)).$$

Note that the assumption on the limit set is also necessary in the case when $X$ is a real tree. Indeed a group $U$ acting on a real tree by isometries either contains a hyperbolic element and has two distinct points in the limit set or the group acts with a fixed point in which case the minimal $U$-invariant subtree is not necessarily unique.

Traditionally the weak convex hull $X_U$ is termed the "convex hull" of $U$ (see for example [31, 42]). It is a very useful and natural geometric object with many interesting applications.

However, as in [43], it turns out that $X_U$ is not the right object for our purposes and we need to consider a bigger $U$-invariant set. Namely, if $U$ acts on a simplicial tree $T$, it is necessary (see [43]) to study the set that contains not only the minimal $U$-invariant subtree but also the points that are fixed under the action of some non-trivial element of $U$. Following this analogy, we introduce the following notions:

**Definition 4.2.** *Let $U \leq G$ be a non-trivial subgroup.*

*We define the* small displacement set *of $U$, denoted $E(U)$, as*

$$E(U) := \{x \in X \mid d(x, gx) \leq 100\delta \text{ for some } g \in G, g \neq 1\}$$

*We put $Z(U) := Conv(\Lambda(U) \cup E(U))$. We now define*

$$X(U) := \overline{Conv(Z(U))}$$

*and call $X(U)$ the* convex hull *of $U$.*

The constant $100\delta$ in the above definition was chosen for computational convenience and can in fact be decreased. Recall that if $g$ is an isometry of a metric space $(X, d)$, the *translation length* $||g||$ of $g$ is defined as

$$||g|| = \inf\{d(x, gx) \mid x \in X\}.$$

An important observation is the following simple lemma:

**Lemma 4.3.** *Let $U \leq G$ be a non-trivial subgroup. Then $X(U)$ is nonempty.*

*Proof.* We will sketch the argument and leave the details to the reader. If $\delta = 0$ then $X$ is an $\mathbb{R}$-tree and the statement is obvious since any isometry of $X$ either fixes a point or acts by translation on a line in $X$ (Hence either $E(U)$ is nonempty or $\Lambda U$ has at least two points and hence $Conv(\Lambda U)$ is nonempty.

Suppose $\delta > 0$. Let $g \in G$ be an isometry of $X$ such that $||g|| \geq 50\delta$. Let $x \in X$ such that $d(x, gx) \leq ||g|| + \delta$. Thus $d(x, gx) \geq 49\delta$. It is easy to see that the bi-infinite path

$$\sigma = \cdots \cup g^{-2}[x, gx] \cup g^{-1}[x, gx] \cup [x, gx] \cup g[x, gx] \cup g^2[x, gx] \cup \ldots$$

is $49\delta$-local $(1, 4\delta)$-quasigeodesic. It then follows that $\sigma$ is a $(K, K)$-quasigeodesic for some constant $K = K(\delta) > 0$ (see for example Lemma 1.1 in [20]). Hence the sequences $(g^n x)_{n \geq 1}$ and $(g^{-n} x)_{n \geq 1}$ both converge at infinity and are not equivalent. Therefore $\Lambda \langle g \rangle$ consists of two distinct points.

Suppose now $U \leq G$ is a non-trivial subgroup. Choose a non-trivial element $u \in U$. If $||u|| \leq 50\delta$ then obviously $E(U)$ is nonempty. If $||u|| \geq 50\delta$ then $\Lambda U$ has at least two distinct points. In any event $Z(U)$ is clearly nonempty and hence $X(U)$ is nonempty as well. □

We shall also need the following simple geometric observation:

**Lemma 4.4.** *Let $U \leq G$ be a non-trivial subgroup and suppose that $\delta > 0$ (Recall that $\delta$ is the hyperbolicity constant of $X$). Then*

(1) *The sets $X(U), X_U, E(U)$ and $Z(U)$ are $U$-invariant. The set $X_U$ is $8\delta$-quasiconvex and the set $X(U)$ is connected, closed and $4\delta$-quasiconvex.*
(2) *Suppose $p \in X$, $u \in U, u \neq 1$ and let $p'$ be a projection of $p$ on $X(U)$. Then the path $[p, p'] \cup [p', up'] \cup u[p', p]$ is $(1, 100\delta)$-quasigeodesic.*
(3) *For any $g \in G$ we have $E(gUg^{-1}) = gE(U)$, $\Lambda(gUg^{-1}) = g\Lambda(U)$, $Z(gUg^{-1}) = gZ(U)$, $X(gUg^{-1}) = gX(U)$ and $X_{gUg^{-1}} = gX_U$.*

*Proof.* It is clear from Definition 4.2 that $X(U)$, $Z(U)$ and $X_U$ are $U$-invariant, since the limit set $\Lambda(U)$ and the small displacement set $E(U)$ are $U$-invariant by construction. The quasiconvexity estimates follow directly from Lemma 3.9. The set $X(U)$ is connected and closed by construction.

We will now show that (2) holds. Note that $up' \in X(U)$ since $X(U)$ is $U$-invariant and $p' \in X_G$. It is obvious that $up'$ is a projection of $gu$ on $X(U)$ and that $d(p, p') = d(up, up')$. Since $X(U)$ is $4\delta$-quasiconvex, the point $p'$ is at most $4\delta$-away from a projection $p''$ of $p$ on $[p', q']$. Similarly, $up'$ is at most $4\delta$-away from a projection of $gp$ on $[p', gp']$. Therefore by Lemma 3.12 either part (2) of Lemma 4.4 holds or there are points $x \in [p, p']$ and $y \in u[p, p']$ such that $d(x, p') = d(y, up') > \delta$ and $d(x, y) \leq 100\delta$. The choice of $y$ implies that in fact $y = ux$. Thus $d(x, ux) \leq 100\delta$ and $u \neq 1$. Therefore $x \in E(U)$ and hence $x \in X(U)$. Since $d(p', x) > \delta$, we have

$$d(p, x) = d(p, p') - d(p', x) < d(p, p') + \delta,$$

which contradicts our assumption that $p'$ is a projection of $p$ on $X(U)$.

Part (3) follows directly from the definitions. □

**Lemma 4.5.** *Let $U \leq G$ be a non-trivial subgroup and assume that $\delta > 0$. Suppose $A \subseteq X$ is a nonempty $C$-quasiconvex and $U$-invariant set such that*

FREELY INDECOMPOSABLE GROUPS ACTING ON HYPERBOLIC SPACES    13(1) *the set $A$ is contained in the $C$-neighborhood of $X(U)$; and*
(2) $E(U) \subseteq A$.

*Then $X(U)$ and $A$ are $(3C + 4\delta)$-Hausdorff close.*

*Proof.* We need to show that $X(U)$ is contained in the $(3C + 4\delta)$-neighborhood of $A$.

By condition (2) we know that $E(U) \subseteq A$. If $\Lambda U = \emptyset$, then $Z(U) = E(U)$, $X(U) = \overline{Conv(E(U))}$ and the conclusion of the lemma is obvious. Assume now that $\Lambda U \neq \emptyset$.

Choose a point $x \in A$. Since $A$ is $U$-invariant, the orbit $Ux$ is contained in $A$. Hence the collections of limits in $\partial X$ of sequences from $A$ contains $\Lambda U = \Lambda(Ux)$, the limit set of $U$.

Let $p \in \Lambda U$ be any point. Consider a geodesic ray $[x, p]$ in $X$ and choose any $y \in [x, p]$. Since $p$ is approximated by elements of $Ux$ and $X$ is $\delta$-hyperbolic, there exists $u \in U$ such that $d(y, y') \leq 2\delta$ for some $y' \in [x, ux]$. Since $A$ is $U$-invariant, $ux \in A$. Since $A$ is $C$-quasiconvex, there is $a \in A$ with $d(y', a) \leq C$. Hence $d(y, a) \leq C + 2\delta$. Since $y \in [x, p]$ was chosen arbitrarily, the geodesic ray $[x, p]$ is contained in the $(C + 2\delta)$-neighborhood of $A$.

Let $p, q \in \Lambda U$ be two distinct points and let $[p, q] \subseteq X$ be a bi-infinite geodesic from $p$ to $q$. Let $y \in [p, q]$ be an arbitrary point. Since both $p$ and $q$ are approximated by points of the orbit $Ux$, there exist $u_p, u_q \in U$ such that for some point $y' \in [u_p x, u_q x]$ we have $d(y, y') \leq 2\delta$. Since $u_p x, u_q x \in A$ and $A$ is $C$-quasiconvex, there is $a \in A$ such that $d(y', a) \leq C$ and hence $d(y, a) \leq C + 2\delta$. Since $y \in [p, q]$ was chosen arbitrarily, this implies that $[p, q]$ is contained in the $(C + 2\delta)$-neighborhood of $A$. We have already seen that $[x, p]$ and $[x, y]$ are contained in the $(C + 2\delta)$-neighborhood of $A$ as well for any $x \in A$. Since $E(U) \subseteq A$ and $A$ is $C$-quasiconvex, this implies that $Z(U) = Conv[\Lambda U \cup E(U)]$ is contained in the $(2C + 2\delta)$-neighborhood of $A$. (This is true even if $\Lambda(U)$ is a single point.)

Since $X$ is $\delta$-hyperbolic and $A$ is $C$-quasiconvex, we conclude that $Conv(Z(U))$ is contained in the $(3C + 3\delta)$-neighborhood of $A$. Hence $X(U) = \overline{Conv(Z(U))}$ is contained in the $(3C + 4\delta)$-neighborhood of $A$ and so the sets $A$ and $X(U)$ are $(3C + 4\delta)$-Hausdorff close, as required. □

In our application to word-hyperbolic groups it turns out that the set $X(U)$ and $X_U$ are Hausdorff close. To see that this is not true in general consider the following two situations:

**Example 4.6.** *(1) Consider the action of the infinite dihedral group $U = D_\infty = \langle x, y | x^2 = 1, y^2 = 1 \rangle$ on the hyperbolic plane $\mathbb{H}^2$ where $x$ and $y$ act as reflection in disjoint geodesic lines. It is clear that we get $\Lambda(U) = \{y_1, y_2\}$ for some $y_1, y_2 \in \partial \mathbb{H}^2$ and that $X_U$ is the geodesic line $[y_1, y_2]$ which is perpendicular to the axes of $x$ and $y$. However both $x$ and $y$ clearly fix points in arbitrary distance from $X_U$.*

*(2) Consider the action of $U = \mathbb{R}$ on $\mathbb{H}^2$ where $U$ acts by translation along a fixed geodesic line $L$. We clearly have $X_U = L$. However since $U$ is connected we have that $X(U) = E(U) = \mathbb{H}^2$.*

*(3) Consider the group $A = \mathbb{Z} \times \mathbb{Z}_2$ with the generating set $S = \{a, b\}$ where $a$ in the generator of $\mathbb{Z}$ and $b$ is the generator of $\mathbb{Z}_2$. Then the Cayley graph $X = \Gamma(A, S)$ is 1-hyperbolic. Then for $U = \langle b \rangle = \mathbb{Z}_2$ we have $X_U = \emptyset$ and $X(U) = E(U) = X$.*

We will describe a situation when $X_U$ and $X(U)$ are Hausdorff close.

**Definition 4.7.** *For an element $g \in G$ we define the* asymptotic translation length $||g||_\infty$ *of $g$ as*
$$||g||_\infty := \liminf_{n \to \infty} \frac{||g^n||}{n}.$$

One can show that for any $x \in X$ we have

$$||g||_\infty = \lim_{n \to \infty} \frac{d(x, g^n x)}{n} \leq ||g|| \leq d(x, gx).$$

It is clear that $\Lambda(\langle g \rangle) = \{y_1, y_2\}$ for some distinct $y_1, y_2 \in \partial X$ if $||g||_\infty > 0$ (since in this case the orbit map $\mathbb{Z} \to X, n \mapsto g^n x$ is a quasi-isometric embedding for any $x \in X$). We refer the reader to [18] for the background information on asymptotic translation length.

We can now state a criterion which allows one to rule out situations similar to (2) in Example 4.6. For any $g \in G$ we denote $E(g) = \{x \in X | d(x, gx) \leq 100\delta\}$ and call $E(g)$ the *small displacement set* of $g$.

**Lemma 4.8.** *Suppose that $\delta > 0$. Let $g \in G$ and $C = \langle g \rangle \leq G$ and suppose that $||g||_\infty = c > 0$. Then $E(g)$ is contained in the $k$-neighborhood of $X_C$ where $k = \frac{5000\delta^2}{c} + 50\delta$.*

*Proof.* We can clearly assume that $||g||_\infty \leq 100\delta$ since otherwise $E(g) = \emptyset$. Let $x \in E(g)$, i.e. $d(x, gx) \leq 100\delta$. Note that by the choice of $g$ the limit set $\Lambda C$ consists of two distinct points $y_1, y_2 \in \partial X$ where $y_1 = \lim_{n \to \infty} g^{-n} x$ and $y_2 = \lim_{n \to \infty} g^n x$. Hence $X_C$ is $8\delta$-quasiconvex and $4\delta$-Hausdorff close to a geodesic $[y_1, y_2]$.

Choose $m$ to be the smallest integer such that $m \geq \frac{100\delta}{c} = \frac{100\delta}{||g||_\infty}$. Note that $m \leq \frac{100\delta}{c} + 1$. It clearly follows immediately from the definition of the asymptotic translation length that $d(y, g^m y) \geq mc \geq 200\delta$ for all $y \in X$.

Let $p_x$ be a projection of $x$ onto $X_C$. Then $d(p_x, g^m p_x) > 100\delta$. Since $X_C$ is $8\delta$-quasiconvex, it now follows from Lemma 3.12 that the path $[x, p_x] \cup [p_x, g^m p_x] \cup [g^m p_x, g^m x]$ is a $(1, 100\delta)$-quasigeodesic. Hence $d(x, g^m x) \geq d(x, p_x) + d(p_x, g^m p_x) + d(g^m p_x, g^m x) - 100\delta \geq d(x, p_x) + 100\delta + d(g^m p_x, g^m x) - 100\delta = d(x, p_x) + d(g^m p_x, g^m x) = 2d(x_p, g^m x_p)$ and therefore $d(x, p_x) \leq \frac{1}{2} d(x, g^m x)$. Since by assumption $d(x, gx) \leq 100\delta$, it follows that $d(x, g^m x) \leq m 100\delta \leq (\frac{100\delta}{c} + 1) 100\delta = \frac{10000\delta^2}{c} + 100\delta$. It follows that

$$d(x, p_x) \leq \frac{1}{2} d(x, g^m x) \leq \frac{1}{2}(\frac{10000\delta^2}{c} + 100\delta) = \cdot \frac{5000\delta^2}{c} + 50\delta = k,$$

as required. $\square$

The following lemma shows that the situation as in part (1) of Example 4.6 does not occur in our applications to hyperbolic groups.

**Lemma 4.9.** *For any integers $n, k, \delta \geq 1$ there exists a constant $k_1 = k_1(n, k, \delta)$ with the following property. Suppose that $G = \langle S | R \rangle$ is an almost torsion-free group with $\delta$-hyperbolic Cayley-graph $X = \Gamma(G, S)$ and such that $S$ consists of $n$ elements. Then for any non-trivial element of finite order $g \in G$ the set $E(g, k) = \{x \in X | d(x, gx) \leq k\}$ has diameter at most $k_1$ in $X$.*

*Proof.* Choose an integer $L = L(\delta) > 0$ and a number $\lambda = \lambda(\delta) > 0$ such that any $L$-local geodesic in $X$ is $(\lambda, \lambda)$-quasigeodesic. We show that the assertion holds for

$$k_1 = k_1(n, k, \delta) = L + L(2n)^{k+5\delta}(2n)^L + 2 = L + L(2n)^{k+5\delta+L} + 2.$$

Assume that the statement of Lemma 4.9 fails for $G$. We will show that then there exists an element $h \in G$ represented by a geodesic word $w = w(h)$ of length at least $L$ such that $w^n$ is a $L$-local geodesic for all integers $n \geq 0$ and such that $h$ commutes with a non-trivial element of finite order in $G$. Such an element $h$ is clearly of infinite order which yields a contradiction to the assumption that centralizers of non-trivial torsion elements are finite.

Denote $k_1' = k_1 - 2$. Let $g \in G$ be a non-trivial torsion element. By the properties of the word-metric on the Cayley graph $X$ of $G$, it suffices to show that for any two vertices $x, y \in E(g, k)$ we have $d(x, y) \le k_1'$.

Suppose, on the contrary, that there exists vertices $x, y \in X = \Gamma(G, S)$ such that $d(x, y) \ge k_1'$ and that $x, y \in E(g, k)$, that is $d(x, gx) \le k$ and $d(y, gy) \le k$. Choose a geodesic $[x, y]$. For each $i = 0, 1, \ldots, (2n)^{k+5\delta+L}$ define $p_i$ to be the point of $[x, y]$ such that $d(x, p_i) = iL$. Note that $p_i$ is a vertex of $X$ since $iL$ is an integer and $x \in G$. Denote by $w_i$ the subsegment $[p_i, p_{i+1}]$ of $[x, y]$. Clearly $w_i$ has length $L$. We also choose for each $i = 0, 1, \ldots, (2n)^{k+5\delta+L}$ a geodesic path $v_i = [p_i, gp_i]$. It follows immediately from the fact that geodesic quadrilaterals in $X$ are $2\delta$-thin that $d(p_i, gp_i) \le k + 5\delta$ for all $i$.

Since $S$ has $n$ elements, for any integer $m \ge 0$ there are at most $(2n)^m$ possibilities for the label of an edge-path of length at most $m$ in $X$. Since the paths $w_i$ are of length $L$ and the paths $v_i$ are of length at most $k + 5\delta$ and since there are $1 + (2n)^{k+5\delta+L}$ choices for $i$ it follows that there exist $i, j$ with $0 \le i < j \le (2n)^{k+5\delta+L}$ such that $w_i$ has the same label as $w_j$ and $v_i$ has the same label as $v_j$. Let $w_{(i,j)}$ be the subsegment $[p_i, p_j]$ of $[x, y]$.

Note that the label of $v_i$ gives an element $g' \in G$ conjugate to $g$. Denote by $h$ the element of $G$ represented by the label of $w_{(i,j)}$. It is clear that $g'h = hg'$ since we have the quadrilateral with vertices $p_i, p_j, gp_i$ and $gp_j$ where opposite sides have the same labels corresponding to either $g'$ or $h$. Moreover, by construction the path $w_{(i,j)}^m$ representing $h^m$ is a $L$-local geodesic for any integer $m \ge 1$. Hence $h$ has infinite order. However, $h$ commutes with a non-trivial element of finite order $g'$, which yields a contradiction. □

We can now show that for a wide class of subgroups of hyperbolic groups weak convex hulls are Hausdorff-close to convex hulls.

**Lemma 4.10.** *For any integers $n \ge 1$ and $\delta \ge 1$ there exists a constant $c = c(n, \delta) \ge 0$ with the following property. Suppose that $G = \langle S | R \rangle$ with $\#S = n$ is an almost torsion-free group with a $\delta$-hyperbolic Cayley graph $X = \Gamma(G, S)$. Let $U < G$ be an infinite subgroup. Then $X_U$ and $X(U)$ are $c$-Hausdorff close.*

*Proof.* It clearly suffices to show that there exists a constant $c(n, \delta)$ such that $E(u, 100\delta) = \{x \in X | d(x, ux) \le 100\delta\}$ is contained in the $c$-neighborhood of $X_U$ for all non-trivial $u \in U$.

For non-torsion elements this follows from Lemma 4.8 since it is well-known that there exists a constant $k(n, \delta) > 0$ such that $|g|_\infty \ge k(n, \delta)$ for all non-torsion elements $g \in G$ (see for example [41]).

Suppose now that $u \in U$ is a non-trivial torsion element and choose $x \in E(u, 100\delta)$, so that $d(x, ux) \le 100\delta$. Let $k_1 = k_1(n, 300\delta, \delta)$ be the constant provided by Lemma 4.9. Since $U$ is infinite, there exits a point $p \in \Lambda U$. Since $\Lambda U$ is $U$ invariant we have $up \in \Lambda U$.

We will first observe that $p \ne up$. Indeed, suppose $up = u$. It follows that $u$ maps $[x, p]$ to $[ux, p] = [ux, up]$. Since $d(x, ux) \le 100\delta$ and since ideal triangles are $4\delta$-thin it follows that for every $z \in [x, p]$ with $d(z, x) \ge 105\delta$ there exists a $y \in [ux, p]$ such that $d(z, y) \le 4\delta$. The triangle inequality gives that $d(x, z) - 100\delta - 4\delta \le d(ux, y) \le d(x, z) + 100\delta + 4\delta$. This clear implies that $d(z, uz) \le d(z, y) + d(y, uz) \le 4\delta + |d(ux, uz) - d(ux, y)| = 4\delta + |d(x, z) - d(ux, y)| \le 4\delta + 104\delta = 108\delta$. This however clearly contradicts Lemma 4.9 since $z$ can be arbitrarily far from $x$.

Let $z \in [x, p]$ be such that $d(x, z) = k_1 + 300\delta$. By the choice of $k_1$ this implies $d(z, uz) \geq 300\delta$.

We will prove that $x$ lies in the $(k_1 + 308\delta)$-neighborhood of any geodesic $[p, up]$ and therefore in the $(k_1 + 308\delta)$-neighborhood of $X_U$. It suffices to show that $z$ lies in the $8\delta$-neighborhood of $[p, up]$. Suppose now that $z$ does not lie in the $8\delta$-neighborhood of $[p, up]$. Consider the ideal geodesic quadrilateral $[x, p] \cup [x, ux] \cup u[x, p] \cup [p, up]$. Since this quadrilateral in $8\delta$-thin, the point $z$ lies in the $8\delta$-neighborhood of $[x, ux] \cup [ux, up]$. Choose $y \in [x, ux] \cup [ux, up]$ such that $d(y, z) \leq 8\delta$. If $y \in [x, ux]$ then by the triangle inequality $d(x, ux) \leq 100\delta$ implies $d(z, x) \leq 108\delta$, contrary to the choice of $z$. Thus $z \in u[x, p]$. Since $d(x, ux) \leq 100\delta$ and $d(z, y) \leq 8\delta$, by the triangle inequality we have $|d(x, z) - d(y, uz)| \leq 108\delta$. Since $d(x, z) = d(ux, uz)$, this implies that $d(uz, y) \leq 108\delta$. It now follows from $d(z, y) \leq 8\delta$ that $d(z, uz) \leq 116\delta$, which contradicts our assumptions. □

**Lemma 4.11.** *There exists a constant $c(n, \delta) \geq 0$ such that the following holds. Let $G = \langle S | R \rangle$ be as in Lemma 4.10. Let $U \leq G$ be a non-trivial finite subgroup. Then the diameter of $X(U)$ is bounded by $c(n, \delta)$.*

*Proof.* By a theorem of V.Gerasimov and O.Bogopolskii [10] (see also [13]), any finite subgroup of $G$ is conjugate to a subgroup contained in the ball of radius $2\delta + 1$ around 1 in the Cayley graph $X = \Gamma(G, S)$. This implies that for any non-trivial finite subgroup $U \leq G$ there exists a vertex $x \in X$ such that $d(x, ux) \leq 2\delta + 1$ for all $u \in U$.

Let $k_1 = k_1(n, 102\delta + 1, \delta)$ be the constant provided by Lemma 4.9. It is clear that

$$E(U) = \bigcup_{u \in U - \{1\}} E(u, 100\delta) \subset \bigcup_{u \in U - \{1\}} E(u, 102\delta + 1).$$

Furthermore $x \in E(u, 102\delta + 1)$ for all non-trivial $u \in U$. It follows that the diameter of $E(U)$ is bounded by $2k_1$. This clearly implies that the diameter of $X(U) = \overline{Conv(E(U))}$ is bounded by $2k_1 + 10\delta$. □

5. Minimal networks

Throughout this section $(X, d)$ is a strongly geodesic 1-hyperbolic $G$-space.

**Definition 5.1** (Bridge). *For two nonempty sets $A, B$ in a $X$ we will say that $[a, b]$ is a* bridge *between $A$ and $B$ if $a \in A$, $b \in B$ and for any $a' \in A, b' \in B$ we have $d(a, b) \leq d(a', b') + 1$, that is $d(a, b) \leq d(A, B) + 1$. Since $1 > 0$, a bridge always exists but is not necessarily unique. We will sometimes denote a bridge $[a, b]$ between $A$ and $B$ by $[A, B]$.*

We will need the following simple lemma which states that if two quasiconvex subsets in a hyperbolic space are sufficiently far apart, then there is an almost unique bridge between these sets.

**Lemma 5.2.** *Let $A, B$ be nonempty 4-quasiconvex subsets of $X$. Let $a \in A$, $b \in B$ be such that $[a, b]$ is a bridge between $A$ and $B$. Suppose also that $d(a, b) \geq 100$, that $a' \in A$ and $b' \in B$. Then the following hold:*
  (1) *The path $[a', a] \cup [a, b] \cup [b, b']$ is $(1, 50)$-quasigeodesic which is $50$-Hausdorff close to any geodesic $[a', b']$.*
  (2) *If $d(a', b') \leq d(A, B) + R$ then $d(a, a') \leq R + 50$ and $d(b, b') \leq R + 50$.*

(3) Let $\gamma = [a', b']$ be a geodesic connecting $a'$ and $b'$. Then there are points $p, q \in [a', b']$ such that:

(i) $d(a, p) \leq 100$ and $d(b, q) \leq 100$;

(ii) The segments $[a', a]$ and $[a', p]$ are 100-Hausdorff close, the segments $[b', q]$ and $[b', b]$ are 100-Hausdorff close and the segments $[p, q]$ and $[a, b]$ are 100-Hausdorff close.

(iii) For any point $x \in [p, q]$ we have $|d(x, A) - d(x, p)| \leq 100$ and $|d(x, B) - d(x, q)| \leq 100$.

*Proof.* (1) Since $[a, b]$ is a bridge between $A$ and $B$, it follows that $b$ is a projection of $a$ on $B$ and that $a$ is a projection of $b$ on $A$. Hence by Lemma 3.11 if $x$ is a point on $[a, b]$ contained in the 1-neighborhood of $[a, a']$ then $d(a, x) \leq 2 + 4 + 2 = 8$. Similarly, if $y$ is a point on $[a, b]$ contained in the 1-neighborhood of $[b, b']$ then $d(b, y) \leq 8$. Now statement (1) easily follows from applying the definition of hyperbolicity since $d(a, b) \geq 100$.

(2) Suppose $a', b'$ are such that $d(a', b') \leq d(A, B) + R$. By part (1) we have
$$d(a', b') \geq d(a, b) + d(a, a') + d(b, b') - 50$$

If at least one of $d(a, a'), d(b, b')$ is greater than $R + 50$, then the above formula implies $d(a', b') > d(a, b) + R = d(A, B) + R$, which contradicts the choice of $a', b'$.

Part (3) follows from (1) and (2) and we leave the details to the reader. □

It is well known that polygons with geodesic sides in hyperbolic metric spaces look like trees. This statement easily generalizes to the case of "polygons" whose "vertices" are quasiconvex sets. Indeed, a geodesic $n$-gon in a 1-hyperbolic space is obviously $(n-2)$-thin, that is each side is contained in the $(n-2)$-neighborhood of the union of $n-1$ remaining sides. This implies that minimal networks in hyperbolic spaces are well-behaved and also look like unions of trees. We refer the reader to [24], [19] and [12] for the background information and basic facts regarding approximating trees in $\delta$-hyperbolic spaces.

If $\Gamma$ is an embedded closed finite subgraph of $(X, d)$ with geodesic edges, we can assign a natural *perimeter* $P(\Gamma)$ to $\Gamma$. Namely, $\Gamma$ can be represented as a union of finitely many geodesic segments with disjoint interiors. We put $P(\Gamma)$ to be the sum of the lengths of these segments.

**Definition 5.3** (Connector). *Let $(X, d)$ be a 1-hyperbolic geodesic metric space. Let $A_1, \ldots, A_n$ be closed connected nonempty subsets of $X$ such that $n \geq 2$ and $A_i \cap A_j = \emptyset$ for $i \neq j$. We define a* connector *for $A_1, \ldots, A_n$ as an embedded closed finite subgraph $Y$ of $X$ such that:*

(a) *The set $Y \cup A_1 \cup \cdots \cup A_n$ is connected.*

(b) *For any proper closed subgraph $\Gamma$ of $Y$ the set $\Gamma \cup A_1 \cup \cdots \cup A_n$ is disconnected.*

We further define $d(A_1, \ldots, A_n) = \inf\{P(Y) | Y \text{ is a connector for } A_1, \ldots, A_n\}$. We say that a connector $\Omega$ for $A_1, \ldots, A_n$ is a minimal connector if $P(\Omega) \leq d(A_1, \ldots, A_n) + 1$ and the closures of the component of $\Omega - (A_1 \cup \ldots \cup A_n)$ are trees such that branching points subdivide them into geodesic segments. A minimal connector will be denoted $\Omega(A_1, \ldots, A_n)$. We term the closures of the components of $\Omega - (A_1 \cup \cdots \cup A_n)$ the tree components *of* $\Omega$.

Also, if $A$ is a closed connected nonempty subset of $X$, we set $d(A) = 0$. If $n \geq 1$ and $A_1, \ldots, A_n$ are closed nonempty connected subsets of $X$ (not necessarily disjoint), we set
$$d(A_1, \ldots, A_n) := d(A'_1, \ldots, A'_p),$$

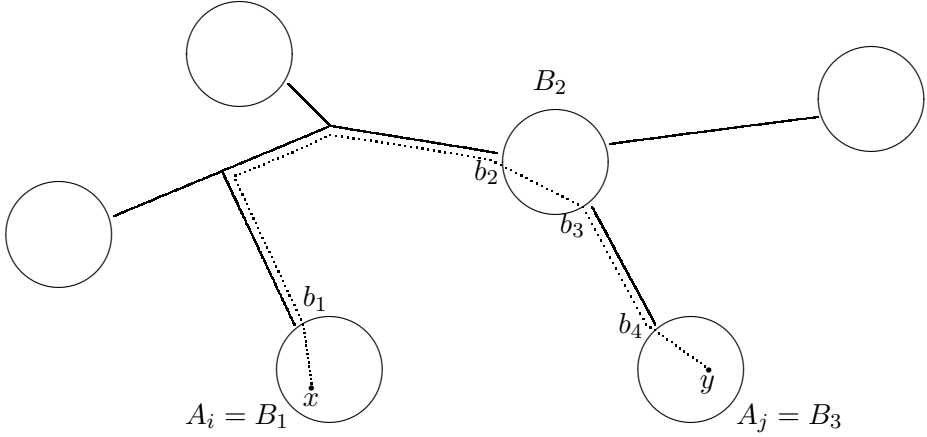

FIGURE 1. A connector $\Omega$ with 3 tree-components and an $\Omega$-geodesic joining $x \in A_i$ with $y \in A_j$.

where $A'_1, \ldots, A'_p$ are the distinct connected components of $A_1 \cup \cdots \cup A_n$.

Clearly minimal connector always exists but may not be unique. However, we will see that if $A_i$ are quasiconvex subsets sufficiently far away from each other, any two minimal connectors are Hausdorff close. Note also that for two connected sets $A_1, A_2$ the definition of a connector coincides with that of a bridge.

The following elementary lemma follows directly from the definition of the minimal connector. For points $x_1, \ldots, x_k \in X$ we denote $\Omega(\{x_1\}, \ldots, \{x_k\})$ by $\Omega(x_1, \ldots, x_k)$.

**Lemma 5.4.** *Let $A_1, \ldots, A_n$ be disjoint connected closed nonempty subsets of $X$ (where $n \geq 2$) and let $\Omega = \Omega(A_1, \ldots, A_n)$ be a minimal connector. Then the following statements hold:*

(1) *If $x \in A_i$, $y \in A_j$ are terminal vertices of a tree-component $T$ of $\Omega$, then $i \neq j$.*
(2) *The total number of terminal vertices of tree-components $T$ of $\Omega$ is at most $2(n-1)$.*
(3) *Let $T$ be a tree-component of $\Omega$ with terminal vertices $b_1, \ldots, b_k$. Then $T = \Omega(b_1, \ldots, b_k)$.*

Suppose that $A_1, \ldots, A_n$ are disjoint connected closed subsets of $X$ and $\Omega = \Omega(A_1, \ldots, A_n)$ is a minimal connector. It is clear that for any $i \neq j$ there exists a unique minimal set of tree-components $T_1, \ldots, T_k$ of $\Omega$ such that $A_i$ and $A_j$ belong to the same connected component of $T_1 \cup \cdots \cup T_k \cup A_1 \cup \cdots \cup A_n$. They define a sequence of sets $A_i = B_1, B_2, \ldots, B_{k+1} = A_j$ (where each $B_s \in \{A_1, \ldots, A_n\}$) and a sequence of points $b_1 \in B_1, b_2 \in B_2, b_3 \in B_2, b_4 \in B_3, b_5 \in B_3 \ldots, b_{2k} \in B_{k+1}$ such that for every $s = 1, \ldots, k$ the points $b_{2s-1}$ and $b_{2s}$ are terminal vertices of $T_k$.

**Definition 5.5.** *Suppose that $x \in A_i$ and $y \in A_j$. Let the sequences $(T_i), (B_i), (b_i)$ be as above.*

*We then call a path of the form*

$$\sigma = [x, b_1] \cup [b_1, b_2]_{T_1} \cup [b_2, b_3] \cup [b_3, b_4]_{T_2} \cup \cdots \cup [b_{2k-1}, b_{2k}]_{T_k} \cup [b_{2k}, y]$$

*an $\Omega$-geodesic from $x$ to $y$ and denote it by $[x, y]_\Omega$. Sub-paths of the above path will also be termed $\Omega$-geodesics. The length of $[x, y]_\Omega$ will be denoted by $d_\Omega(x, y)$.*

Clearly, any two points $x, y \in \Omega \cup A_1 \cup \cdots \cup A_n$ can be joined by an $\Omega$-geodesic. It is also easy to see that $[x, y]_\Omega$ is unique if $x, y$ belong to the same tree-component of $\Omega$. In general $[x, y]_\Omega$ is unique up to the choices of appropriate geodesic segments in the sets $A_i$.

The following proposition follows from the standard properties of approximating trees and thin geodesic polygons in hyperbolic spaces.

**Proposition 5.6.** *There is a constant $D = D(n)$ with the following property. Let $a_1, \ldots, a_n$ be distinct points of $X$. Let $\Omega = \Omega(a_1, \ldots, a_n)$ be a minimal connector for these points. Then the following hold:*
  (1) *For any $i \neq j$ the $\Omega$-geodesic $[a_i, a_j]_\Omega$ is a $(1, D)$-quasigeodesic in $X$ which is $D$-Hausdorff close to any $X$-geodesic $[a_i, a_j]$.*
  (2) *Let $P$ be the union of all $X$-geodesics with endpoints in the set $\{a_1, \ldots, a_n\}$. Then $P$ and $\Omega$ are $D$-Hausdorff close.*

**Lemma 5.7.** *For any integer $n \geq 2$ there exists a constant $K_1 = K_1(n) > D(n) > 1$ with the following property.*

*Suppose $A_1, \ldots, A_n$ are connected closed nonempty subsets of $X$ which are $4$-quasiconvex in $X$. Suppose further that $d(A_k, A_l) \geq K_1(n)$ for $k \neq l$. Let $\Omega = \Omega(A_1, \ldots, A_n)$. Let $T$ be a tree component of $\Omega$ and let $a_i \in A_i$, $a_j \in A_j$ be two distinct terminal vertices of $T$. Let $[a, b] = [A_i, A_j]$ be a bridge.*

*Then $|d(a_i, a_j) - d(a, b)| \leq K_1$, $d(a_i, a) \leq K_1$ and $d(a_j, b) \leq K_1$. Thus $[a_i, a_j]$ is almost the bridge between $A_i$ and $A_j$.*

*Proof.* Let $D = D(n)$ be the constant provided by Proposition 5.6. Put $K_1(n) = 100 + 2n(2D + 210) + 2D$. Suppose that $d(A_k, A_l) \geq K_1(n)$ for $k \neq l$.

Notice that since $T$ has at most $n$ terminal vertices, it has at most $n - 2$ branch points. We may assume (after re-labeling) that in (1) $i = 1$, $j = 2$ and $a_i = a_1 \in A_1$, $a_j = a_2 \in A_2$. It is obvious that $T$ is a minimal connector between those $A_k$ which contain the terminal vertices of $T$. It is also clear that $T$ is a minimal connector for the set of terminal vertices of $T$.

Consider the $X$-geodesic segment $[a_1, a_2]$ and a bridge $[a, b] = [A_1, A_2]$ where $a \in A_1$ and $b \in A_2$. Also let $\alpha = [a_1, a_2]_\Omega$ be the $\Omega$-geodesic from $a_1$ to $a_2$. By Proposition 5.6 the paths $\alpha$ and $[a_1, a_2]$ are $D$-close. By Lemma 5.2 there are point $p, q \in [a_1, a_2]$ such that $d(p, a) \leq 100$, $d(q, b) \leq 100$, the segments $[a_1, p]$ and $[a_1, a]$ are $100$-close, the segments $[a_2, b]$ and $[a_2, q]$ are $100$-close and the segments $[p, q]$ and $[a, b]$ are $100$-close. By Proposition 5.6 there a point $p' \in [a_1, a_2]_\Omega$ such that $d(p, p') \leq D(n)$. Denote the length of the $\Omega$-geodesic $[a_1, p']_\Omega \subset [a_1, a_2]_\Omega$ by $L$. Suppose $L > n(2D + 210)$.

Since $T$ has at most $n - 2$ branching points, these branching points subdivide the segment of $[a_1, a_2]_\Omega$ from $a_1$ to $p'$ in at most $(n-1)$ segments. Hence the length of one of these segments is at least $L/n > 2D + 210$. Denote this segment by $[s, t]_\Omega$. In particular either $s = a_1$ or $s$ is a branch point of $T$ and either $t = p'$ or $t$ is a branch point of $T$. Then there exist point $s', t' \in [a_1, a]$ such that $d(s, s') \leq D + 100$ and $d(t, t') \leq D + 100$. Since $A_1$ is $4$-quasiconvex, there are points $s'', t'' \in A_1$ such that $d(s, s'') \leq D + 104$ and $d(t, t'') \leq D + 104$. Recall that the $\Omega$-geodesic $[s, t]_\Omega$ has no branch points of $T$ in its interior. Let $\Omega'$ be obtained from $\Omega$ by removing the interior of $[s, t]_\Omega$ and adding either the segment $[s, s'']$ or the segment $[t, t'']$ depending on which makes $\Omega'$ a connector of $A_1, \ldots, A_n$. (In the figure this is $[t, t'']$.)

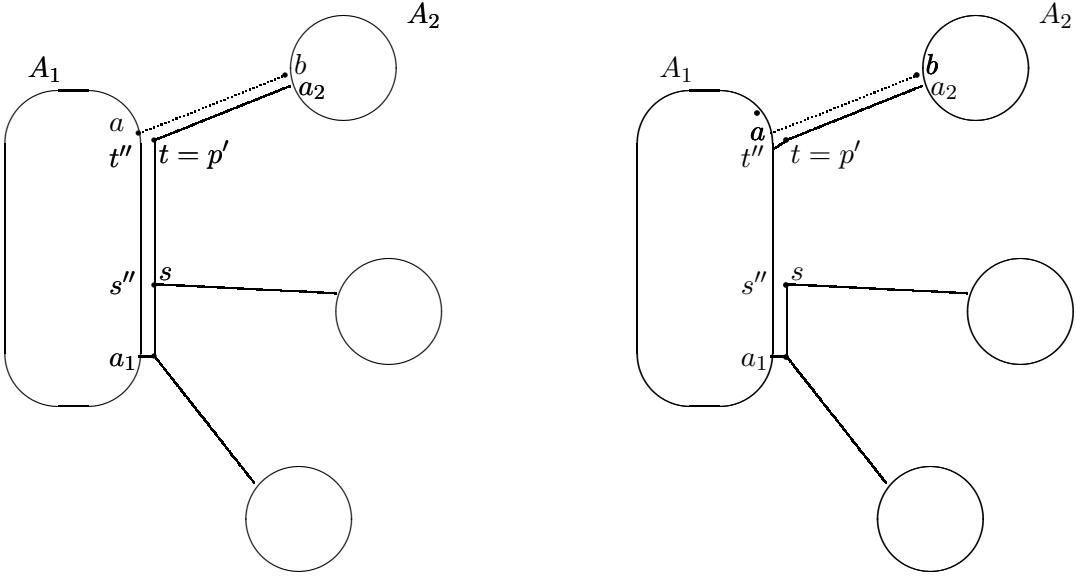

Figure 2. The new connector is shorter, yielding a contradiction to the minimality of $\Omega$.

The perimeter of $\Omega'$ can be estimated as:
$$P(\Omega') \leq P(\Omega) - d_\Omega(s, t) + 2D + 210 < P(\Omega) - 1,$$

which contradicts our assumption that $\Omega$ is a minimal connector of $A_1, \ldots, A_n$. Thus $L \leq n(2D + 210)$ and hence $d(a_1, p) \leq n(2D + 210) + d(p', p) \leq n(2D + 210) + D$. A symmetric argument shows that $d(a_2, q) \leq n(2D + 210) + D$. Hence by Lemma 5.2

$$d(a_1, a) \leq n(2D+210)+D+100 \leq K_1, \quad d(a_2, b) \leq n(2D+210)+D+100 \leq K_1$$

and
$$|d(a, b) - d(a_1, b_1)| \leq 100 + 2n(2D + 210) + 2D = K_1$$

as required. □

**Lemma 5.8.** *Let $n > 0$ and $K_1(n) > D(n)$ be the constants provided by Lemma 5.7 and Proposition 5.6.*

*There exists a constant $K_2 = K_2(n) > K_1(n) > 0$ with the following property. Let $A_1, \ldots, A_n$ be 4-quasiconvex closed connected subsets of $X$ such that $d(A_i, A_j) \geq K_1(n)$ for $i \neq j$. Let* $\Omega = \Omega(A_1, \ldots, A_n)$ be a minimal connector. Let $p, q \in A_i$ be terminal vertices of distinct tree-components $T_p$ and $T_q$ of $\Omega$. Let $p' \ne p$ be another terminal vertex of $T_p$ and let $q' \ne q$ be another terminal vertex of $T_q$.

Then the following hold:

(1) If $x \in [p, p']$ is a point such that $d(x, [q, q']) \le 2$ then $d(p, x) \le K_2$.
(2) The path $[p', p]_\Omega \cup [p, q] \cup [q, q']_\Omega$ is a $(K_2, K_2)$-quasigeodesic.

*Proof.* We first argue that part (2) follows from part (1). By considering a 2-thin geodesic quadrilateral with sides $[p', p]$, $[p, q]$, $[q, q']$ and $[q', p']$ part (1) of Lemma 5.8 and Lemma 5.2 imply that the path $[p', p] \cup [p, q] \cup [q, q']$ is a $(K_2', K_2')$-quasigeodesic for some constant $K_2' = K_2'(n)$. Since by Proposition 5.6 the paths $[p', p]_\Omega$ and $[q, q']_\Omega$ are $(1, D)$-quasigeodesic, statement (2) of Lemma 5.8 follows. Thus it suffices to establish part (1) of Lemma 5.8.

Suppose $x \in [p, p']$ is such that for some $y \in [q, q']$ we have $d(x, y) \le 2$. Assume that $d(p, x) > n(4D + 2K_1 + 31) + D$.

By Proposition 5.6 there are points $x' \in [p, p']_\Omega \subseteq T_p$ and $y' \in [q, q']_\Omega \subseteq T_q$ such that $d(x, x') \le D$ and $d(y, y') \le D$. Therefore $d(x', y') \le 2D + 2$ and

$$d_\Omega(p, x') \ge d(p, x') \ge d(p, x) - D > n(4D + 2K_1 + 31).$$

We claim that the segment $[p, q]$ is short. Indeed, by Lemma 5.7 if $p_1$ is a projection of $x'$ on $A_i$ and if $q_1$ is a projection of $y'$ on $A_i$ then $d(p, p_1) \le K_1 = K_1(n)$ and $d(q, q_1) \le K_1 = K_1(n)$. Let $p_2 \in [p_1, q_1]$ be the farthest from $p_1$ point which is contained in the 2-neighborhood of $[p_1, x']$. Similarly, let $q_2 \in [p_1, q_1]$ be the farthest from $q_1$ point which is contained in the 2-neighborhood of $[q_1, y']$. Since $A_i$ is 2-quasiconvex, Lemma 3.11 implies that $d(p_1, p_2) \le 10$ and $d(q_1, q_2) \le 10$. Since the geodesic quadrilateral with vertices $p_1, q_1, x', y'$ is 2-thin, either $d(p_1, q_1) \le 20$ or the segment $[p_2, q_2]$ is contained in the 2-neighborhood of $[x', y']$. Since $d(x', y') \le 2D+2$, this implies that either $d(p_2, q_2) \le 2D + 6$ or $d(p_1, q_1) \le 20$. In any event we see that $d(p_1, q_1) \le 2D + 26$ and hence $d(p, q) \le 2D + 2K_1 + 26$. Since $d(x', y') \le 2D + 2$, the definition of 1-hyperbolicity implies that $[p, x']$ and $[q, x']$ are $2D + 2K_1 + 28$-Hausdorff close. By Lemma 5.6 $[p, x']$ and $[p, x']_\Omega$ are $D + 1$-close. Similarly, $[q, y']$ and $[q, y']_\Omega$ are $D + 1$-close. Therefore $[p, x']_\Omega$ and $[q, y']_\Omega$ are $4D + 2K_1 + 30$-close.

Recall that $T_p$ and $T_q$ are distinct tree-components of $\Omega$ and that $d(p, x')_\Omega > n(4D+2K_1+31)$. The tree $T_q$ has at most $n - 2$ branching points and therefore the interior of $[p, x']_\Omega$ has at most $n - 2$ branching points. Let $\Omega'$ be obtained from $\Omega$ by removing the interior of $[p, x']_\Omega$ and for each $\Omega$-branching point $b$ in the interior of $[p, x']_\Omega$ adding to the result a geodesic path $[b, b']$ from $b$ to $b' \in [q, y']$ with $d(b, b') \le 4D + 2K_1 + 30$. Then $\Omega' \cup A_1 \cup \ldots A_n$ is clearly connected, but

$$P(\Omega') \le P(\Omega) - d_\Omega(p, x') + n(4D + 2K_1 + 30) < P(\Omega) - 1$$

which contradicts the assumption that $\Omega$ is a minimal connector. Thus $d(p, x) \le n(4D + 2K_1 + 31) + D$ and part (1) of the lemma holds with $K_2 = n(4D + 2K_1 + 31) + D$. □

**Lemma 5.9.** *There exists a constant $K_3 = K_3(n) > K_2(n) > 0$ with the following property.*

*Suppose $n \ge 2$ and $A_1, \ldots, A_n$ are connected closed nonempty 4-quasiconvex sets in $X$ such that $d(A_i, A_j) \ge K_3(n)$ for $i \ne j$. Let $\Omega = \Omega(A_1, \ldots, A_n)$ be a minimal connector. Then for any $a, b \in \Omega \cup A_1 \cup \cdots \cup A_n$ any $\Omega$-geodesic $[a, b]_\Omega$ is a $(K_3, K_3)$-quasigeodesic.*

*Proof.* Recall that $\Omega$-geodesics have the explicit form described in Definition 5.5. The statement of Lemma 5.9 now follows directly from Lemma 5.8 by the "pasting of quasigeodesics" argument, that is by Lemma 3.7. $\square$

## 6. Proof of the main result in the purely elliptic case

In this section we will consider $G$-tuples of the type $M = (U_1, \ldots, U_n; (-))$, where the hyperbolic set $H$ is empty. Our goal is to establish Proposition 2.4 for $G$-tuples of this type. Recall that by the definition of a $G$-tuple all $U_i$ are non-trivial and therefore $X(U_i) \neq \emptyset$. Till the end of this section we assume that $(X, d)$ is a 1-hyperbolic strongly geodesic $G$-space.

**Definition 6.1.** *For a $G$-tuple $M = (U_1, \ldots, U_n; (-))$ with $n \geq 1$ we define*
$$d(M) := d(X(U_1), \ldots, X(U_n)).$$

**Definition 6.2** (Minimal tuple). *Let $M = (U_1, \ldots, U_n; (-))$ be a $G$-tuple with $n \geq 1$.*

*We say that $M$ is* minimal *if for every $G$-tuple $M'$ that is equivalent to $M$ we have $d(M) \leq d(M') + 1$.*

Once again, it is clear that every $G$-tuple with $n \geq 1$ and empty hyperbolic component is equivalent to a minimal tuple. The main result of this section is:

**Proposition 6.3.** *For any integer $n \geq 2$ there exists a constant $C = C(n) > 0$ with the following property. Suppose $(X, d)$ is a strongly geodesic 1-hyperbolic $G$-space. Suppose $(U_1, \ldots, U_n; (-))$ is a minimal $G$-tuple with $U_i \neq 1$ for $1 \leq i \leq n$ and $d(X_{U_i}, X_{U_j}) \geq C$ for $1 \leq i < j \leq n$. Then*

(1) *The subgroup $U$ generated by $U_1 \cup \cdots \cup U_n$ is the free product*
$$U = U_1 * \cdots * U_n$$

(2) *Denote $X_i = X(U_i)$ for $i = 1, \ldots, n$. Let $\Omega = \Omega(X_1, \ldots, X_n)$. The set $X(U)$ is $C$-Hausdorff close to the set*
$$\bigcup_{u \in U} u \cdot (\Omega \cup X_1 \cup \cdots \cup X_n)$$

**Convention 6.4.** *Till the end of this section, unless specified otherwise, let $U_1, \ldots, U_n \leq G$ be non-trivial subgroups of $G$ such that the $G$-tuple $M = (U_1, \ldots, U_n, (-))$ is minimal. Denote $X_i = X(U_i)$ for $i = 1, \ldots, n$. We assume that $d(X_i, X_j) \geq K_3(n)$, where $K = K_3(n)$ is the constant provided by Lemma 5.9. Recall that the sets $X_i$ are 4-quasiconvex, closed and connected since $X$ is 1-hyperbolic. Recall also that each $X_i$ is $U_i$-invariant. Let $\Omega = \Omega(X_1, \ldots, X_n)$ be a minimal connector.*

The following lemma will play a crucial role in our argument.

**Lemma 6.5.** *There exists a constant $K_4(n) > K_3(n) > 0$ with the following property.*

*Suppose that $T$ and $T'$ are tree-components of $\Omega$, and that $p \in X_i$ is a terminal vertex of $T$ and $p' \in X_i$ is a terminal vertex of $T'$. Let $q \in T$ and $q' \in T'$ be terminal vertices of $T$ and $T'$, respectively, such that $q, q' \notin X_i$. Let $u \in U_i$ be a non-trivial element.*

*Then the path*
$$\sigma = [q, p]_\Omega \cup [p, up'] \cup [up', uq']$$
*is a $(K_4, K_4)$-quasigeodesic.*

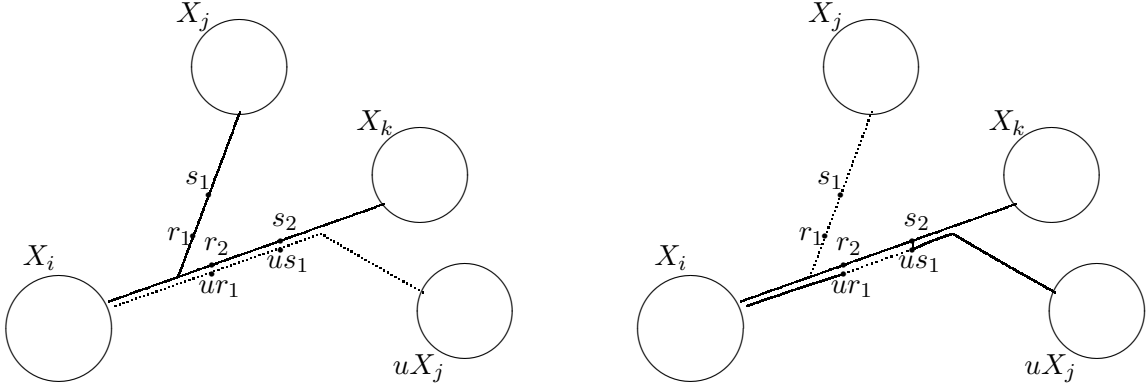

FIGURE 3. The case $T = T'$. Left side: $\Omega = \Omega(X_i, X_j, X_k)$.
Right side: $\Omega' = \Omega(X_i, X_j, uX_j)$. Clearly $P(\Omega') < P(\Omega)$.

*Proof.* Recall that $T'$ is a tree-component of $\Omega$ with terminal vertices $p' \in X_i$ and $q' \in X_k$ (where $k \neq i, k \neq j$). Note that it is possible that $p = p'$ or even that $T = T'$. Clearly it is enough to prove the following:

**Claim.** There is a constant $K'' > 1$ such that if $t \in [p', q']$ and $s \in [p, q]$ satisfy $d(t, us) \leq 2$, then $d(p', t) \leq K''$ and $d(p, s) \leq K''$.

Suppose now that $t \in [p', q']$ and $s \in [p, q]$ are such that $d(t, us) \leq 2$. Put $L = 24D + 1000$ and suppose that

$$d(p, s) > K'(nL + 2D + 2) + K'.$$

Recall by Proposition 5.6 the point $p'$ is $D$-close to a projection $p_1$ of $t$ on $X_i$. Similarly, $p$ is $D$-close to a projection $p_2$ of $s$ on $X_i$. Thus $[s, p]$ and $[t, p']$ are almost "orthogonal" to $X_i$. Namely, by Lemma 3.12 either $d(p', up) \leq 2D + 100$ or the path $\tau' = [t, p'] \cup [p', up] \cup u[p, s]$ is $(1, 20+4D)$-quasigeodesic. Since $d(t, us) \leq 2$, the former possibility implies that $[p', t]$ and $u[p, s]$ are $2D+100$-Hausdorff close. Suppose now the latter occurs and $\tau'$ is $(1, 20+4D)$-quasigeodesic. Since the endpoints $t$ and $us$ of $\tau'$ are at most 2-apart, this implies that $d(p', t) \leq 22 + 4D$ and $d(up, us) = d(p, s) \leq 22 + 4D$, contrary to our assumption.

Thus $d(p', up) \leq 2D+100$ and the paths $[p', t]$ and $u[p, s]$ are $2D+100$-Hausdorff close. Recall that by assumption $d(p', t) > nL$.

By Proposition 5.6 there are points $t' \in T'$ and $s' \in T$ such that $d(t, t') \leq D$ and $d(s, s') \leq D$ and hence $d(t', us') \leq 2D+2$. If $p = p'$, let $r \in T$ be such that $[p, r]_\Omega = [p, s']_\Omega \cap [p, t']_\Omega$. If $p \neq p'$ let $r = p$. By part (1) of the lemma the path $[r, p]_\Omega \cup [p, up] \cup u[p, s]$ is a $(K', K')$-quasigeodesic.

Suppose first $p = p'$ and $d(r, t') \leq nL$. Then $d(r, us) \leq d(r, t') + d(t', us') \leq nL + 2D + 2$. Therefore the length of the path $[r, p]_\Omega \cup [p, up] \cup u[p, s]$ is at most $K'(nL+2D+2)+K'$. Hence $d(p, s) \leq K'(nL+2D+2)+K'$ and $d(p, s) \leq K'(nL+2D+2)+K'$, contrary to our assumptions.

Thus either $p = p'$ and $d(r, t') \geq nL$ and or $p \neq p'$. Assume first that $p = p'$ and $d(r, t') \geq nL$. Then $d_\Omega(r, t) \geq nL$. Since the paths $[p, t]$ and $u[p, s]$ are $2D + 100$-Hausdorff close, we conclude that $|d(p, t) - d(p, s)| \leq 4D + 200$ and hence $|d(p, s') - d(p, t')| \leq 6D + 200$. This implies that $|d(r, s') - d(r, t')| \leq 12D + 400$ and so $d_\Omega(r, s') \geq d(r, s') \geq nL - 12D - 400$. If $p \neq p'$ then $r = p$. In this case $d(p, s) > nL$ implies $d(r, s') \geq nL - D$. Thus in both cases $d(r, s') \geq nL - 12D - 400$.

The tree $T$ has at most $n - 2$ branch-points. Therefore there is a segment $[r_1, s_1]_\Omega \subseteq [r, s']_\Omega$ such that the interior of $[r_1, s_1]_\Omega$ has no branch-points and such that

$$d_\Omega(r_1, s_1) \geq (nL - 12D - 400)/n \geq L/2$$

(where the last inequality follows from the choice of $L$). Since $[p, t']_\Omega$ and $u[p, s]_\Omega$ are $4D + 100$-Hausdorff close, there are points $r_2, s_2 \in [p, t]_\Omega$ such that $d(ur_1, r_2) \leq 4D + 100$ and $d(us_1, s_2) \leq 4D + 100$.

Let $J$ be the collection of all indices $j'$ such that an $\Omega$-geodesic from $X_{j'}$ to $X_i$ passes through the interior of $T$ (so that $j \in J$). Let $\Omega_1$ be the union of all tree-components of $\Omega$ whose interiors intersect non-trivially $\Omega$-geodesics from $X_{j'}$ to $X_i$ for $j' \in J$. Thus $T$ is contained in $\Omega_1$ and $\Omega_1$ is a minimal connector for the collection $X_i \cup \bigcup_{j' \in J} X_{j'}$.

Let $\Omega'$ be obtained from $\Omega$ as follows:

Remove $\Omega_1$ from $\Omega$ and replace it by $u\Omega_1$. Then remove the segment $u[r_1, s_1]_\Omega$ and add geodesic segments $[ur_1, r_2]$ and $[us_1, s_2]$ (which by the previous remarks have length at most $4D + 100$ each). Then by construction

$$P(\Omega') \leq P(\Omega) + 8D + 200 - d(r_1, s_1) \leq P(\Omega) + 8D + 200 - L/2 < P(\Omega) - 1.$$

Put $V_k = U_k$ for $k \notin J$. Also put $V_k = uU_k u^{-1}$ for $k \in J$ (so that $X(V_k) = uX_k$). Then $(V_1, \ldots, V_n; (-))$ is equivalent to $(U_1, \ldots, U_n; (-))$ and $\Omega'$ is such that

$$\Omega' \cup X(V_1) \cup \cdots \cup X(V_n)$$

is connected. The fact that $P(\Omega') < P(\Omega) - 1$ contradicts the assumption that $(U_1, \ldots, U_n; (-))$ is a minimal $G$-tuple.

Thus we in fact have $d(p, s) \leq K'(nL + 2D + 2) + K'$. We have shown in the argument above that either $d(p', up) \leq 2D + 100$ or the path $\tau' = [t, p'] \cup [p', up] \cup u[p, s]$ is $(1, 20 + 4D)$-quasigeodesic. If $d(p', up) \leq 2D + 100$ then by the triangle inequality we have

$$d(p', t) \leq d(p', up) + d(up, us) + d(us, t) \leq K'(nL + 2D + 2) + K' + 2D + 102.$$

If $\tau'$ is a $(1, 20 + 4D)$-quasigeodesic then $d(t, us) \leq 2$ implies that the length of $\tau$ is at most $22 + 4D$. Hence $d(p', t) = \leq 22 + 4D$. Thus in any event we see that

$$d(p, s) \leq K'', \quad d(p', t) \leq K''$$

where $K'' = K'(nL + 2D + 2) + K' + 4D + 102$. Thus the Claim is verified and Lemma 6.5 is proved. □

**Corollary 6.6.** *There exists a constant $K_5(n) > K_4(n) > 0$ with the following property.*

*Let $p', q' \in (\Omega \cup X_1 \cup \ldots \cup X_n) - X_i$ and $p, q \in X_i$ such that $d_\Omega(X_i, p') = d_\Omega(p, p')$ and $d_\Omega(X_i, q') = d_\Omega(q, q')$. Then the path*

$$\sigma = [q, p]_\Omega \cup [p, up'] \cup [up', uq']$$

*is a $(K_5, K_5)$-quasigeodesic.*

*Proof.* This follows directly from Lemma 6.5 and Lemma 5.9.  □

*Proof of Proposition 6.3.* Let $n \geq 2$ and $K_5(n) > K_4(n) > K_3(n) > K_2(n) > K_1(n) > D(n)$ be constants provided by Corollary 6.6, Lemma 6.5, Lemma 5.9, Lemma 5.8, Lemma 5.7 and Proposition 5.6 accordingly.

Let $K_6(n) > K_5(n)$ be such that any $K_6$-local $(K_5, K_5)$-quasigeodesic is a global $(K_6, K_6)$-quasigeodesic. Let $L > 100$ be such that any two $(K_6, K_6)$-quasigeodesics in $X$ with common endpoints are $L$-Hausdorff close.

Suppose now that $M = (U_1, \ldots, U_n; (-))$ is a minimal $G$-tuple such that $d(X(U_i), X(U_j)) > K_6^2 + 10K_6 + K_6(K_6 + 10L + 20)$ for $i \neq j$. Denote $X_i = X(U_i)$ for $i = 1, \ldots, n$ and let $\Omega = \Omega(X_1, \ldots, X_n)$ be a minimal connector.

(a) We will first establish that $U$ is indeed a free product $U = U_1 * U_2 \ldots U_n$. Choose a point $z \in \Omega \cup X_1 \cup \cdots \cup X_n$.

Let $u = u_1 u_2 \ldots u_k$ be a strictly alternating product, where $k \geq 1$, $u_j \in U_{i_j}, u_j \neq 1$ and where $i_j \neq i_{j+1}$.

For each $j = 1, \ldots, k-1$ let $p_j \in X_{i_j}$ be the $d_\Omega$-closest to $X_{i_{j+1}}$ point of $X_{i_j}$. Similarly, for each $j = 2, \ldots, k$ let $q_j \in X_{i_j}$ be the $d_\Omega$-closest to $X_{i_{j-1}}$ point of $X_{i_j}$. Thus $p_j, q_j \in X_{i_j}$ are terminal vertices of $\Omega$. Since $u_j \in U_{i_j}$ we have $u_j X_{i_j} = X_{i_j}$ and $u_j p_j, u_j q_j \in X_{i_j}$. Note that since $i_j \neq i_{j+1}$ we have $d(X_{i_j}, X_{i_{j+1}}) > K_6^2 + 10K_6$. Hence $d_\Omega(p_j, q_{j+1}) \geq d(p_j, q_{j+1}) > K_6^2 + 10K_6$.

Consider the path

$$\sigma = [z, u_1 p_1]_\Omega \cup u_1[p_1, q_2]_\Omega \cup u_1[q_2, u_2 p_2] \cup u_1 u_2[p_2, q_3]_\Omega \cup u_1 u_2[q_3, u_3 p_3] \cup \ldots$$
$$\ldots \cup u_1 u_2 \ldots u_{k-1}[p_{k-1}, u_k q_k] \cup u_1 u_2 \ldots u_{k-1} u_k[q_k, z]_\Omega$$

It is clear that $\sigma$ is a path from $z$ to $uz$. By Corollary 6.6 each path $[p_{i-1}, q_i]_\Omega \cup [q_i, u_i p_i] \cup u_i[p_i, q_i]_\Omega$ is a $(K_5, K_5)$-quasigeodesic. Similarly, the paths $[z, u_1 p_1]_\Omega \cup u_1[p_1, q_2]_\Omega$ and $[p_{k-1}, u_k q_k] \cup u_k[q_k, z]_\Omega$ are $(K_5, K_5)$-quasigeodesic. Since $d_\Omega(p_j, q_{j+1}) \geq 10K_6$ for $j = 1, \ldots, k-1$, the path $\sigma$ is a $K_6$-local $(K_5, K_5)$-quasigeodesic. Hence by the choice of $K_6$ the path $\sigma$ is a global $(K_6, K_6)$-quasigeodesic. If $k = 1$ then $u = u_1 \neq 1$ by the choice of $u$. If $k \geq 2$ then the length of $\sigma$ is at least $K_6^2 + 10K_6$ and hence

$$d(z, uz) \geq (K_6^2 + 10K_6)/K_6 - K_6 > 0$$

Hence $u \neq 1$. Thus we have established that $U$ is indeed a free product $U = U_1 * U_1 * \ldots U_n$.

(b) Let

$$A = \bigcup_{u \in U} u(\Omega \cup X_1 \cup \cdots \cup X_n).$$

We can now show that $A$ is quasiconvex in $X$. Recall that by Lemma 5.9 for any $a_1, a_2 \in A$ any $\Omega$-geodesic $[a_1, a_2]_\Omega$ is $(K_3, K_3)$-quasigeodesic. Since $K_6 > K_3$, the path $[a_1, a_2]_\Omega$ is also $(K_6, K_6)$-quasigeodesic. Since $A$ is obviously $U$-invariant, to see that $A$ is $C$-quasiconvex it suffices to show that for any $z, z' \in A$ and $u \in U$ a geodesic $[z, uz']$ is contained in a $C$-neighborhood of $A$. Assume $u \neq 1$. Consider the 1-thin geodesic triangle $[z, uz] \cup u[z, z'] \cup [z, uz']$. Let $\sigma$ be a path from $z$ to $uz$ as in (1). Since $\sigma$ is a $(K_6, K_6)$-quasigeodesic, by the choice of $L$ the paths $\sigma$ and $[z, uz]$ are $L$-close. Since $X_i$ are 4-quasiconvex, the construction of $\sigma$ implies that $\sigma$ is contained in the 4-neighborhood of $A$. Thus $[z, uz]$ is contained in $(L+4)$-neighborhood of $A$. Since $[z, z']_\Omega$ is a $(K_6, K_6)$-quasigeodesic, the path $u[z, z']$ is contained in $(L+4)$-neighborhood of $A$. Hence by the definition of hyperbolicity $[z, uz']$ is contained in $(L+5)$-neighborhood of $A$. This is also obviously true for $u = 1$. Thus we have proved that $A$ is $(L+5)$-quasiconvex.

(c) We will now observe that $A$ is contained in a bounded neighborhood of $X(U)$.

Since $U_i \leq U$, we have $X_i = X(U_i) \subseteq X(U)$. Therefore $\bigcup_{u \in U} u(X_1 \cup \cdots \cup X_n)$ is contained in $X(U)$. For any two terminal vertices $p, q$ of a tree-component $T$ of $\Omega$ the path $[p, q]_\Omega$ is $D$-close to a geodesic $[p, q]$. Since $X(U)$ is 4-quasiconvex, and $p, q \in X(U)$ this implies that $[p, q]_\Omega$ is contained in the $(D+4)$-neighborhood of $X(U)$. Hence $\Omega$ is contained in the $(D+4)$-neighborhood of $X(U)$ and therefore so is $u\Omega$ for any $u \in U$. Thus $A$ is contained in the $(D+4)$-neighborhood of $X(U)$, as required.

(d) We will now observe that points outside of $A$ are moved by elements of $U$ by a substantial distance. Namely, we will show that $E(U) \subseteq A$. Since $E(U_i) \subseteq X(U_i) = X_i \subseteq A$, we know that for any $x \in X - A$ and any $u \in U_i, u \neq 1$ we have $d(x, ux) > 100$. Hence for any $u' \in U$ and any $y \in X - A$ we also have $d(y, (u')^{-1}uu'y) > 100$.

Suppose now $u \in U, u \notin U_i$, $x \in X - A$ and $d(x, ux) \leq 100$. Let $z$ be a projection of $x$ on $A$. By conjugation of $u$ in $U$ we may assume that $z \in \Omega \cup X_1 \cup \cdots \cup X_n$. We may further assume that $u$ has syllable length at least 2, that is $u = u_1 \ldots u_k$, where $k \geq 2$ where $k \geq 1$, $u_j \in U_{i_j}, u_j \neq 1$ and where $i_j \neq i_{j+1}$. Let $\sigma$ be a path from $z$ to $uz$ constructed as in the proof of (1). Then the length of $\sigma$ is at least $d_\Omega(p_1, q_2)$. By assumption

$$d_\Omega(p_1, q_2) \geq d(p_1, q_2) \geq d(X_{i_1}, X_{i_2}) \geq K_6(K_6 + 10L + 20).$$

Since $\sigma$ is a $(K_6, K_6)$-quasigeodesic, this implies that

$$d(z, uz) \geq K_6(K_6 + 10L + 20)/K_6 - K_6 = 10L + 20.$$

Since $A$ is $L$-quasiconvex, $z$ is a projection of $x$ on $A$ and $uz$ is a projection of $ux$ on $A$, Lemma 3.11 implies that $d(x, ux) \geq 2L$. Since $L$ was chosen $L > 100$, this implies $d(x, ux) > 200$ and hence $x \notin E(U)$. Thus we have established that $E(U) \subseteq A$.

We now know that $A$ is $U$-invariant, $(L+5)$-quasiconvex, contains $E(U)$ and is contained in $(D+4)$-neighborhood of $X(U)$. Therefore by Lemma 4.5 the sets $A$ and $X(U)$ are $(3L+3D+26)$-Hausdorff close.

Thus we have verified that Proposition 6.3 holds with

$$C(n) = K_6^2 + 10K_6 + K_6(K_6 + 10L + 20) + (3L + 3D + 26)$$

$\square$

*Proof of Proposition 6.3.* Let $n \geq 2$ and $K_4(n) > K_3(n) > K_2(n) > K_1(n) > D(n)$ be constants provided by Lemma 6.6, Lemma 5.9, Lemma 5.8, Lemma 5.7 and Proposition 5.6 accordingly.

Let $K_5(n) > K_4(n)$ be such that any $K_5$-local $(K_4, K_4)$-quasigeodesic is a global $(K_5, K_5)$-quasigeodesic. Let $L > 100$ be such that any two $(K_5, K_5)$-quasigeodesics in $X$ with common endpoints are $L$-Hausdorff close.

Suppose now that $M = (U_1, \ldots, U_n; (-))$ is a minimal $G$-tuple such that $d(X(U_i), X(U_j)) > K_5^2 + 10K_5 + K_5(K_5 + 10L + 20)$ for $i \neq j$. Denote $X_i = X(U_i)$ for $i = 1, \ldots, n$ and let $\Omega = \Omega(X_1, \ldots, X_n)$ be a minimal connector.

(a) We will first establish that $U$ is indeed a free product $U = U_1 * U_2 \ldots U_n$. Choose a point $z \in \Omega \cup X_1 \cup \cdots \cup X_n$.

Let $u = u_1 u_2 \ldots u_k$ be a strictly alternating product, where $k \geq 1$, $u_j \in U_{i_j}, u_j \neq 1$ and where $i_j \neq i_{j+1}$.

For each $j = 1, \ldots, k-1$ let $p_j \in X_{i_j}$ be the $d_\Omega$-closest to $X_{i_{j+1}}$ point of $X_{i_j}$. Similarly, for each $j = 2, \ldots, k$ let $q_j \in X_{i_j}$ be the $d_\Omega$-closest to $X_{i_{j-1}}$ point of $X_{i_j}$. Thus $p_j, q_j \in X_{i_j}$ are terminal vertices of $\Omega$. Since $u_j \in U_{i_j}$ we have $u_j X_{i_j} = X_{i_j}$ and $u_j p_j, u_j q_j \in X_{i_j}$. Note that since $i_j \neq i_{j+1}$ we have $d(X_{i_j}, X_{i_{j+1}}) > K_5^2 + 10K_5$. Hence $d_\Omega(p_j, q_{j+1}) \geq d(p_j, q_{j+1}) > K_5^2 + 10K_5$.

Consider the path

$$\sigma = [z, u_1 p_1]_\Omega \cup u_1[p_1, q_2]_\Omega \cup u_1[q_2, u_2 p_2] \cup u_1 u_2[p_2, q_3]_\Omega \cup u_1 u_2[q_3, u_3 p_3] \cup \ldots$$
$$\ldots \cup u_1 u_2 \ldots u_{k-1}[p_{k-1}, u_k q_k] \cup u_1 u_2 \ldots u_{k-1} u_k[q_k, z]_\Omega$$

It is clear that $\sigma$ is a path from $z$ to $uz$. By Lemma 6.6 each path $[p_{i-1}, q_i]_\Omega \cup [q_i, u_i p_i] \cup u_i[p_i, q_i]_\Omega$ is a $(K_4, K_4)$-quasigeodesic. Similarly, the paths $[z, u_1 p_1]_\Omega \cup u_1[p_1, q_2]_\Omega$ and $[p_{k-1}, u_k q_k] \cup u_k[q_k, z]_\Omega$ are $(K_4, K_4)$-quasigeodesic. Since $d_\Omega(p_j, q_{j+1}) \geq 10K_5$ for $j = 1, \ldots, k-1$, the path $\sigma$ is a $K_5$-local $(K_4, K_4)$-quasigeodesic. Hence by the choice of $K_5$ the path $\sigma$ is a global $(K_5, K_5)$-quasigeodesic. If $k = 1$ then $u = u_1 \neq 1$ by the choice of $u$. If $k \geq 2$ then the length of $\sigma$ is at least $K_5^2 + 10K_5$ and hence

$$d(z, uz) \geq (K_5^2 + 10K_5)/K_5 - K_5 > 0$$

Hence $u \neq 1$. Thus we have established that $U$ is indeed a free product $U = U_1 * U_1 * \ldots U_n$.

(b) Let

$$A = \bigcup_{u \in U} u(\Omega \cup X_1 \cup \cdots \cup X_n).$$

We can now show that $A$ is quasiconvex in $X$. Recall that by Lemma 5.9 for any $a_1, a_2 \in A$ any $\Omega$-geodesic $[a_1, a_2]_\Omega$ is $(K_3, K_3)$-quasigeodesic. Since $K_5 > K_3$, the path $[a_1, a_2]_\Omega$ is also $(K_5, K_5)$-quasigeodesic. Since $A$ is obviously $U$-invariant, to see that $A$ is $C$-quasiconvex it suffices to show that for any $z, z' \in A$ and $u \in U$ a geodesic $[z, uz']$ is contained in a $C$-neighborhood of $A$. Assume $u \neq 1$. Consider the 1-thin geodesic triangle $[z, uz] \cup u[z, z'] \cup [z, uz']$. Let $\sigma$ be a path from $z$ to $uz$ as in (1). Since $\sigma$ is a $(K_5, K_5)$-quasigeodesic, by the choice of $L$ the paths $\sigma$ and $[z, uz]$ are $L$-close. Since $X_i$ are 4-quasiconvex, the construction of $\sigma$ implies that $\sigma$ is contained in the 4-neighborhood of $A$. Thus $[z, uz]$ is contained in $(L+4)$-neighborhood of $A$. Since $[z, z']_\Omega$ is a $(K_5, K_5)$-quasigeodesic, the path $u[z, z']$ is contained in $(L+4)$-neighborhood of $A$. Hence by the definition of hyperbolicity $[z, uz']$ is contained in $(L+5)$-neighborhood of $A$. This is also obviously true for $u = 1$. Thus we have proved that $A$ is $(L+5)$-quasiconvex.

(c) We will now observe that $A$ is contained in a bounded neighborhood of $X(U)$.

Since $U_i \leq U$, we have $X_i = X(U_i) \subseteq X(U)$. Therefore $\bigcup_{u \in U} u(X_1 \cup \cdots \cup X_n)$ is contained in $X(U)$. For any two terminal vertices $p, q$ of a tree-component $T$ of $\Omega$ the path $[p, q]_\Omega$ is $D$-close to a geodesic $[p, q]$. Since $X(U)$ is 4-quasiconvex, and $p, q \in X(U)$ this implies that $[p, q]_\Omega$ is contained in the $(D+4)$-neighborhood of $X(U)$. Hence $\Omega$ is contained in the $(D+4)$-neighborhood of $X(U)$ and therefore so is $u\Omega$ for any $u \in U$. Thus $A$ is contained in the $(D+4)$-neighborhood of $X(U)$, as required.

(d) We will now observe that points outside of $A$ are moved by elements of $U$ by a substantial distance. Namely, we will show that $E(U) \subseteq A$. Since $E(U_i) \subseteq X(U_i) = X_i \subseteq A$, we know that for any $x \in X - A$ and any $u \in U_i, u \neq 1$ we have $d(x, ux) > 100$. Hence for any $u' \in U$ and any $y \in X - A$ we also have $d(y, (u')^{-1} u u' y) > 100$.

Suppose now $u \in U, u \notin U_i$, $x \in X - A$ and $d(x, ux) \leq 100$. Let $z$ be a projection of $x$ on $A$. By conjugation of $u$ in $U$ we may assume that $z \in \Omega \cup X_1 \cup \cdots \cup X_n$. We may further assume that $u$ has syllable length at least 2, that is $u = u_1 \ldots u_k$, where $k \geq 2$ where $k \geq 1$, $u_j \in U_{i_j}, u_j \neq 1$ and where $i_j \neq i_{j+1}$. Let $\sigma$ be a path from $z$ to $uz$ constructed as in the proof of (1). Then the length of $\sigma$ is at least $d_\Omega(p_1, q_2)$. By assumption

$$d_\Omega(p_1, q_2) \geq d(p_1, q_2) \geq d(X_{i_1}, X_{i_2}) \geq K_5(K_5 + 10L + 20).$$

Since $\sigma$ is a $(K_5, K_5)$-quasigeodesic, this implies that

$$d(z, uz) \geq K_5(K_5 + 10L + 20)/K_5 - K_5 = 10L + 20.$$

Since $A$ is $L$-quasiconvex, $z$ is a projection of $x$ on $A$ and $uz$ is a projection of $ux$ on $A$, Lemma 3.11 implies that $d(x, ux) \geq 2L$. Since $L$ was chosen $L > 100$, this implies $d(x, ux) > 200$ and hence $x \notin E(U)$. Thus we have established that $E(U) \subseteq A$.

We now know that $A$ is $U$-invariant, $(L+5)$-quasiconvex, contains $E(U)$ and is contained in $(D+4)$-neighborhood of $X(U)$. Therefore by Lemma 4.5 the sets $A$ and $X(U)$ are $(3L+3D+26)$-Hausdorff close.

Thus we have verified that Proposition 6.3 holds with

$$C(n) = K_5^2 + 10K_5 + K_5(K_5 + 10L + 20) + (3L + 3D + 26)$$

□

## 7. The case of one elliptic subgroup

As before, in this and the next sections we assume that $(X, d)$ is a 1-hyperbolic strongly geodesic $G$-space.

Suppose that $g \in G$ and that $V \leq G$ is a non-trivial subgroup. We denote by $l_V(g)$ the distance

$$l_V(g) := d(X(V), gX(V)).$$

We will refer to $l_V(g)$ as $V$-*length of* $g$, or the *length of $g$ relative $V$*. It follows from the definition of $l_V$ that $l_V(g) = l_V(v_1 g v_2)$ for any $g \in G$ and $v_1, v_2 \in V$ since $d(X(V), v_1 g v_2 X(V)) = d(X(V), v_1 g X(V)) = d(v_1^{-1} X(V), v_1^{-1}(v_1 g X(V))) = d(X(V), g X(V))$.

**Definition 7.1** (Minimality). *Let $M = (V; H)$ be a $G$-pair with $H = (h_1, \ldots, h_m)$, where $m \geq 1$ and $V$ is a non-trivial subgroup of $G$. We define the $V$-length $|H|_V$ of $H$ as*

$$|H|_V := l_V(h_1) + \cdots + l_V(h_n).$$

*We will also refer to $|H|_V$ as the* complexity *of the $G$-tuple $M = (V; H)$. We say that a pair $M = (V; H)$ is* minimal *if for any other pair $M = (V'; H')$ equivalent to $M$ we have $|H|_V \leq |H'|_{V'} + 1$.*

Minimal pairs clearly exist in every equivalence class. The main result in this section is the following:

**Proposition 7.2.** *For any integer $m \geq 1$ there exists a constant $R(m) > 0$ with the following property.*

*Suppose $(X, d)$ is a 1-hyperbolic strongly geodesic $G$-space. Let $V \leq G$ be a non-trivial subgroup and $H = (h_1, \ldots, h_m) \in G^m$ (where $m \geq 1$) be such that the $G$-pair $M = (V; H)$ is minimal. Let*

$U = \langle V \cup \{h_1, \ldots, h_m\}\rangle \leq G$. Then either $U = V * F(H)$ or after replacing $M$ by an equivalent pair $M' = (V'; H' = (h'_1, \ldots, h'_m))$ one of the following holds:
  (1) $l_{V'}(h'_i) \leq R(m)$ for some $i \in \{1, \ldots, m\}$.
  (2) There exists a $x \in X$ such that $d(x, h'_i x) \leq R(m)$ for some $i \in \{1, \ldots, m\}$.

**Convention 7.3** (Essential part with respect to a subgroup $V \leq G$). Let $V \subset G$ be a non-trivial subgroup. We denote by $Y(V)$ the 10-neighborhood of $X(V)$ in $X$. Since $X(V)$ is 4-quasiconvex and $V$-invariant and $X$ is 1-hyperbolic, it is clear that $Y(V)$ is also $V$-invariant and 12-quasiconvex.

Choose a base-point $x \in Z(V) \subset X(V) \subset Y(V)$. Recall that $X(V)$ contains all geodesic segments joining points of $Z(V)$.

Let $g \in G$ be such that $l_V(g) > 40$. We choose a geodesic segment $W_g = [x, gx]$ and look at the set $C_g = W_g - (Y(V) \cup gY(V))$. It is clear that $C_g$ is non-empty and that some component of $C_g$ is of length at least $l_V(g) - 20$. It is also easy to see that at most one component of $C_g$ is of length greater than 20. Thus there is a unique component of $C_g$ of length greater than 20. We call this component the *essential part* of $W_g$ relative $x$ and denote it by $E_g^x$.

It is easy to see that the definition of $E_g^x$ does not depend very much on the choice of $x$ and that the length of the stable part is close to the $V$-length of $g$:

**Lemma 7.4.** *Let $V \leq G$ be a non-trivial subgroup and let $g \in G$ be such that $l_V(g) > 40$. Let $[x_1, x_2] = [X(V), gX(V)]$ be a bridge.*
  (1) *For any $x, y \in Z(V)$ the stable parts $E_g^x$ and $E_g^y$ are 10-Hausdorff close.*
  (2) *For any $x \in Z(V)$ we have $l_V(g) - 20 \leq l(E_g) \leq l_V(g) + 20$.*

*Proof.* Choose any $y_1, y_2, y_3, y_4 \in Z(V)$ and any two geodesic segments $W_1 = [y_1, gy_2]$ and $W_2 = [y_3, gy_4]$. Define $E_i$ to be the unique component of length greater than 20 of $W_i - (Y(V) \cup gY(V))$ for $i = 1, 2$. Choose a bridge $[x_1, x_2] = [X(V), gX(V)]$. Since geodesic quadrilaterals are 2-thin, it follows that $E_1$ and $E_2$ lie in the 2 neighborhood of $[x_1, x_2]$. Moreover, it is easy to see that both $E_1$ and $E_2$ are 5-Hausdorff close to $[x_1, x_2] - (X(V) \cup gX(V))$ since the bridge minimizes distance to $X(V)$ and $gX(V)$ (up to the additive constant 1). Therefore $E_1$ and $E_2$ are 10-Hausdorff close. These observations easily imply the statement of the lemma. □

Because of Lemma 7.4 we will often omit the reference to the base-point $x \in Z(V)$ and denote the stable part $E_g^x$ by $E_g$.

The following lemma gives us some information about essential parts of the products of type $g_1 v g_2$ where $g_1, g_2 \in G$ and $v \in V$.

**Lemma 7.5.** *For any $g_1, g_2 \in G$ with $l_V(g_1) \geq 100$ and $l_V(g_2) \geq 100$ the following hold:*
  (1) *There exists at most one $v \in V$ such that $l_V(g_1 v g_2) \leq l_V(g_1) + l_V(g_2) - 200$.*
  (2) *If $l_V(g_1 v g_2) > l_V(g_1) + l_V(g_2) - 200$ then the initial segment of $E_{g_1}$ of length $l_V(g_1) - 100$ and the terminal segment of $g_i v E_{g_2}$ of length $l_V(g_2) - 100$ lie in the 10-neighborhood of $E_{g_1 v g_2}$*
  (3) *If $l_V(g_1 v g_2) \leq l_V(g_1) + l_V(g_2) - 200$ then the initial segment of $E_{g_1}$ of length $l_V(g_1) - \frac{1}{2}(l_V(g_1) + l_V(g_2) - l_V(g_1 v g_2))$ and the terminal*

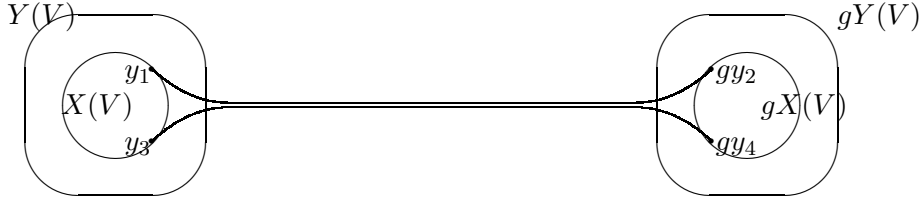

FIGURE 4. The essential part of an element $g$ is almost independent of the choice of $x \in Z(V)$.

*segment of $g_i v E_{g_2}$ of length $l_V(g_2) - \frac{1}{2}(l_V(g_1) + l_V(g_2) - l_V(g_1 v g_2))$ lie in the 10-neighborhood of $E_{g_1 v g_2}$. Furthermore the terminal segment of $E_{g_1}$ of length 100 and the initial segment of $g_1 v E_{g_2}$ of length 100 are 10-Hausdorff close.*

*Proof.* It suffices to prove (1), as (2) and (3) then follow easily from the definitions of hyperbolicity and of the essential part.

Let $v \in V$ be such that $l_V(g_1 v g_2)$ is almost minimal, that is to say for any $v' \in V$ we have $l_V(g_1 v g_2) \le l_V(g_1 v' g_2) + 1$. After replacing $g_1$ with $g_1 v$ we can assume that in fact $v = 1$.

If $l_V(g_1 g_2) > l_V(g_1) + l_V(g_2) - 200$ then there is nothing to prove.

Suppose that $l_V(g_1 g_2) \le l_V(g_1) + l_V(g_2) - 200$. Since geodesic triangles in $X$ are 1-thin, it follows that $W_{g_1 g_2}$ lies in 1-neighborhood of the union of the initial segment of $W_{g_1}$ of length $l(W_{g_1}) - k$ and the terminal segment of $g_1 W_{g_2}$ of length $l(W_{g-2}) - k$ for some positive integer $k$ and the terminal segment of $W_{g_1}$ and the initial segment of $g_1 W_{g_2}$ of length $k$ are 1-Hausdorff close.

It follows as in the proof of Lemma 7.4 that the terminal segment $[x_1, x_2]$ of $E_{g_1}$ of length 100 and the initial segment $[y_2, x_2]$ of $g_1 E_{g_2}$ of length 100 are 10-Hausdorff close. Thus we have a picture as in Figure 5.

Now suppose that there exists an element $u \in V - 1$ such that $l_V(g_1 u g_2) \le l_V(g_1) + l_V(g_2) - 200$. The same argument as above shows that the terminal segment $[x_1, y_1]$ of $E_{g_1}$ of length 100 and the initial segment $[y_2', x_2']$ of $g_1 u E_{g_2}$ of length 100 are 10-Hausdorff close. It is clear that $x_2' = (g_1 u g_1^{-1}) x_2$ and that $y_2' = (g_1 u g_1^{-1}) y_2$, that is we have the picture illustrated in Figure 6

The above statements clearly imply that $d(x_2, x_2') \le 40$. However, it follows from Lemma 4.4 that $[x_2, y_2] \cup [y_2, (g_1 u g_1^{-1}) y_2] \cup [(g_1 u g_1^{-1}) y_2, (g_1 u g_1^{-1}) x_2] = [x_2, y_2] \cup [y_2, y_2'] \cup [y_2', x_2']$ is a $(1, 100)$-quasigeodesic. Therefore $d(x_2, x_2') \ge 100$ since the path $[x_2, y_2] \cup [y_2, y_2'] \cup [y_2', x_2']$ is at least of length 200. This yields a contradiction. □

The following lemma in a certain sense mimics the lexicographical part of the order in the Nielsen method for free groups and groups acting on simplicial trees.

**Lemma 7.6.** *For any integer $m \ge 1$ there exists a constant $c = c(m) \ge 0$ such that the following holds.*

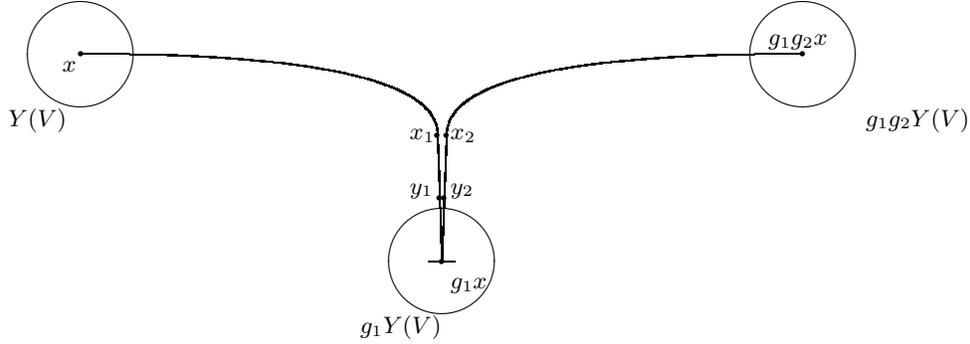

FIGURE 5. The path of a product of type $g_1 g_2$

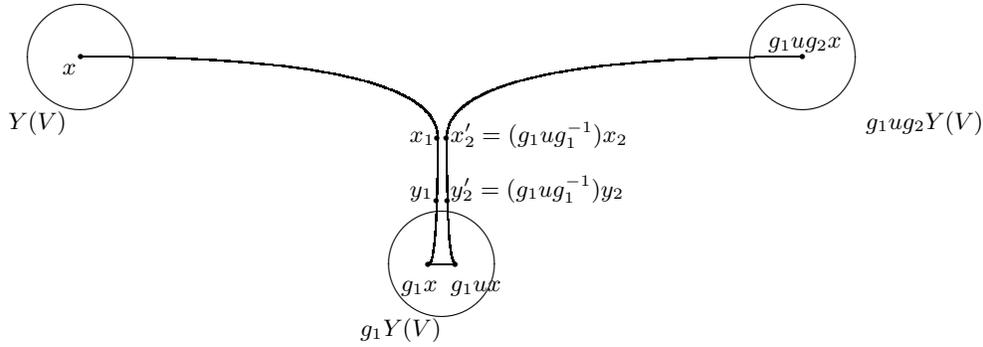

FIGURE 6. The path of a product of type $g_1 u g_2$

Let $V' \leq G$ be a non-trivial subgroup and let $M' = (V'; H')$ be a $G$-tuple with $H' = (h'_1, \ldots, h'_m)$.

Then $M$ is equivalent to a minimal tuple $M = (V; H = (h_1, \ldots, h_m))$ such that
$$l_V(hv\bar{h}) \geq l_V(h) + l_V(\bar{h}) - c$$
for any $h, \bar{h} \in \{h_1, \ldots, h_m\}^{\pm 1}$ and any $v \in V - 1$.

Furthermore the initial segment of $E_h$ of length $l_V(h) - \frac{1}{2}c$ and the terminal segment of $hvE_{\bar{h}}$ of length $l_V(\bar{h}) - \frac{1}{2}c$ lie in the 10-neighborhood of $E_{hv\bar{h}}$ if $l_V(hv\bar{h}) > 0$.

*Proof.* We first replace the $G$-tuple $(V'; H')$ by a minimal $G$-tuple $(V; H)$. In the course of the proof we modify the tuple $H$ (while keeping the same notation) by replacing an element $h_i$ by an element $v_1 h_i v_2$ for some $v_1, v_2 \in V$ and $i \in \{1, \ldots, m\}$. Since $l_V(h) = l_V(v_1 h v_2)$, it follows that these changes preserve the minimality. We will show that $c(m) := 100(m+2)$ satisfies the requirements of Lemma 7.6. The proof is conducted by induction on $m$. For $m = 0$ the statement is trivial.

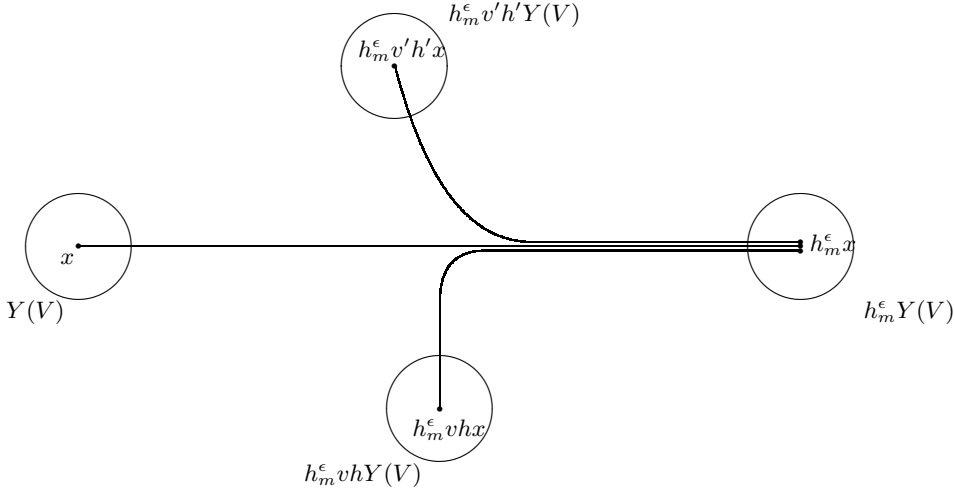

FIGURE 7. The paths of the products $h_m^\epsilon vh$ and $h_m^\epsilon v'h'$

Suppose that $m \geq 1$, that is $H = (h_1, \ldots, h_m)$. Note that $c(m) \geq 300$ and that $c(m) = c(m-1)+100$. Denote $H_{m-1} = (h_1, \ldots, h_{m-1})$. By induction we can assume that the conclusion of Lemma 7.6 holds for $H_{m-1}$ and the constant $c(m-1)$. Choose $v \in V$ such that $l_V(h_m v h_m) \leq l_V(h_m v' h_M) + 1$ for all $v' \in V$ and replace $h_m$ by $h_m v$. Note that by Lemma 7.5 this choice of $v$ is unique if there exists a $v \in V$ such that $l_V(h_m v h_m) < 2 l_V(h_m) - 200$. We denote the new element again by $h_m$.

We first show that the conclusion holds for products of type $h_m^{\epsilon_1} v h_m^{\epsilon_2}$ with $v \in V - 1$ and $\epsilon_1, \epsilon_2 \in \{-1, 1\}$. For products of type $h_m{}^\epsilon v h_m{}^{-\epsilon}$ with $\epsilon \in \{-1, 1\}$ this follows from Lemma 7.5 and the fact that the unique $v$ such that $l_V(h_m{}^\epsilon v h_m{}^{-\epsilon}) \leq l_V(h_m) + l_V(h_m) - 200$ is clearly $v = 1$. For products of type $h_m{}^\epsilon v h_m{}^\epsilon$ the conclusion of Lemma 7.6 follows from Lemma 7.5 and from our choice of $v$ in the replacement of $h$ above.

It remains to consider products of type $h_m^\epsilon vh$ with $\epsilon \in \{-1, 1\}$, $h \in H_{m-1}^{\pm 1} = H_{m-1} \cup H_{m-1}^{-1}$ and $v \in V - 1$. If $l_V(h_m^\epsilon vh) \geq l_V(h_m) + l_V(h) - c(m)$ for all products of this type, there is nothing to prove.

Suppose that there exist $v \in V - 1$, $\epsilon \in \{-1, 1\}$ and $h \in H_{m-1}^{\pm 1}$ such that $l_V(h_m^\epsilon vh) < l_V(h_m) + l_V(h) - (c(m-1) + 50)$.

**Claim.** Suppose that there exist other elements $v' \in V$ and $h' \in H_{m-1}^{\pm 1}$ such that $l_V(h_m^\epsilon v'h') < l_V(h_m) + l_V(h') - (c(m-1) + 50)$. Then $v = v'$ and, in particular, $v' \neq 1$.

The inequality $l_V(h_m^\epsilon v'h') < l_V(h_m)+l_V(h')-(c(m-1)+50)$ implies that the terminal segment of $E_{h_m^\epsilon}$ and the initial parts of $h_m^\epsilon v' E_{h'}$ of length $\frac{1}{2}(c(m-1)+50)$ are 20-close. Analogously it follows from $l_V(h_m^\epsilon vh) < l_V(h_m)+l_V(h)-(c(m-1)+50)$ that the terminal segment of $E_{h_m^\epsilon}$ and the initial parts of $h_m^\epsilon v E_h$ of length $\frac{1}{2}(c(m-1)+50)$ are 20-close. Together this implies that the initial parts of $h_m^\epsilon v E_h$ and of $h_m^\epsilon v' E_{h'}$ of length $\frac{1}{2}(c(m-1)+50)$ are 40-close, see Figure 7. This implies that

$$l_V((h^{-1}v^{-1}h_m^{-\epsilon})(h_m^\epsilon v'h')) = l_V(h^{-1}(v^{-1}v')h') \leq l_V(h) + l_V(h') - c(m-1),$$

which contradicts the induction hypothesis unless $v^{-1}v' = 1$, that is unless $v = v'$. This proves the claim.

We distinguish the cases when $l_V(h_m h_m) \geq 2l_V(h_m) - c(m)$ and when $l_V(h_m h_m) < 2l_V(h_m) - c(m)$. Recall that $v$ is the unique element such that $l_V(h_m^\epsilon vh) < l_V(h_m) + l_V(h) - (c(m-1)+50)$ for some $h \in H_{m-1}^{\pm 1}$.

**Case 1:** Suppose that $l_V(h_m h_m) \geq 2l_V(h_m) - c(m)$.

Note that the choice of $h_m$ above implies that $l_V(h_m u h_m) \geq 2l_V(h_m) - c(m)$ for all $u \in V$. We replace $h_m$ by $h_m v$ if $\epsilon = 1$ and by $v^{-1}h_m$ if $\epsilon = -1$ and denote it again by $h_m$. We clearly have $l_V(h_m u h_m) \geq 2l_V(h_m) - c(m)$ for all $u \in U$. It follows from the claim that $l_V(h_m^\epsilon uh) \geq l_V(h_m) + l_V(h) - (c(m-1)+50) \geq l_V(h_m) + l_V(h) - c(m)$ for all $u \in V - 1$ and $h \in H_{m-1}^{\pm 1}$. Again it follows from Lemma 7.5 that $l_V(g_m^\epsilon v g_m^{-\epsilon}) \geq 2l_V g_m - c(1) \geq 2l_V(g_m) - c(m)$ for $v \in V - 1$ since $v = 1$ is the only element of $V$ such that $l_V(g_m^\epsilon v g_m^{-\epsilon}) \leq 2l_V(h_m) - 200$ and since $c(m) \geq 300$. After making the same change to $h_m$ that correspond to products of type $h_m^{-\epsilon} vh$ the assertion follows.

**Case 2:** Suppose that $l_V(h_m h_m) \leq 2l_V(h_m) - c(m)$.

We replace $h_m$ by $v^{-1}h_m v$. It follows as above that $l_V(h_m^\epsilon vh) \geq l_V(h_m) + l_V(h) - (c(m-1)+50)$ for all $v \in V - 1$ and $h \in H_{m-1}^{\pm 1}$. Since $l_V(v^{-1}h_m vuh) = l_V(h_m vuh)$ it follows as in Case 1 that after this replacement $l_V(h_m uh) \geq l_V(h_m) + l_V(h) - c(m)$ for all $u \in V - 1$ and $h \in H_{m-1}^{\pm 1}$. Note that conjugation of $h_m$ by $v$ preserves the fact that $l_V(h_m h_m) \leq l_V(h_m v h_m) + 1$ for all $v \in V$. Hence it follows from Lemma 7.5 that $l_V(h_m^\epsilon v h_m^{\epsilon'}) \geq 2l_V(h_m) - c(m)$ for all $v \in V - 1$ and $\epsilon, \epsilon' \in \{-1, 1\}$. It remains to show that $l_V(h_m^{-\epsilon} vh) \geq l_V(h_m) + l_V(h) - c(m)$ for all $v \in V - 1$ and $h \in H_{m-1}^{\pm 1}$.

Suppose that there exist $v \in V - 1$ and $h \in H_{m-1}^{\pm 1}$ such that $l_V(h_m^{-\epsilon} vh) < l_V(h_m) + l_V(h) - c(m)$. In this case the same argument as in the proof of the claim shows that $l_V(h_m^\epsilon vh) < l_V(h_m) + l_V(h) - (c(m)-50) < l_V(h_m) + l_V(h) - c(m-1)$ which yields a contradiction to the induction hypothesis. □

In order to give the proof of Proposition 7.2 we need the following statement, which is a relative version of Theorem 2 of [32]. Recall that $x \in X(V)$ and that $E_g \subset W_g = [x, gx]$.

**Proposition 7.7.** *Let $m \geq 1$ be an integer. There exist constants $c_1 = c_1(m)$ and $c_2 = c_2(m)$ with the following property.*

*Suppose $(X, d)$ is a strongly geodesic 1-hyperbolic $G$-space. Let $V \subset G$ be a non-trivial subgroup and let $H = (h_1, \ldots, h_m) \in G^m$ be an m-tuple of elements of $G$. Then one of the following holds:*

(1) *The tuple $H$ is Nielsen-equivalent to a tuple $\bar{H} = (\bar{h}_1, \ldots, \bar{h}_m)$ such that for some $i \in \{1, \ldots, m\}$ either $l_V(\bar{h}_i) \leq c_1$ or that $d(y, \bar{h}_i y) \leq c_1$ for some $y \in X$.*

(2) *The subgroup $U := \langle H \rangle \leq G$ is free on $(h_1, \ldots, h_m)$ and for any freely reduced product $u = h_{i_1}^{\epsilon_1} \cdots h_{i_k}^{\epsilon_k}$ the following hold:*

   (a) $l_V(u) \geq \frac{1}{2}(l_V(h_{i_1}) + l_V(h_{i_k})) - c_2$

(b) *The initial segment of $E_{h_{i_1}^{\epsilon_1}}$ of length $\frac{1}{2}l_V(h_{i_1}) - c_2$ lies in the 1-neighborhood of $[x, ux]$.*

(c) *The terminal segments of $h_{i_1}^{\epsilon_1} \cdots h_{i_{k-1}}^{\epsilon_{k-1}} E_{h_{i_k}^{\epsilon_k}}$ of length $\frac{1}{2}l_V(h_{i_k}) - c_2$ lies in the 1-neighborhood of $[x, ux]$.*

We postpone the proof of Proposition 7.7 to the next section and proceed with the proof of Proposition 7.2.

*Proof of Proposition 7.2.* Let $c = c(m)$ be the constant provided by Lemma 7.6 and let $c_1, c_2$ be the constants $c_1(m), c_2(m)$ of Proposition 7.7. Further choose $L > 200$ such that any $L$-local $(1, 100)$-quasigeodesic in a 1-hyperbolic space is a $(\frac{1}{2}, 200)$-quasigeodesic. We will show that Proposition 7.2 holds with the constant $R(m) = \max(c_1, c_2, 2L + 2c + 2c_2 + 500)$.

Let $M = (V; H)$ be a minimal $G$-tuple, where $H = (h_1, \ldots, h_m)$.

We can assume that the statements (a)–(c) of Proposition 7.7 (2) hold for $H$. Indeed, if not then by Proposition 7.7 after passing to an equivalent $G$-tuple either case (1) or (2) of Proposition 7.2 occurs.

We can further assume that $l_V(h_i) \geq 2L + 2c + 2c_2$, since otherwise we are already in case (2) of Proposition 7.2.

The following is an immediate consequence of Proposition 7.7:

**Observation:** For any freely reduced product $u = h_{i_1}^{\epsilon_1} \cdots h_{i_k}^{\epsilon_k}$ we can find a $(1, 5)$-quasigeodesic $\omega$ joining $x$ and $ux$ such that $\omega$ contains the initial segment of $E_{h_{i_1}^{\epsilon_1}}$ of length $c + L$ and the terminal segment of $p = h_{i_1}^{\epsilon_1} \cdots h_{i_{k-1}}^{\epsilon_{k-1}} E_{h_{i_k}^{\epsilon_k}}$ or length $c + L$. Moreover, the path $\omega$ lies in the 1-neighborhood of $[x, ux]$.

We will show that any element $g \in U * F(H) - U$ acts non-trivially on $X$, which clearly implies that $\langle U, H \rangle = U * F(H)$. In the following we will not distinguish an element of $U * F(H)$ from its image in $G$. We show that $gx \neq x$. After conjugation by an element of $U$ we can assume that $g = w_1 u_1 w_2 \cdots w_{k-1} u_{k-1} w_k u_k$ where $w_i \in F(H) - 1$ for $1 \leq i \leq k$, $u_i \in U - 1$ for $1 \leq i \leq k-1$, $u_k \in U$ and $k \geq 1$. In the following we denote $w_i u_i$ by $a_i$ for $1 \leq i \leq k$, $a_1 \cdots a_i = w_1 u_1 \cdots w_i u_i$ by $g_i$ for $1 \leq i \leq k$ and define $g_0 = 1$. In particular we have $g_k = g$.

We look at the products of type $a_i a_{i+1}$ for $1 \leq i \leq k-1$. Recall that $g_{i-1} a_i = g_i$ and that $g_i a_{i+1} = g_{i+1}$.

It follows easily from the observation above and Lemma 7.6 that we can find a $L$-local $(1, 100)$-quasigeodesic $[g_{i-1} x, p_i] \cup [p_i, g_{i+1} x]$ from $g_{i-1} x$ to $g_{i+1} x$ such that $p_i$ lies in the 2-neighborhood of $[g_{i-1} x, g_i x]$ and of $[g_{i-1} x, g_i x] = [g_i x, g_i g_{i+1} x]$ and such that $[g_{i-1} x, p_i]$ contains $g_{i-1} E_{w_i}$ and $[p_i, g_{i+1} x]$ contains $g_i E_{w_{i+1}}$ with their initial and terminal segments of length $c$ removed.

It is clear that $d(p_i, p_{i+1}) \geq 2L$ and that the path $[p_i, p_{i+1}] \cup [p_{i+1}, p_{i+2}]$ is a $(1, 100)$-quasigeodesic. It follows that the path that $[x, p_1] \cup [p_1 \cup p_2] \cup \cdots \cup [p_{k-2}, p_{k-1}] \cup [p_{k-1}, gx]$ is a $L$-local $(1, 100)$-quasigeodesic and therefore a $(\frac{1}{2}, 200)$-quasigeodesic. Since this path of length more than 200, it follows that $x \neq gx$ and therefore $g \neq 1$. □

## 8. The proof of Proposition 7.7

In this section we will indicate how to prove Proposition 7.7 using the ideas and techniques of [32]. The only change from the proof of the main theorem of [32] is that we consider a different length function. In [32] we were dealing with the distance $|g|_x = d(x, gx)$, that is the length with respect to a base point. Here we are dealing with the length function $l_V(g) = d(X(V), gX(V))$, that is the length with respect to a quasiconvex subset. We will not repeat all the details in the new setting but rather indicate why all arguments go through without change. Recall that $x \in X(V)$. We show the following:

**Proposition 8.1.** *Let $m$ and $c$ be positive numbers. Then there exist numbers $d_1 = d_1(m,c)$, $d_2 = d_2(m)$, $d_3 = d_3(m)$, $d_4 = d_4(m,c)$ and $d_5 = d_5(m,c)$ such that the following hold.*

*Any minimal m-tuple $H = (h_1, \ldots, h_m) \in G^m$ is either Nielsen-equivalent to an m-tuple $H' = (h'_1, \ldots, h'_m)$ such that*

$$l_V(h'_1) \leq d_1 \text{ or } d(y, h'_1 y) \leq d_1 \text{ for some } y \in X$$

*or $U = \langle H \rangle \leq G$ is freely generated by $H$ and the following hold:*

(1) *For any $u = h_{i_1}^{\epsilon_1} \cdots h_{i_k}^{\epsilon_k} \in U$, the segment $[x, ux]$ lies in the $d_2$-neighborhood of $[x, h_{i_1}^{\epsilon_1} x] \cup h_{i_1}^{\epsilon_1}[x, h_{i_2}^{\epsilon_2} x] \cup \cdots \cup h_{i_1}^{\epsilon_1} \cdots h_{i_{k-1}}^{\epsilon_{k-1}}[x, h_{i_k}^{\epsilon_k} x]$. This implies in particular that $[x, ux]$ is contained in the a-neighborhood of $UX(V)$ where $a = \max_{i=1,\ldots,m} (\frac{1}{2} l_V(h_i) + d_2)$.*

(2) *For any freely reduced $u = h_{i_1}^{\epsilon_1} \cdots h_{i_k}^{\epsilon_k} \in U$ we have that $l_V(u) \geq l_V(h_{i_j}) - d_3$.*

(3) *For any freely reduced product $u = h_{i_1}^{\epsilon_1} \cdots h_{i_k}^{\epsilon_k} \in U$ we have that the initial segment of the essential part $E_{h_{i_1}^{\epsilon_1}}$ of length $\frac{1}{2} l_V(h_{i_1}) - d_4$ lies in the 1-neighborhood of $[x, ux]$.*

(4) *If $S$ is a subsegment of geodesic segment $[x, ux]$ for some $u \in U$, where the length of $S$ is greater than $d_5$, then $S$ intersects non-trivially the b-neighborhood of $UX(V)$ with $b = \max_{i=1,\ldots,m} (\frac{1}{2} l_V(h_i) - c)$.*

Since any $G$-tuple $(V, H)$ is equivalent to a minimal tuple, it is clear that Proposition 7.7 is an immediate consequence of Proposition 8.1.

As in [32] we can prove Proposition 8.1 by induction on $m$. The case $m = 0$ yields the empty statement and is therefore true. For the remainder of this section we fix two constants $N_1$ and $L$ such that $N_1 > L > 0$ and that any $N_1$-local $(1, 500)$-quasigeodesic is a $(L, L)$-quasigeodesic. The same argument as in [32] shows that Proposition 8.1 follows from the following statement.

**Proposition 8.2.** *Let $T \geq N_1$ be a constant and suppose that Proposition 8.1 holds for $m - 1$. Then there exist numbers $c_1 = c_1(m, T)$, $c_2 = c_2(m)$, $c_3 = c_3(m)$ and $c_4 = c_4(m)$ such that the following hold.*

*Any minimal m-tuple $M = (h_1, \ldots, h_m) \in G^m$ with $l_V(h_i) \leq l_V(h_m)$ for $1 \leq i \leq m-1$ is either Nielsen-equivalent to a tuple $M' = (h'_1, \ldots, h'_m)$ such that*

$$l_V(h'_1) \leq c_1 \text{ or } d(y, h'_1 y) \leq c_1 \text{ for some } y \in X$$

*or $U = \langle H \rangle$ is freely generated by $H$ and the following hold:*

(1) *For any $u = h_{i_1}^{\epsilon_1} \cdots h_{i_k}^{\epsilon_k} \in U$, the segment $[x, ux]$ lies in the $d_2$-neighborhood of $[x, h_{i_1}^{\epsilon_1} x] \cup h_{i_1}^{\epsilon_1}[x, h_{i_2}^{\epsilon_2} x] \cup \cdots \cup h_{i_1}^{\epsilon_1} \cdots h_{i_{k-1}}^{\epsilon_{k-1}}[x, h_{i_k}^{\epsilon_k} x]$.*

(2) *For any freely reduced $u = h_{i_1}^{\epsilon_1} \cdots h_{i_k}^{\epsilon_k} \in U$ we have that $l_V(u) \geq l_V(h_{i_j}) - c_3$.*

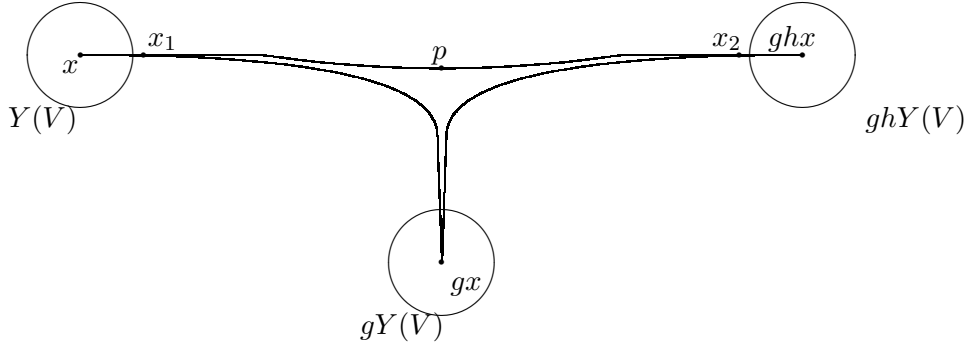

Figure 8. The path corresponding to a product of type $gh$ with cancellation

(3) *For any freely reduced product $u = h_{i_1}^{\epsilon_1} \cdots h_{i_k}^{\epsilon_k} \in U$ the initial segment of the essential part $E_{h_{i_1}^{\epsilon_1}}$ of length $\frac{1}{2} l_V(h_{i_1}) - c_4$ lies in the 1-neighborhood of $[x, ux]$.*

(4) *For any product of type $w = g_1 h_m^{\eta_1} g_2 h_m^{\eta_2} \cdots h_m^{\eta_{l-1}} g_l h_m^{\eta_l} g_{l+1}$ with $\eta_i \in \{-1, 1\}$ for $1 \leq i \leq l$, $g_i \in \langle h_1, \ldots, h_{m-1} \rangle$ for $1 \leq i \leq l+1$ and $g_i \neq 1$ if $\eta_i = -\eta_{i+1}$ for $2 \leq i \leq l$ the following holds:*
*The $2L + 100$-neighborhood of every subsegment of $[x, wx]$ of length at least $10T$ contains a subsegment of length at least $T$ of either a segment of type $g_1 h_m^{\eta_1} g_2 h_m^{\eta_2} \cdots h_m^{\eta_{i-1}} g_i [x, h_m^{\eta_i} x]$ or $g_1 h_m^{\eta_1} g_2 h_m^{\eta_2} \cdots h_m^{\eta_{i-1}} [x, g_i x]$ for some $i$.*

The proof of Proposition 8.2 is identical to the proof of Proposition 11 of [32]. Namely, we need to study cancellations (relative $V$) in products of type $gh$ and $ghf$, where $g, h, f \in G$. We also need to define *stable parts* of products of type $h_m^\eta g$ and $g h_m^\eta$ where $g \in \langle h_1, \ldots, h_{m-1} \rangle$ and $\eta \in \{-1, 1\}$.

We first study cancellation in products of length two.

**Lemma 8.3.** *Suppose that $g, h \in G$ and that $l_V(g), l_V(h), l_V(gh) \geq 20$. Then the product $gh$ has one of the following types:*

(1) *If $l_V(gh) < l_V(g) + l_V(h)$, we have the picture as in Figure 8. That is, we can write $E_{gh} = [x_1, p] \cup [p, x_2]$ so that the 10-neighborhood of $[x_1, p]$ contains the initial segment of $E_g$ of length $l_V(g) - \frac{1}{2}(l_V(g) + l_V(h) - l_V(gh)) - 10 = \frac{1}{2}(l_V(g) - l_V(h) + l_V(gh)) - 10$*
*and the 10-neighborhood of $[p, x_2]$ contains the terminal segment of $gE_h$ of length $l_V(h) - \frac{1}{2}(l_V(g) + l_V(h) - l_V(gh)) - 10 = \frac{1}{2}(l_V(h) - l_V(g) + l_V(gh)) - 10$. In this case we say that the product is of type 1.*

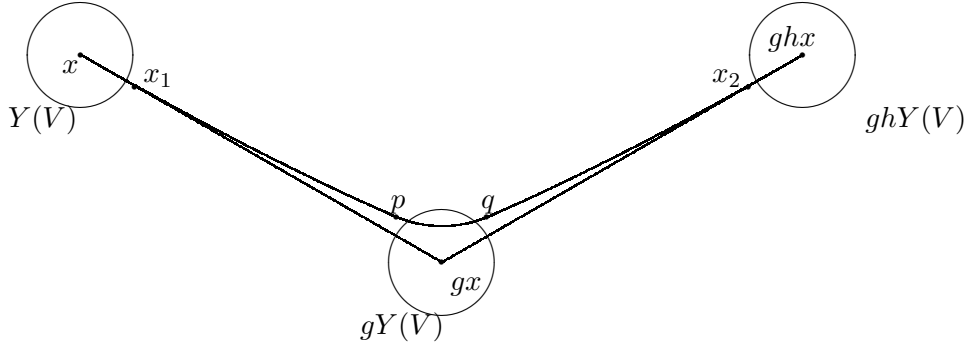

FIGURE 9. The path corresponding to a product of type $gh$ without cancellation

(2) If $l_V(gh) \geq l_V(g) + l_V(h)$, we have a picture as in Figure 9. That is, we can write $E_{gh} = [x_1, p] \cup [p, q] \cup [q, x_2]$ so that the 10-neighborhood of $[x_1, p]$ contains $E_g$ and the 10-neighborhood of $[q, x_2]$ contains $gE_h$. In this case we say that the product is of type 2.

As a consequence of Lemma 8.3 we obtain the following two lemmas which correspond to the two parts of Lemma 17 of [32]. They will guarantee that the stable parts defined below always exist.

**Lemma 8.4.** *Suppose that $g \in G$ with $l_V(g) \geq 20$. Let $N \geq N_1$. Then either $d(y, gy) \leq 20$ for some $y \in X$ or $l_V(g) \leq 2N + 2d_2 + 100$ or we have $E_{g^2} = [x_1, p] \cup [p, q] \cup [q, x_2]$ where $[x_1, p]$ lies in the 10-neighborhood of $E_g$ and $[q, x_2]$ lies in the 10-neighborhood of $gE_g$ and $d(x_1, p) \geq \frac{1}{2}l_V(g) + N + d_2$ and $d(q, x_2) \geq \frac{1}{2}l_V(g) + N + d_2$.*

*Proof.* Suppose that $l_V(g) > 2N + 2d_2 + 100$. If the product $gg$ is of Type 2, we clearly get the last conclusion. Suppose that the last conclusion does not hold. It follows that the product $gg$ is of Type 1. In this case it is clear that the midpoint $y$ of $E_g$ and the midpoint $y' = gy$ of $gE_g$ are at distance at most 20 from each other, which gives the first conclusion. □

**Lemma 8.5.** *Suppose that $H = (h_1, \ldots, h_m)$ is minimal, that $l_V(h_i) \leq l_V(h_m)$ for $1 \leq i \leq m-1$, that $g \in \langle h_1, \ldots, h_{m-1} \rangle$, that $\eta \in \{-1, 1\}$ and that the conclusion of Proposition 8.1 holds for $m-1$ with $d_2 = d_2(m-1)$. Then we have*

$$E_{gh_m^\eta} = [x_1, p] \cup [p, q] \cup [q, x_2]$$

*where*

(1) *$[x_1, p]$ lies in the 10-neighborhood of $E_g$ and $d(x_1, p) \geq \frac{1}{2}l_V(g) - 20$.*
(2) *$[q, x_2]$ lies in the 10-neighborhood of $gE_{h_m^\eta}$ and $d(q, x_2) \geq \frac{1}{2}l_V(h_m) - (20 + d_2)$.*

*Proof.* The existence of the path $E_{gh_m^\eta} = [x_1, p] \cup [p, q] \cup [q, x_2]$ such that $[x_1, p]$ lies in the 10-neighborhood of $E_g$ and that $[q, x_2]$ lies in the 10-neighborhood of $gE_{h_m^\eta}$ follows from the discussion of the paths associated to products of length 2; here we set $p = q$ if the path is of Type 1. Now $d(x_1, p) \geq \frac{1}{2}l_V(g) - 20$ since otherwise $l_V(gh_m^\eta) < l_V(h_m) - 1$ which contradicts the minimality of $H$. The last statement follows precisely as the proof of Lemma 17 (2) of [32] since by assumption (1) of Proposition 8.2 holds for $m-1$. □

The main component of the proof in [32] that has to be mimicked is the definition of the *stable part*. Note that as in [32] the stable part $S_w$ of $w = gh_m^\eta$ applies to the product $w = gh_m^\eta$ and not simply to the element $w$.

**Definition 8.6.** *[Stable part]* Let $N \geq N_1$. Suppose that $w$ is the product $w = gh_m^\eta$ with $g \in \langle h_1, \ldots, h_{m-1} \rangle$ and $\eta \in \{-1, 1\}$ and that $l_V(h_m) \geq 4N + 2d_2(m-1) + 2d_4(m-1, d_2(m-1) + 20) + 100$. We assign to $w$ the stable part $S_w$ of $w$ relative $N$ as follows:

(1) If the terminal segment of length $\frac{1}{2}l_V(h_m) + N + d_2(m-1)$ of $gE_{h_m^\eta}$ lies in the 2-neighborhood of $E_{hg_m^\eta}$, then we put

$$S_w := [s, t] \subset gE_{h_m^\eta}$$

where $[s, t] \subset gE_{h_m^\eta} = [x_1, x_2]$ are chosen such that $d(x_2, s) = \frac{1}{2}l_V(h_m) + N + d_2(m-1)$ and $d(x_2, t) = \frac{1}{2}l_V(h_m) + d_2(m-1)$. (In Figure 10 $m$ is the midpoint of $gE_{h_m^\eta} = [x_1, x_2]$, $d(m, t) = d_2(m-1)$ and $d(s, t) = N$.)

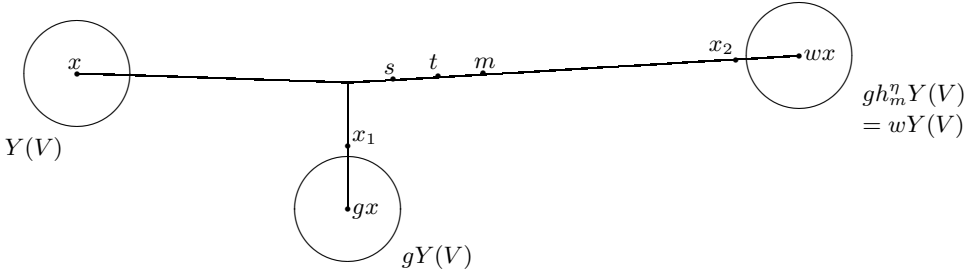

FIGURE 10. Stable part: first case.

(2) If the terminal segment of length $\frac{1}{2}l_V(h_m) + N + d_2(m-1)$ of $gE_{h_m^\eta}$ does not lie in the 2-neighborhood of $E_{hg_m^\eta}$ then we must clearly have a product of Type 1 and we choose the point $p$ accordingly. We then put

$$S_w := [s, t] \subset E_g$$

where $[s, t] \in E_g = [y_1, y_2]$ are chosen such that $d(s, p) = d_4(m-1, d_2(m-1) + 20)$ and $d(t, p) = d_4(m-1, d_2(m-1) + 20) + N$ (see Figure 11). The fact that the segment $[s, t]$ exists follows precisely as in [32] from the minimality of and the fact that $l_V(h_m) \geq 4N + 2d_2(m-1) + 2d_4(m-1, d_2(m-1) + 20) + 100$.

We define the stable part relative $N$, denoted $S_v$, for a product $v = h_m^\eta g$ relative $N$ by considering $v$ as the inverse of the product $w = v^{-1} = h_m^{-\eta}g^{-1}$: Since $S_w = [s, t]$ lies in the 2-neighborhood of $[x, wx]$ for some geodesic $[x, wx]$, it follows that $w^{-1}[s, t]$ lies in 2-neighborhood of $[w^{-1}x, x] = [vx, x]$ for some geodesic segment $[vx, x]$. We define $S_v := w^{-1}[t, s]$, where $[t, s]$ is the geodesic segment $[s, t]$ traveled from $t$ to $s$.

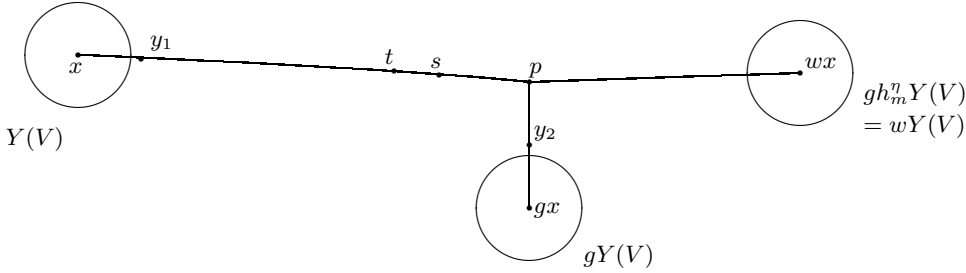

FIGURE 11. Stable part:second case.

Informally speaking, the stable part now represents part of the segment corresponding to a product of length 2 that survives as a subsegment in longer products. The same argument as in [32] gives us the following local stability result.

**Lemma 8.7.** Let $N \geq N_1$. There exists a constant $k = k(N, m)$ with the following property: Suppose that $H = (h_1, \ldots, h_m)$ is as in Proposition 8.2 and that $l_V(h_m) \geq 4N + 2d_2(m-1) + 2d_4(m-1, d_2(m-1) + 20) + 100$. Then either the first conclusion of Proposition 8.2 holds or for any $g, g_1, g_2 \in \langle h_1, \ldots, h_{m-1} \rangle$ and $\eta, \eta_1, \eta_2 \in \{-1, 1\}$ the following hold:

(1) We can write $E_{g_1 h_m^\eta g_2} = [x_1, p] \cup [p, x_2]$ so that the 20-neighborhood of $[x_1, p]$ contains $S_{g_1 h_m^\eta}$ and the 20-neighborhood of $[p, x_2]$ contains $g_1 S_{h_m^\eta g_2}$.
(2) We can write $E_{h_m^{\eta_1} g h_m^{\eta_2}} = [x_1, p] \cup [p, x_2]$ so that the 20-neighborhood of $[x_1, p]$ contains $S_{h_m^{\eta_1} g}$ and the 20-neighborhood of $[p, x_2]$ contains $h_m^{\eta_1} S_{g h_m^{\eta_2}}$.

We now define a path $\sigma_N(w)$ which is a quasigeodesic close to $S_w$, where $W$ is a product of type
$$w = g_1 h_m^{\eta_1} g_2 h_m^{\eta_2} \cdots h_m^{\eta_{l-1}} g_l h_m^{\eta_l} g_{l+1}$$
with $\eta_i \in \{-1, 1\}$ for $1 \leq i \leq l$, $g_i \in \langle h_1, \ldots, h_{m-1} \rangle$ for $1 \leq i \leq l+1$ and $g_i \neq 1$ if $\eta_i = -\eta_{i+1}$ for $2 \leq i \leq l$. First define $w_i = g_i h_m^{\eta_i}$ and $v_i = h_m^{\eta_i} g_{i+1}$ for $1 \leq i \leq l$ and choose $s_i, t_i, s'_i$ and $t'_i$ such that $S_{w_i} = [s_i, t_i]$ and $S_{v_i} = [s'_i, t'_i]$. Furthermore, denote $S_i := w_1 \cdots w_{i-1} s_i$, $T_i := w_1 \cdots w_{i-1} t_i$, $S'_i = w_1 \cdots w_{i-1} g_i s'_i$ and $T'_i = w_1 \cdots w_{i-1} g_i t'_i$. This implies that $w_1 \cdots w_{i-1} S_{w_i} = [S_i, T_i]$ and $w_1 \cdots w_{i-1} h_i S_{v_i} = [S'_i, T'_i]$.

As in [32], we now conclude that the path

$$\sigma_N(w) = [x, S_1] \cup [S_1, T_1] \cup [T_1, S'_1] \cup [S'_1, T'_1] \cup [T'_1, S_2] \cup [S_2, T_2] \cup \cdots \cup [T_2, S'_2] \cup [S'_2, T'_2] \cup \cdots \cup [S'_l, T'_l] \cup [T'_l, wx]$$

is a $N$-local $(1, 500)$-quasigeodesic and therefore $L$-Hausdorff close to any geodesic segment $[x, wx]$.

The following simple fact follows immediately from the definition of hyperbolicity:

**Lemma 8.8.** *Let $(X,d)$ be a 1-hyperbolic geodesic space and suppose that $[x,y]$ is a geodesic segment and that $c = d(z,[x,y])$. Suppose further that $d(x,z) \geq 10+c$. Then the initial segment of $[x,z]$ of length $d(x,z) - (c+10)$ lies in the 1-neighborhood of $[x,y]$.* □

*Proof of Proposition 8.2* All conclusions of Proposition 8.2 except (3) follow exactly as in [32]. We therefore only have to verify the third claim of Proposition 8.2.

We will show that part (3) of Proposition 8.2 holds for $c_4(m) = \max(d_4(m-1) + d_3(m-1) + d_4(m-1, d_2(m-1)+20) + L + 31, d_2(m-1) + L + 10)$. Let $u = h_{i_1}^{\epsilon_1} \cdots h_{i_k}^{\epsilon_k} \in U$ be a freely reduced product. We rewrite $u$ as the product $g_1 h_m^{\eta_1} g_2 h_m^{\eta_2} \cdots h_m^{\eta_{l-1}} g_l h_m^{\eta_l} g_{l+1}$ by joining all sub-products with factors in $\{h_1, \ldots, h_{m-1}\}^{\pm 1}$ to element $g_i \in \langle h_1, \ldots, g_{m-1} \rangle$. We choose the same notations as above for the quasigeodesic $\sigma_N(w)$. We distinguish the cases that $g_1 \neq 1$ and that $g_1 = 1$.

**Case 1:** $g_1 = 1$. This means that $g_1 h_m^{\eta_1} = h_m^{\eta_1} = h_{i_1}^{\epsilon_1}$. In this case the stable part $[S_1, T_1]$ of the product $g_1 h_m^{\eta_i}$ is by definition of the stable part a subsegment of $E_{h_m^{\eta_1}} = E_{h_{i_1}^{\epsilon_1}}$ which lies at most in distance $d_2(m-1)$ from the midpoint of $E_{h_m^{\eta_1}}$. If follows that the 1-neighborhood of $[x, S_1] \cup [S_1, T_1]$ contains the initial segment of $E_{h_m^{\eta_i}}$ of length at least $\frac{1}{2}l_V(h_m) - d_2(m-1)$. It follows from Lemma 8.8 that the initial segment of $E_{h_m^{\eta_i}} = E_{h_{i_1}^{\epsilon_1}}$ of length at least $\frac{1}{2}l_V(h_m) - d_2(m-1) - L - 10 = \frac{1}{2}l_V(h_1) - d_2(m-1) - L - 10$ lies in the 1-neighborhood of $[x, ux]$ since $[x, S_1] \cup [S_1, T_1]$ lies in the $L$-neighborhood of $[x, ux]$. This proves the claim since by assumption $c_4 \geq d_2(m-1) + L + 10$

**Case 2:** $g_1 \neq 1$. In this case the stable part $[S_1, T_1]$ of $g_1 h_m^{\eta_1}$ is either a subsegment of $[x, g_1 x]$ or of $g_1[x, h_m^{\eta_1} x]$. In both cases it follows from the definition of the stable part and from Lemma 8.5 that the path $[x, S_1] \cup [S_1, T_1]$ contains $[x, y_1] \cup [y_1, y_2] \subset [x, g_1 x]$ where $[y_1, y_2]$ the initial segment of $E_{g_1}$ of length $\frac{1}{2} l_V(g_1) - d_4(m-1, d_2(m-1)+20) - 20$. By induction we know that the initial part of $E_{h_{i_1}^{\epsilon_1}}$ of length $\frac{1}{2}l_V(h_{i_1}) - d_4(m-1)$ lies in the 1-neighborhood of $[x, g_1 x]$. Since also by induction $l_V(g_1) \geq l_V(h_{i_1}) - d_3(m-1)$ it follows that the initial segment of $E_{h_{i_1}^{\epsilon_1}}$ of length $\frac{1}{2}l_V(h_{i_1}) - d_4(m-1) - d_3(m-1) - d_4(m-1, d_2(m-1)+20) - 20$ lies in the 1 neighborhood of the $[x, y_1] \cup [y_1, y_2]$. It follows that the initial segment of $E_{h_{i_1}^{\epsilon_1}}$ of length $\frac{1}{2}l_V(h_{i_1}) - d_4(m-1) - d_3(m-1) - d_4(m-1, d_2(m-1)+20) - 20$ lies in the 2-neighborhood of $[x, S_1] \cup [S_1, T_1]$. As in the first case we argue that the initial segment of $E_{h_{i_1}^{\epsilon_1}}$ of length $\frac{1}{2}l_V(h_{i_1}) - d_4(m-1) - d_3(m-1) - d_4(m-1, d_2(m-1)+20) - L - 31$ lies in the 1-neighborhood of $[x, ux]$ which proves the assertion because of the choice of $c_4$.

## 9. Final deduction

In this section we finish the proof of Theorem 2.4. We generalize the definition of minimality to $G$-tuples with non-trivial hyperbolic set.

**Definition 9.1** (Minimal $G$-tuples). *Let $(X,d)$ be a 1-hyperbolic strongly geodesic $G$-space. Let $M = (U_1, \ldots, U_n; H)$ be a $G$-tuple, where $n \geq 1$, $m \geq 0$, $H = (h_1, \ldots, h_m) \in G^m$ and $U_i$ are non-trivial subgroups of $G$.*

(1) *We say that $M$ is* weakly minimal *if*
$$d(X(U_1), \ldots, X(U_n)) \leq d(X(U_1'), \ldots, X(U_n')) + 1.$$

for any $G$-tuple $M' = (U'_1, U'_2, \ldots, U'_n; H')$ which is equivalent to $M$.

(2) We say that $M$ is minimal *if $M$ is weakly minimal and if*

$$|H|_V \leq |H'|_{V'} + 1$$

for any weakly minimal $G$-tuple $M' = (U_1, U_2, \ldots, U_n, H)$ which is equivalent to $M$. Here $V = \langle U_1 \cup \cdots \cup U_n \rangle$ and $V' = \langle U'_1 \cup \cdots \cup U'_n \rangle$.

Note that this notion of minimality coincides with the definition given in Section 6 if $m = 0$ and $n \geq 1$ and with the definition given in Section 7 if $n = 1$. Although for the purposes of this paper it is not necessary to define minimality for the case $n = 0, m \geq 1$, this case was handled in [32]. Once again, since $\delta > 0$, minimal elements exist in every equivalence class.

We can now obtain the proof of our main technical result.

*Proof of Theorem 2.4.* It suffices to prove Theorem 2.4 for $\delta = 1$ since if $(X, d)$ is $\delta$-hyperbolic (with $\delta > 0$), then the space $(X, d/\delta)$ is 1-hyperbolic. Thus we will assume that $(X, d)$ is a 1-hyperbolic strongly geodesic $G$-space.

Let $C(k)$ be the constant from Proposition 6.3 and let $R(k)$ be the constant provided by Proposition 7.2. Let $C'(k)$ be the constant provided by Theorem 1 in [32]. Put $K(k) := 2k(2C(k) + R(k) + 2) + 2C(k) + R(k) + C'(k)$.

If $n = 0$ and $M = (; H)$, the statement of Theorem 2.4 follows from the main result of [32]. If $m = 0$ and $H = (-)$ then the statement of Theorem 2.4 follows from Proposition 6.3. If $n = 1$ the statement of Theorem 2.4 follows from Proposition 7.2.

Therefore we will assume that $m \geq 1$ and $n \geq 2$.

Suppose that $M = (U_1, \ldots, U_n; H)$ is a minimal $G$-tuple, where $n \geq 2$, $H = (h_1, \ldots, h_m) \in G^m$, $m \geq 1$ and $m + n \leq k$. Recall that by the definition of $G$-tuple $U_i$ are non-trivial subgroups of $G$.

Denote $X_i = X(U_i)$ and let $\Omega = \Omega(X_1, \ldots, X_n)$ be a minimal connector. Put $V = \langle U_1 \cup \cdots \cup U_n \rangle \leq G$. By Proposition 6.3 we may assume that

$$V = U_1 * \cdots * U_n$$

and that $d(X_i, X_j) \geq \delta C(k)$ for $i \neq j$ and that

$$A := \bigcup_{u \in U} u(\Omega \cup X_1 \cup \cdots \cup X_n)$$

is $C(k)$-close to $X(V)$.

We claim that the $G$-tuple $(V; H)$ is minimal. Indeed, suppose not. Then there is a pair $(V'; H')$ equivalent to $(V; H)$ such that $|H'|_{V'} < |H|_V + 1$. The definition of equivalence implies that there is a $G$-tuple $M' = (U'_1, \ldots, U'_n; H')$ equivalent to $M$ such that for some $g \in G$ we have $V' = \langle U'_1 \cup \cdots \cup U'_n \rangle$ and $U'_i = gU_ig^{-1}$, $V' = gVg^{-1}$. Thus $X(U'_i) = gX(U_i)$ and therefore $d(X(U'_1), \ldots, X(U'_n)) = d(X(U_1), \ldots, X(U_n))$. This means that the $G$-tuple $M'$ is also weakly minimal. But $|H'|_{V'} < |H|_V + 1$, which contradicts the minimality of $M$.

Thus $(V; H)$ is a minimal $G$-tuple. Therefore by Proposition 7.2 either $U = V * F(H) = U_1 * \cdots * U_n * F(H)$ or $(V; H)$ is equivalent to $(V'; H')$ such that one of the two cases of Proposition 7.2 applies to $(V'; H')$. Once again, the definitions of equivalence imply that there is $g \in G$ such $V' = gVg^{-1}$ and that for $U'_i = gU_ig$ the $G$-tuples $M$ and $M' = (U'_1, \ldots, U'_n, H')$ are equivalent. Again, we see that $d(X(U_1), \ldots, X(U_n)) = d(X(U'_1), \ldots, X(U'_n))$, so $M'$ is weakly minimal. We will assume that $g = 1$ and it will be clear from the proof that the same argument applies in general. Thus $U_i' = U_i$ and $V' = V$ (but $H$ and $H'$ are not necessarily the same).

Suppose first that case (b) of Proposition 7.2 applies to $(V'; H') = (V; H')$. Then there is $h \in H'$ and $x \in X$ with $d(x, hx) \leq R(k)$ and the conclusion of Theorem 2.4 obviously holds. Suppose now that case (a) of Proposition 7.2 applies to $(V; H')$. Thus there are $x_1, x_2 \in X(V)$ and $h \in H'$ such that $d(x_1, x_2) \leq R(k)$. Since $A$ and $X(V)$ are $C(k)$-close, there are points $x, y \in B := \Omega \cup X_1 \cdots \cup X_n$ and $v_1, v_2 \in V$ such that $d(v_1 x, h v_2 y) \leq 2C(k) + R(k)$. Replacing $h \in H'$ by $v_1^{-1} h v_2$ (which preserves the equivalence classes of both $M'$ and $(V; H')$ and the weak minimality of $M'$) we may further assume that $v_1 = v_2 = 1$. Thus there are $x, y \in B$ such that $d(x, hy) \leq 2C(k) + R(k)$.

Denote $L = (2C(k) + R(k) + 2)k$. If both $x, y$ are in the $L$-neighborhood of $X_1 \cup \cdots \cup X_n$ then there are points $x' \in X_i$ and $y' \in X_j$ with $1 \leq i, j \leq n$ and $d(x', hy') \leq 2L + 2C(k) + R(k) = K(k)$. If $i = j$ then conclusion (2) of Theorem 2.4 holds. If $i \neq j$ the conclusion (1) of Theorem 2.4 holds after replacing $U_j$ by $h U_j h^{-1}$ and therefore $X(U_i)$ by $hX(U_i)$. Thus we may assume that $x \in \Omega$ and that $d(x, X_i) > L$ for all $i = 1, \ldots, n$. Recall that $y \in B$ and $d(x, hy) \leq 2C(k) + R(k)$ for some $h \in H'$.

Consider the $\Omega$-geodesic $[x, y]_\Omega$. Let $z \in [x, y]_\Omega$ be such that $d_\Omega(x, z) = L$. Thus $z \in \Omega$ and moreover, $z$ and $x$ belong to the same tree-component $T$ of $\Omega$. Since $n \leq k$, the forest $\Omega$ has at most $k - 2$ branching points. Hence there is a subsegment $[x', z']_\Omega \subseteq [x, z]_\Omega$ such that $d_\Omega(x', z') = L/k = 2C(k) + R(k) + 2$ and the segment $[x', z']_\Omega$ contains no branching points of $\Omega$.

Removing $z'$ from $B$ disconnects $B$ into two connected components: component $B_{x'}$ containing $x'$ and component $B_y$ containing $y$. Denote $D_{x'} := (B_{x'} \cap \Omega) - [x', z']_\Omega$ and $D(y) := B_y \cap \Omega$.

Put $\Omega'' = D_{x'} \cup h^{-1} D_y \cup [x', h^{-1} y]$. Then

$$P(\Omega'') \leq P(\Omega) - (2C(k) + R(k) + 2) + 2C(k) + R(k) < P(\Omega) - 1.$$

Let $I$ be the set of all $i$ such that $X_i \subseteq B_{x'}$ and $J$ be the collection of all $i$ such that $X_i \subseteq B_y$. Put $U_i'' = U_i$ if $i \in I$ and put $U_i'' = h^{-1} U_i h$ if $i \in J$. Put $M'' = (U_1'', \ldots, U_n''; H')$. Then $M$ and $M''$ are clearly equivalent $G$-tuples. Denote $X_i'' = X(U_i'')$. Since for $i \in J$ we have $X_i'' = h^{-1} X_i$, the construction of $\Omega''$ implies that $\Omega'' \cup X_1'' \cup \cdots \cup X_p''$ is connected. We have already seen that $P(\Omega'') < P(\Omega) - 1$, which contradicts the minimality of $M$. This completes the proof of Theorem 2.4. $\square$

**Remark 9.2.** *The assumption in Theorem 2.4 that $(X, d)$ be a strongly geodesic space can be replaced by requiring $(X, d)$ to be geodesic. We imposed the "strongly geodesic" assumption in order to talk about boundaries of hyperbolic spaces in terms of equivalence classes of geodesic rays, as it is traditionally done. However, following M.Gromov, one can define $\partial X$ to be the set of equivalence classes of "sequences tending to infinity" (see the definition of sequential boundary in Section 3). Then we are no longer guaranteed that two distinct points in the boundary can be connected by a geodesic or that a point in the boundary can be connected by a geodesic to a point in the space. However, it is still true that any two distinct points in $X \cup \partial X$ can be connected by $(1, 20\delta)$-quasigeodesic. Thus one should define convex hulls as unions of $(1, 20\delta)$-quasigeodesics rather than geodesics, and the rest of the argument will go through.*

## 10. Applications to hyperbolic groups

If $G$ is a group and $Z = (g_1, \ldots, g_n)$ is an $n$-tuple of elements of $G$, we will call $n$ the *length of $Z$* and denote $n = L(Z)$. We shall also say that the empty tuple () has zero length.

**Convention 10.1** (Partitioned tuple)**.** *Let $G$ be a group. If $Y$ is an $n$-tuple of elements of $G$ we will say that $n$ is the length of $Y$. We will say that $M = (Y_1, \ldots, Y_p; H)$ is a partitioned tuple in $G$ if $p \geq 0$ and $Y_1, \ldots, Y_p, H$ are finite tuples of elements of $G$ such that at least one of these tuples has positive length and such that for any $i \geq 1$ the tuple $Y_i$ has positive length. We will call the sum of the lengths of $L(Y_1) + \cdots + L(Y_p) + L(H)$ the length of $M$ and denote it by $L(M)$.*

We will now apply Theorem 2.4 to almost torsion-free word-hyperbolic groups. First we observe that for an almost torsion-free word-hyperbolic group $G = \langle S|R\rangle$ acting on its Cayley graph $X = \Gamma(G, S)$ the set $X(U)$ is in fact close to $U$ for a quasiconvex subgroup $U$.

We will say that the group $G = \langle S|R\rangle$ is $\delta$-*hyperbolic* if the Cayley-graph $\Gamma(G, S)$ is $\delta$-hyperbolic. The constant $\delta$ only depends on $G$ and $S$ and can in fact be computed using the algorithm of P.Papasoglu [36]. The generating set $S$ will always be finite. We denote the metric in the Cayley graph by $d_S$ and define $|g|_S = d_S(1, g)$. Often we simply use $d$ instead of $d_S$ and $|g|$ instead of $|g|_S$.

**Convention 10.2.** *Till the end of this section, unless specified otherwise, we shall assume that $G = \langle S|R\rangle$ is an almost torsion-free $\delta$-hyperbolic group, where $\delta > 0$. We will also denote the Cayley graph $\Gamma(G, S)$ of $G$ by $X$.*

**Lemma 10.3.** *There exists a constant $c_1 = c_1(G, S, \delta, \epsilon) > 0$ such that the following holds:*

*Let $U$ be an $\epsilon$-quasiconvex infinite subgroup of $G$, where $\epsilon > 0$ is an integer. Then $U$ and $X(U)$ are $c_1$-Hausdorff close (where $X(U)$ is defined relative $\delta$).*

*Proof.* Note first that $U$ has a non-trivial element $u_0 \in U, u_0 \neq 1$ such that $d(1, u_0) \leq 3\epsilon$. Indeed, since $U$ is infinite, there is $u' \in U$ such that $d(1, u') > 3\epsilon$. Consider a geodesic $[1, u']$ and let $x \in [1, u']$ be such that $d(1, x) = 2\epsilon$. Since $U$ is $\epsilon$-quasiconvex in $\Gamma(G, S)$, there is $u_0 \in U$ with $d(x, u) \leq \epsilon$. Then by the triangle inequality

$$0 < \epsilon \leq d(1, u_0) \leq 3\epsilon$$

Thus $u_0 \neq 1, u_0 \in U$ and $d(1, u_0) \leq 3\epsilon$, as required.

Recall that the set $X(U)$ is $4\delta$-quasiconvex. Let $p \in X(U)$ be such that $d(1, p) = d(1, X(U))$. Then by Lemma 4.4 the path $[1, p] \cup [p, u_0 p] \cup u_0[p, 1]$ is $(1, 100\delta)$-quasigeodesic. Since $d(1, u_0) \leq 3\epsilon$, we conclude that $d(1, p) \leq 3\epsilon + 100\delta$.

Then for any $u \in U$ we have $up \in X(U)$ and $d(u, up) = d(1, p) \leq 3\epsilon + 100\delta$. Therefore $U$ is contained in the $(3\epsilon + 100\delta)$-neighborhood of $X(U)$.

Let $t' \in X(U)$. By Lemma 4.10 the sets $X(U)$ and $X_U$ are $c$-close and hence there is $t \in X_U$ such that $d(t, t') \leq c$. Here $c = c(\#S, \delta)$ is the constant provided by Lemma 4.10.

Then $t$ lies on a biinfinite geodesic $[a, b]$ for some $a, b \in \Lambda U \subseteq \partial G$. Thus both $a$ and $b$ are approximated by elements of the $U$-orbit of an arbitrary point of $\Gamma(G, S)$, for example the $U$-orbit of 1. Therefore $a$ and $b$ are approximated by elements of $U$. Recall that $t \in [a, b]$. Hence there exist $u_1, u_2 \in U$ such that for some $y \in [u_1, u_2]$ we have $d(y, t) \leq 2\delta$. Since $U$ is $\epsilon$-quasiconvex, there is $u \in U$ with $d(y, u) \leq \epsilon$. Hence $d(t, u) \leq 2\delta + 2\epsilon$ and $d(t', u) \leq c + 2\delta + 2\epsilon$. Since $t' \in X(U)$ was chosen arbitrarily, this implies that $X(U)$ is contained in the $(c+2\delta+2\epsilon)$-neighborhood of $U$.

Since we already know that $U$ is contained in the $(3\epsilon+100\delta)$-neighborhood of $X(U)$, the sets $U$ and $X(U)$ are $(c+100\delta+3\epsilon)$-close, that is the claim holds for $c_1 = c+100\delta+3\epsilon$. $\square$

**Lemma 10.4.** *For any $C \geq 0$ there is a constant $c_2 = c_2(G,S,\delta,C)$ such that for any non-trivial finite subgroup $U \leq G$ contained in the $C$-ball around $1$ in $\Gamma(G,S)$ the sets $U$ and $X(U)$ are $c_2$-Hausdorff close (where $X(U)$ is defined relative $\delta$).*

*Proof.* Choose $\delta$ such that $G = \langle S|R \rangle$ is $\delta$-hyperbolic. Note first that since $G$ is a word-hyperbolic group, it has only finitely many conjugacy classes of subgroups of finite order. By a theorem of V.Gerasimov and O.Bogopolskii [10] any subgroup of finite order in $G$ is conjugate to a subgroup contained in the ball of radius $2\delta + 1$ around $1 \in G$ in the Cayley graph $\Gamma(G,S)$. Let $U \leq G$ be a non-trivial finite subgroup of $G$ such that $|u|_S \leq 2\delta + 1$ for any $u \in U$. Let $g \in G$ be such that for some $u \in U, u \neq 1$ we have $d(g, ug) \leq 100\delta$. Since the number of such subgroups is finite, and for each such subgroup the diameter of $X(U)$ is finite by Lemma 4.11, there is a constant $c' = c'(G,\delta,S)$ such that $U$ and $X(U)$ are $c'$-close

Suppose now $U$ is an arbitrary finite subgroup of $G$ contained in the ball of radius $C$ around $1$. Then $U = hVh^{-1}$ for some $h \in G$ and some finite subgroup $V$ of $G$ contained in the ball of radius $2\delta + 1$ and with $|h|_S \leq K$, where $K = K(C,\delta,G,S)$ is some constant. Then $X(U) = hX(V)$. Hence $X(U)$ is contained in the ball or radius $K + c'$ around $1$. Also $U$ is contained in the ball of radius $2K + 4\delta + 2$. Therefore $U$ and $X(U)$ are $c_2$-Hausdorff close, where $c_2 = 2K + 4\delta + 2 + 2c'$. $\square$

**Lemma 10.5.** *For any integers $C \geq 0$ and $n \geq 1$ there is a constant $c_3 = c_3(G,S,\delta,n,C)$ with the following property. Suppose $U \leq G$ is a non-trivial quasiconvex subgroup $U \leq G$ generated by a set $Y$ with at most $n$ elements such that $Y$ is contained in the $C$-ball around $1$ in $\Gamma(G,S)$. Then the sets $U$ and $X(U)$ are $c_3$-Hausdorff close.*

*Proof.* Let $E_n$ be the maximum of quasiconvexity constants of infinite quasiconvex subgroups generated by subsets with at most $n$-elements from the $C$-ball around $1$ in $\Gamma(G,S)$. Put $c_3 := c_1(G,S,\delta,E_n) + c_2(G,S,\delta,C)$ The statement of Lemma 10.5 now follows directly from Lemma 10.3 and Lemma 10.4. $\square$

We establish two lemmas that correspond to two different cases in the inductive step of the proof of Theorem 1.3.

**Lemma 10.6.** *For any integers $K \geq 0$ and $n_1, n_2 \geq 1$ there exists a constant $c_4 = c_4(G,S,\delta,n_1,n_2,K)$ with the following property:*

*Suppose $U_1 = \langle Y_1 \rangle$ and $U_2 = \langle gY_2g^{-1} \rangle$ are two quasiconvex subgroups of $G$, where $g \in G$, $Y_1 = (y_1,\ldots,y_{n_1}) \in G^{n_1}$, $Y_2 = (y'_1,\ldots,y'_{n_2}) \in G^{n_2}$, $|y_i| \leq K$ for $1 \leq i \leq n_1$ and $|y'_i| \leq K$ for $1 \leq i \leq n_2$. Suppose that $d(X(U_1), X(U_2)) \leq K$.*

*Then the $(n_1+n_2)$-tuple $(y_1,\ldots,y_{n_1}, gy'_1g^{-1},\ldots,gy'_{n_2}g^{-1})$ is Nielsen-equivalent to a tuple conjugate in $G$ to $(y_1,\ldots,y_{n_1+n_2})$ where $|y_i| \leq c_4$ for $1 \leq i \leq n_1+n_2$.*

*Proof.* Note that $X(U_2) = gX(\bar{U}_2)$ where $\bar{U}_2 = \langle Y_2 \rangle$. Put $n = \max(n_1, n_2)$ and choose $c_3 = c_3(G,S,\delta,n,K)$ to be the constant provided by Lemma 10.5. Thus $U_1$ and $X(U_1)$ are $c_3$-Hausdorff close and $\bar{U}_2$ and $X(\bar{U}_2)$ are $c_3$-Hausdorff close. Since $d(X(U_1), X(U_2)) \leq K$ and $X(U_2) = gX(\bar{U}_2)$, there exist $x_1 \in X(U_1)$ and $x_2 \in X(\bar{U}_2)$ such that $d(x_1, gx_2) \leq K$. Since $U_1$ and $X(U_1)$ are $c_3$-Hausdorff close there exist $u_1 \in U_1$ such that $d(u_1, x_1) \leq c_3$. Analogously we find a $u_2 \in \bar{U}_2$ such that $d(u_2, x_2) \leq c_3$.

Put $x'_i = u_i^{-1} x_i$ and $g' = u_1^{-1} g u_2$, so that $d(1, x'_i) \leq c_3$. Note that $x_1 \in X(U_1)$ since $X(U_1)$ is $U_1$-invariant and that $x'_2 \in X(\bar{U}_2)$ since $X(\bar{U}_2)$ is $\bar{U}_2$-invariant. Also

$$d(x'_1, g'x'_2) = d(u_1^{-1}x_1, (u_1^{-1}gu_2)u_2^{-1}x_2) = d(u_1^{-1}x_1, u_1^{-1}gx_2) = d(x_1, gx_2) \leq K$$

and hence

$$d(1, g') \leq d(1, x'_1) + d(x'_1, g'x'_2) + d(g'x'_2, g') \leq c_3 + K + c_3 = 2c_3 + K.$$

It follows that

$$|g'y'_j(g')^{-1}| \leq 4c_3 + 3K.$$

It remains to show that the tuple $M_1 = (y_1, \ldots, y_{n_1}, gy'_1 g^{-1}, \ldots, gy'_{n_2} g^{-1})$ after conjugation is Nielsen equivalent to the tuple $(y_1, \ldots, y_{n_1}, g'y'_1(g')^{-1}, \ldots, g'y'_{n_2}(g')^{-1})$. This would imply that the statement of the lemma follows with $c_4 = 4c_3 + 3K$.

Since $gu_2 g^{-1} \in \langle gy'_1 g^{-1}, \ldots, gy'_{n_2} g^{-1} \rangle$ and since $(gu_2 g^{-1})(gy'_i g^{-1})(gu_2^{-1} g^{-1}) = gu_2 y'_i u_2^{-1} g^{-1}$ it follows that $M_1$ Nielsen equivalent to the tuple $M_2 = (y_1, \ldots, y_{n_1}, gu_2 y'_1 u_2^{-1} g^{-1}, \ldots, gu_2 y'_{n_2} u_2^{-1} g^{-1})$. Moreover, $u_1 \in \langle y_1, \ldots, y_{n_1} \rangle$ implies that $M_2$ is Nielsen equivalent to the tuple

$$M_3 = (y_1, \ldots, y_{n_1}, u_1^{-1} gu_2 y'_1 u_2^{-1} g^{-1} u_1, \ldots, u_1^{-1} gu_2 y'_{n_2} u_2^{-1} g^{-1} u_1).$$

Since $g' = u_1^{-1} g u_2$ the conclusion follows. □

**Lemma 10.7.** *Let $G = \langle S | R \rangle$ be an almost torsion-free word-hyperbolic group. Then for any integers $K \geq 0$ and $n \geq 1$ there exists a constant $c_5 = c_5(G, S, \delta, n, K)$ with the following property:*

*Let $U = \langle Y \rangle$ where $Y = (y_1, \ldots, y_n) \in G^n$ and $|y_i| \leq K$ for $1 \leq i \leq n$. Suppose that $U$ is quasiconvex in $G$ and that $d(X(U), gX(U)) \leq K$.*

*Then the tuple $(y_1, \ldots, y_n, g)$ is Nielsen-equivalent to a tuple $(y_1, \ldots, y_n, y_{n+1})$ such that $|y_i| \leq c_5$ for $1 \leq i \leq n+1$.*

*Proof.* Choose $c_3 = c_3(G, S, \delta, n, K)$ as provided by Lemma 10.5. Hence $U$ and $X(U)$ are $c_3$-Hausdorff-close. Since $d(X(U), gX(U)) \leq K$ there exist points $x_1, x_2 \in X(U)$ such that $d(x_1, gx_2) \leq K$. Since $U$ and $X(U)$ are $c_3$-Hausdorff close there exist $u_1, u_2 \in U$ such that $d(x_1, u_1) \leq c_3$ and that $d(x_2, u_2) \leq c_3$.

Put $y_n = u_1^{-1} gu_2$, $x'_1 = u_1^{-1} x_1$ and $x'_2 = u_2^{-1} x_2$. Since $X(U)$ is $U$-invariant, $x'_1, x'_2 \in X(U)$. Also, $d(x_i, u_i) \leq c_3$ implies that $d(1, x'_i) \leq c_3$ for $i = 1, 2$.

Moreover,

$$d(x'_1, y_n x'_2) = d(u_1^{-1}x_1, (u_1^{-1}gu_2)u_2^{-1}x_2) = d(u_1^{-1}x_1, u_1^{-1}gx_2) = d(x_1, gx_2) \leq K.$$

By the triangle inequality we have

$$|y_n| = d(1, y_n) \leq d(1, x'_1) + d(x'_1, y_n x'_2) + d(y_n x'_2, y_n) \leq 2c_3 + K.$$

The tuple $(y_1, \ldots, y_n)$ is obviously Nielsen-equivalent to $(y_1, \ldots, y_{n-1}, g)$. Hence the statement of Lemma 10.7 holds with $c_5 = 2c_3 + K$. □

We now have all tools for the proof of the main result of this section, which coincides with Theorem 1.3 from the introduction.

**Theorem 10.8.** *Let $G = \langle S | R \rangle$ be an almost torsion-free $\delta$-hyperbolic group, where $\delta > 0$. Suppose $l \geq 0$ is an integer such that all $l$-generated subgroups of $G$ are quasiconvex in $G$ and let $k \geq l + 1$ be any integer. Then there is a constant $C = C(G, S, \delta, k, l) > 0$ with the following property.*

*Suppose $U$ is a non-elementary freely indecomposable subgroup of $G$ generated by the elements $g_1, \ldots, g_k$. Then there exist an element $g \in G$ and a tuple $(f_1, \ldots, f_k)$ Nielsen-equivalent to $(g_1, \ldots, g_k)$ such that $|g^{-1} f_i g|_S \leq C$ for $1 \leq i \leq l + 1$.*

*Proof.* Denote $X = \Gamma(G, S)$, so that $(X, d)$ is $\delta$-hyperbolic. We will prove Theorem 10.8 by induction on $l$. Suppose first $l = 0$. (Note that the only 0-generated subgroup of $G$ is the trivial subgroup which is always quasiconvex). Suppose $U = \langle g_1, \ldots, g_k \rangle \leq G$ is non-elementary and freely indecomposable. Put $H = (g_1, \ldots, g_k)$ and consider the $G$-tuple $M = (; H)$. Then by Theorem 2.4 $M$ is equivalent to $\bar{M} = (; \bar{H})$, where $\bar{H} = (f_1, \ldots, f_k)$, such that one of the cases (1), (2),(3) of Theorem 2.4 applies to $\bar{M}$. Obviously, cases (1) and (2) are not applicable, and hence case (3) occurs. Thus for some $1 \leq i \leq k$ and some $x \in X$ we have $d(x, f_i x) \leq K(k)\delta$, where $K(k)$ is the constant provided by Theorem 2.4. Therefore by the definition of the word-metric on $X$ there is $g \in G$ such that $|g^{-1} f_i g|_S \leq K(k)\delta$. Hence the conclusion of Theorem 10.8 holds for $l = 0$ with $C = C(G, S, \delta, k, 0) = K(k)\delta$.

Suppose now that $l > 0$ and that Theorem 10.8 has been established for all $l' < l$.

Let $K(k)$ be the constant from Theorem 2.4. We define $R(i) = R(G, S, \delta, k, i)$ inductively for $1 \leq i \leq l$. First put $R(1) := \delta K(k)$. Then define $R(i) = R(G, S, k, i)$ for $1 < i \leq l$ as

$$R(i) = \max \left( R(i-1), c_5(G, S, \delta, i-1, R(i-1)), \max_{\substack{p, j \geq 1 \\ p+j=i}} c_4(G, S, \delta, p, j, R(i-1)) \right).$$

Note that $\delta K(k) = R(1) \leq R(2) \leq \cdots \leq R(l)$ by construction.

We say that a partitioned tuple $M = (Y_1, \ldots, Y_n; H)$ with $n \geq 0$ and $Y_i = (y_{i,1}, \ldots, y_{i,n_i})$ for $1 \leq i \leq n$ is *good* if

(1) $n_i \leq l$ for $1 \leq i \leq n$.
(2) For every $i \in \{1, \ldots, n\}$ there exists a $g_i \in G$ such that $|g_i y_{i,j} g_i^{-1}| \leq R(n_i)$ for $1 \leq j \leq n_i$.
(3) $U_i = \langle Y_i \rangle \neq 1$ for $1 \leq i \leq n$.

Denote $\mathbb{N} := \{0, 1, 2, 3, \ldots\}$, the set of non-negative integers.

We define the complexity of a partitioned tuple $M = (Y_1, \ldots, Y_n; H)$ with $H = (h_1, \ldots, h_m)$ to be the pair $(m, n) \in \mathbb{N}^2$. We define an order on $\mathbb{N}^2$ by setting $(m, n) \leq (m', n')$ if $m < m$ or if $m = m'$ and $n \leq n'$. This clearly gives a well-ordering on $\mathbb{N}^2$.

Note that for any tuple $(g_1, \ldots, g_k) \in G^k$ there exists a good partitioned tuple $M$ with the underlying tuple Nielsen-equivalent to that has $(g_1, \ldots, g_k)$ as the underlying tuple, namely, $M = (-; H)$, where $H = (g_1, \ldots, g_k)$.

We now show that the conclusion of Theorem 10.8 hold for

$$C = C(G, S, \delta, k, l) := \max \left( c_5(G, S, \delta, l, R(l)), \max_{1 \leq p, j \leq l} (c_4(G, S, p, j, R(l-1)), C(G, S, \delta, k, l-1) \right).$$

Let $H_0 = (g_1, \ldots, g_k)$ be a tuple generating a non-elementary freely indecomposable subgroup $U$ of $G$, where $k \geq l+1$, $l \geq 0$. Recall that by assumption $G$ is almost torsion-free and that all $l$-generated subgroups of $G$ are quasiconvex in $G$. If $(g_1, \ldots, g_k)$ is Nielsen-equivalent to a tuple of the form $(g'_1, \ldots, g'_{k-1}, 1)$ then the statement of Theorem 10.8 holds for $H_0$ by the inductive hypothesis applied to the tuple $(g'_1, \ldots, g'_{k-1})$. Thus we may assume that no tuple Nielsen-equivalent to $H_0$ contains an entry equal to 1. In particular, this means that for any tuple Nielsen-equivalent to $H_0$ the entries of this tuple are pairwise distinct and non-trivial.

Note that the partitioned tuple $M_0 = (; H_0)$ is good. Choose a partitioned good tuple $M = (Y_1, \ldots, Y_n; H)$ with $H = (h_1, \ldots, h_m)$ of minimal complexity such that the underlying tuple of $M$ is Nielsen-equivalent to $H_0 = (g_1, \ldots, g_k)$. We will show that there exists another partitioned tuple $\bar{M} = (\bar{Y}_1, \ldots, \bar{Y}_{\bar{n}}; H)$ such that the underlying tuple is Nielsen equivalent to $H_0$ and that $\bar{Y}_1 = (\bar{y}_1, \ldots, \bar{y}_{\bar{l}})$ is a $\bar{l}$-tuple with $\bar{l} \geq l + 1$ and $|\bar{y}_i| \leq C(G, S, k, l)$ for $1 \leq i \leq \bar{l}$. This clearly proves Theorem 10.8.

Denote $U_j := \langle Y_j \rangle$ for $1 \leq j \leq n$. The subgroups $U_j$ are non-trivial for $j \geq 1$ since we assume that the partitioned tuple $M$ is good.

Since $U$ is freely indecomposable, by Theorem 2.4 we can find a $G$-tuple $N' = (U'_1, \ldots, U'_n; H')$ with $H' = (h'_1, \ldots, h'_m)$ equivalent to the $G$-tuple $N = (U_1, \ldots, U_n; H)$ and such that one of the three cases from Theorem 2.4 applies to $N'$. By definition of equivalence of $G$-tuples each $U'_j$ is conjugate to $U_j$. Therefore $U'_j$ is generated by a tuple $Y'_j$ where $Y'_j$ has the form

$$Y'_j = b_j Y_j b_j^{-1}$$

for some $b_j \in G$. In particular the resulting partitioned tuple $M' = (Y'_1, \ldots, Y'_n; H')$ is good. Denote $X_j = X(U'_j)$.

Suppose case (1) of Theorem 2.4 applies to $N'$. After re-labeling the indices this means that $d(X_1, X_2) \leq \delta K(k)$. After conjugating the partitioned tuple $M'$ and still denoting the result $M'$ we can assume that $Y'_1 = (y'_{1,1}, \ldots, y'_{1,n_1})$ with $|y'_{1,i}| \leq R(n_1)$ for $1 \leq i \leq n_1$ and $Y'_2 = (y'_{2,1}, \ldots, y'_{2,n_2})$ with $|g_2 y'_{2,i} g_2^{-1}| \leq R(n_2)$ for $1 \leq i \leq n_2$ and some $g_2 \in G$. Recall that since $M'$ is a good tuple, we have $n_1, n_2 \leq l$ and hence the subgroups $U'_1, U'_2$ are quasiconvex in $G$ by the assumption on $G$. Recall that $R(i)$ is non-decreasing and that $R(i) \geq \delta K(k)$.

By Lemma 10.6 the $(n_1+n_2)$-tuple $(y'_{1,1}, \ldots, y'_{1,n_1}, g_2 y'_{2,1} g_2^{-1}, \ldots, g_2 y'_{2,n_2} g_2^{-1})$ is Nielsen-equivalent after a conjugation to a tuple $(y_1, \ldots, y_{n_1+n_2})$ with $|y_i| \leq c_4(G, S, \delta, n_1, n_2, R(\max(n_1, n_2)))$. If $n_1+n_2 \geq l+1$ then this gives the assertion of Theorem 10.8 since $c_4(G, S, \delta, n_1, n_2, R(\max(n_1, n_2))) \leq C$ by the choice of $C$. If $n_1 + n_2 \leq l$ then we can replace the two tuples $Y'_1$ and $Y'_2$ in $M'$ by a conjugate of the tuple $Y = (y_1, \ldots, y_{n_1+n_2})$. This produces a new partitioned tuple

$$M'' = (gYg^{-1}, Y'_3, \ldots, Y'_n; H')$$

whose underlying tuple is Nielsen-equivalent to $H_0$. Since $n_1, n_2 \leq n_1 + n_2 - 1$ and by definition of $R(i)$

$$c_4(G, S, \delta, n_1, n_2, R(n_1 + n_2 - 1)) \leq R(n_1 + n_2),$$

it follows that the partitioned tuple $M''$ is good. This however contradicts the minimality of the complexity assumption since the new partitioned tuple $M''$ clearly has smaller complexity than $M'$.

Suppose now that case (2) of Theorem 2.4 applies to $N'$. After re-labeling and possible conjugation of $N'$ we may assume that we have $d(X_1, h'_1 X_1) \leq \delta K(k)$ where $Y'_1 = (y_{1,1}, \ldots, y_{1,n_1})$ and $|y_i| \leq R(n_1)$ for $i = 1, \ldots, n_1$. Recall that $R(n_1) \geq \delta K(k)$.

Thus by Lemma 10.7 the $(n_1 + 1)$-tuple $(y_{1,1}, \ldots, y_{1,n_1}, h'_1)$ is Nielsen-equivalent to a tuple $(y_1, \ldots, y_{n_1+1})$ with $|y_i| \leq c_5(G, S, \delta, n_1, R(n_1))$ for $1 \leq i \leq n_1 + 1$. We construct a new partitioned tuple $M''$ from $M'$ by replacing $H'$ with the $(m-1)$-tuple $(h'_2, \ldots, h'_m)$ in $M'$ replacing $Y'_1$ by $(y_1, \ldots, y_{n_1+1})$. If $n_1 = l$ this implies the assertion of Theorem 10.8 since $c_5(G, S, \delta, l, R(l)) \leq C$ by the definition of $C$. If $n_1 < l$, it follows that the new partitioned tuple $M''$ is good since $c_5(G, S, \delta, n_1, R(n_1)) \leq R(n_1 + 1)$. Since the new partitioned tuple $M''$ has smaller complexity than $M'$, this yields a contradiction as before.

Suppose now that case (3) of Theorem 2.4 applies for $N'$. Thus after re-labeling the elements of $H'$ we have $h'_1 = ghg^{-1}$ with $|h| \leq \delta K(k)$. Note that $h' \neq 1$ by our assumption on $H_0$.

In this case we construct a new partitioned tuple by removing $h'$ from $H'$ and adding the tuple $Y_{p+1} = (h') = (ghg^{-1})$. The underlying set of the new partitioned tuple $M''$ is clearly Nielsen-equivalent to the original set $H_0$. Moreover, $M''$ is good since $h'$ is conjugate to an element of length at most $\delta K(k) = R(1)$. This again contradicts the minimality assumption on the complexity of the partitioned tuple $M'$.

This completes the proof of Theorem 10.8. □

**Remark 10.9** (A note on computability). *Note that by the result of P.Papasoglu [36], given a finite presentation $G = \langle S|R \rangle$ of a word-hyperbolic group, one can effectively find some integer $\delta > 0$ such that the Cayley graph $\Gamma(G, S)$ is $\delta$-hyperbolic.*

*We observe that the constant $C = C(G, S, \delta, k, l)$ from Theorem 10.8 is effectively computable. Indeed, a careful analysis of the proof of Theorem 2.4 shows that $K = K(k)$ is a recursive function of $k$.*

*There are two points in the proof of Theorem 10.8 where computability of constants is not obvious.*

*First, the computability of the constant $c(n, \delta)$ in Lemma 4.10 depends on Lemma 4.8 and the existence of a uniform non-zero lower bound for the asymptotic translation lengths for non-torsion elements in a $\delta$-hyperbolic group $G = \langle s_1, \ldots, s_n | r_1, \ldots, r_m \rangle$. It is well-known that non-torsion elements in hyperbolic groups have positive asymptotic translation lengths. Moreover, for a fixed finite presentation of a hyperbolic group these length are rational numbers with bounded denominators. As follows from the results of E.Swenson [42] and T.Delzant [22], the smallest denominator $b$ can be effectively computed for a given finite presentation of a hyperbolic group $G$. This implies that $||g||_\infty \geq 1/b$ for any non-torsion element $g \in G$.*

*Second, we need to show that the constant $c_3 = c_3(G, S, \delta, n, C)$ in Lemma 10.5 is computable provided it is known that all $n$-generated subgroups in $G$ are quasiconvex. Computing constant $c_3$, apart from dealing with finite subgroups, essentially involves the following calculation.*

*For a given integer $C$ find the maximum $E$ of quasiconvexity constants among the subgroups of $G$ generated by an $n$-tuple of elements contained in the ball of radius $C$ (provided all $n$-generated subgroups of $G$ are known to be quasiconvex). The number of such $n$-tuples is at most $(m^C)^n = m^{Cn}$ where $m$ is the number of elements in $S$. As was proved by I.Kapovich in [27], there is an algorithm which, given a tuple of elements generating a quasiconvex subgroup, produces the quasiconvexity constant of this subgroup. Thus $E$ is effectively computable. It is worth noting that the existence of $c_3$, provided by Lemma 10.5, is the only place in the proof of Theorem 10.8 where we use the assumption that the centralizers in $G$ are cyclic.*

## 11. Applications to 3-manifolds

In this section we will prove Theorem 1.7 from the introduction. Let $M$ be a closed hyperbolic 3-manifold so that the universal cover $\tilde{M}$ is the hyperbolic space $\mathbb{H}^3$. Let $G = \pi_1 M$ so that $G$ acts on $\mathbb{H}^3$ freely, isometrically and discretely with $\mathbb{H}^3/G = M$. Recall that a finitely generated subgroup $K \le G$ is said to be *topologically tame* if $\mathbb{H}^3/K$ is homeomorphic to the interior of a compact 3-manifold.

*Proof of Theorem 1.7.* If $M$ is as in part (1) of Theorem 1.7 then, as observed by G.Swarup [40], $G$ is torsion-free word-hyperbolic and locally quasiconvex. Thus the statement of Theorem 1.7 follows directly from Theorem 1.4.

Suppose now that $M$ is as in part (2) of Theorem 1.7. Thus $M$ is a closed hyperbolic 3-manifold fibering over a circle and such that all finitely generated subgroups of $G = \pi_1(M)$ are topologically tame. Since $M$ is a closed hyperbolic manifold, the universal cover of $M$ is $\mathbb{H}^3$ and $G$ acts on $\mathbb{H}^3$ freely and discretely by isometries with $\mathbb{H}^3/G = M$. Recall also that by a result of G.Swarup [40] a finitely generated subgroup $K$ of $G$ is quasiconvex in $G$ if and only if $K$ is geometrically finite with respect to the induced action of $H$ on $\mathbb{H}^3$.

We will say that a finitely generated group $K \le G$ is a *virtual fiber group* if $H$ is commensurable with a fiber group of some finite cover of $M$. Clearly, any virtual fiber group is not quasiconvex in $G$. A result of R.Canary [16] implies that under our tameness assumption on $G$ the converse is also true:

*A finitely generated subgroup $K \le G$ is quasiconvex in $G$ if and only if $K$ is not a virtual fiber group.*

Let $k \ge 2$ and suppose $K$ is a freely indecomposable subgroup of infinite index in $G$ such that the rank of $K$ is equal to $k$. That is $K$ can be generated by $k$ elements but cannot be generated by fewer than $k$ elements. Let $Y = \{g_1, \dots, g_k\} \subseteq G$ be a set with $k$-elements such that $\langle g_1, \dots, g_k \rangle = K$. Note that any proper subset of $Y$ generates a quasiconvex subgroup of $G$.

Indeed, suppose first that $K$ is a virtual fiber group, that is $K$ is the fundamental group of a closed hyperbolic surface. Then all subgroups of finite index in $K$ have higher rank than $K$ and therefore any proper subset of $Y$ generates a subgroup of infinite index in $K$. Hence this subgroup is free and so not a virtual fiber group and therefore quasiconvex in $G$.

Suppose that $K$ is not a virtual fiber group, so that $K$ is quasiconvex. If $K$ does not contain a virtual fiber group then all finitely generated subgroups of $K$ are quasiconvex in $G$ and hence any proper subset of $Y$ generates a quasiconvex subgroup of $G$. Suppose now that $K$ contains a virtual fiber group $K'$. Then it is easy to see that $K$ is either commensurable with $K'$ or $K$ has finite index in $G$, both of which contradict our assumptions on $K$.

Exactly the same argument implies that for any tuple $(g'_1, \dots, g'_k)$ Nielsen-equivalent to $(g_1, \dots, g_k)$ any proper subset of $Y' = \{g'_1, \dots, g'_k\}$ generates a quasiconvex subgroup of $G$.

We now define good partitioned tuples exactly as in the proof of Theorem 10.8 and follow precisely the proof of Theorem 10.8. The same argument shows that $H_0 := (g_1, \dots, g_k)$ (applied to $l = k - 1$) shows that $H_0$ is Nielsen-equivalent to a tuple of the form $(gf_1g^{-1}, \dots, gf_kg^{-1})$ where $|f_i| < C = C(G, S, \delta, k, k-1)$ for $i = 1 \dots, k$ (and where $C$ is defined exactly as in the proof of Theorem 10.8). Indeed, in the proof of Theorem 10.8 we did not need the full strength of the assumption that all $l$-generated subgroups of $G$ be quasiconvex. It was sufficient to know that for any tuple $H$ Nielsen-equivalent to $H_0$ any $l$ entries of $H$ generate a quasiconvex subgroup of $G$. This condition is clearly satisfied in the present case for any $l \le k - 1$.

This implies the statement of Theorem 1.7. □

## References

[1] J.Alonso, T.Brady, D.Cooper, V.Ferlini, M.Lustig, M.Mihalik, M.Shapiro and H.Short. *Notes on hyperbolic groups.* in "Group theory from a geometric viewpoint", Proc. ICTP. Trieste, World Scientific, Singapore, 1991, 3–63

[2] G.N.Arzhantseva, *Generic properties of finitely presented groups and Howson's theorem,* Comm. Algebra **26** (1998), no. 11, 3783–3792

[3] G.N.Arzhantseva and A.Yu.Ol'shanskii, *Generality of the class of groups in which subgroups with a lesser number of generators are free* (Russian), Mat. Zametki **59** (1996), no. 4,489–496

[4] W. Ballmann, *Lectures on spaces of nonpositive curvature*, Birkhäuser Verlag, Basel, 1995, With an appendix by Misha Brin.

[5] H. Bass, *Covering theory for graphs of groups*, J. Pure Appl. Algebra **89** (1993), no. 1-2, 3–47.

[6] G. Baumslag, C.F. Miller III and H. Short, *Unsolvable problems about small cancellation and word hyperbolic groups*, Bull. London Math. Soc. **26**, 1994, 91–101.

[7] N. Benakli and I. Kapovich, *Boundaries of hyperbolic groups,* "Combinatorial and Geometric Group Theory" (R.Gilman et al, editors), Contemporary Mathematics, **296**, pp. 39–93, American Mathematical Society, 2002

[8] M. Bestvina and M. Feighn, *A combination theorem for negatively curved groups*, J. Differential Geom. **35** (1992), no. 1, 85–101.

[9] ———, *Addendum and correction to: "A combination theorem for negatively curved groups" [J. Differential Geom. **35** (1992), no. 1, 85–101; MR 93d:53053]*, J. Differential Geom. **43** (1996), no. 4, 783–788.

[10] O. Bogopolskii and V. Gerasimov, *Finite subgroups of hyperbolic groups,* Algebra i Logika **34** (1995), no. 6,619–622,

[11] O. Bogopolski and R. Weidmann *On the uniqueness of factors of amalgamated products*, J. Group Theory **5** (2002), no. 2, 233–240

[12] B. H. Bowditch, *Notes on Gromov's hyperbolicity criterion for path-metric spaces,* in "Group theory from a geometrical viewpoint (Trieste, 1990)", 64–167, World Sci. Publishing, 1991

[13] N. Brady, *Finite subgroups of hyperbolic groups,* Internat. J. Algebra Comput. **10** (2000), no. 4, 399–405

[14] M. Bridson and A. Haefliger, *Metric spaces of non-positive curvature*, Springer-Verlag, Berlin, 1999.

[15] I. Bumagin, *On small cancellation $k$-generated groups with $(k-1)$-generated subgroups all free,* Internat. J. Algebra Comput. **11** (2001), no. 5, 507–524

[16] R. Canary, *A covering theorem for hyperbolic 3-manifolds and its applications,* Topology **35** (1996), no. 3, 751–778

[17] J. W. Cannon, *The theory of negatively curved spaces and groups*, Ergodic theory, symbolic dynamics, and hyperbolic spaces (Trieste, 1989), Oxford Univ. Press, New York, 1991, pp. 315–369.

[18] G. R. Conner, *Translation numbers of groups acting on quasiconvex spaces*, Computational and geometric aspects of modern algebra (Edinburgh, 1998), Cambridge Univ. Press, Cambridge, 2000, pp. 28–38.

[19] M. Coornaert, T. Delzant, and A. Papadopoulos, *Géométrie et théorie des groupes. Les groupes hyperboliques de Gromov*, Lecture Notes in Mathematics **111**, Springer-Verlag, Berlin, 1990, .

[20] T. Delzant, *Sous-groupes a deux generateurs des groups hyperboliques*, in "Group theory from a geometric viewpoint", Proc. ICTP. Trieste, World Scientific, Singapore, 1991, 177–192.

[21] T. Delzant, *L'image d'un groupe dans un groupe hyperbolique*, Comment. Math. Helv. **70**, No.2, 267-284 (1995).

[22] T. Delzant, *Sous-groupes distingus et quotients des groupes hyperboliques,* Duke Math. J. **83** (1996), no. 3, 661–682

[23] D. B. A. Epstein, J. W. Cannon, D. F. Holt, S. V. F. Levy, M. S. Paterson, and W. P. Thurston, *Word processing in groups*, Jones and Bartlett Publishers, Boston, MA, 1992

[24] M. Gromov, *Hyperbolic groups*, Essays in group theory, editor S.M. Gersten, Springer-Verlag, MSRI Publications 8, 1985, 75–263.

[25] É. Ghys and P. de la Harpe (eds.), *Sur les groupes hyperboliques d'après Mikhael Gromov*, Birkhäuser Boston Inc., Boston, MA, 1990, Papers from the Swiss Seminar on Hyperbolic Groups held in Bern, 1988.

[26] C. Hruska and D. Wise *Towers, ladders, and the B.B.Newman spelling theorem,* Jour. Aust. Math. Soc. **71** (2001), no. 1, 53–69

[27] I. Kapovich, *Detecting quasiconvexity: algorithmic aspects,* in "Geometric and computational perspectives on infinite groups (Minneapolis, MN and New Brunswick, NJ, 1994)", DIMACS Ser. Discrete Math. Theoret. Comput. Sci., 25, Amer. Math. Soc., Providence, RI, 1996, 91–99

[28] I. Kapovich, *The Combination Theorem and quasiconvexity,* International Journal of Algebra and Computation, **11** (2001), no. 2, 185–216

[29] M. Kapovich, *Hyperbolic manifolds and discrete groups,* Progress in Mathematics, **183**, Birkhäuser, Boston, 2001

[30] I. Kapovich and P. Schupp, *Bounded rank subgroups of Coxeter groups, Artin groups and one-relator groups with torsion,* preprint, 2002

[31] I. Kapovich and H. Short, *Greenberg's theorem for quasiconvex subgroups of word hyperbolic groups*, Canad. J. Math. **48** (1996), no. 6, 1224–1244.

[32] I. Kapovich and R. Weidmann, *Nielsen Methods for groups Acting on hyperbolic spaces*, to appear in Geom. Dedicata.

[33] I. Kapovich and R. Weidmann, *Acylindrical accessibility for groups acting on $\mathbb{R}$-trees*, preprint.

[34] R. C. Lyndon and P. E. Schupp, *Combinatorial group theory*, Springer-Verlag, Berlin, 1977, Ergebnisse der Mathematik und ihrer Grenzgebiete, Band 89.

[35] J. McCammond and D. Wise *Coherence, Local Quasiconvexity, and the Perimeter of 2-Complexes,* preprint, 1999

[36] P. Papasoglu, *An algorithm detecting hyperbolicity* Baumslag, Gilbert (ed.) et al., Geometric and computational perspectives on infinite groups. Proceedings of a joint DIMACS/Geometry Center workshop, 1994, DIMACS, Ser. Discrete Math. Theor. Comput. Sci. 25, 193-200 (1996).

[37] E. Rips, *Subgroups of small cancellation groups*, Bull. Lond. Math. Soc. **14**, 45-47, 1982.

[38] P. Schupp, *Coxeter groups, perimeter reduction and subgroup separability,*, to appear in Geom. Dedicata

[39] Z. Sela, *Structure and rigidity in (Gromov) hyperbolic groups and discrete groups in rank* 1 *Lie groups. II,* Geom. Funct. Anal. **7** (1997), no. 3, 561–593

[40] G. A. Swarup, *Geometric finiteness and rationality,* J. Pure Appl. Algebra **86** (1993), no. 3, 327–333.

[41] E. L. Swenson, *Hyperbolic elements in negatively curved groups*, Geom. Dedicata **55** (1995), no. 2, 199–210.

[42] ———, *Quasi-convex groups of isometries of negatively curved spaces*, Topol. Appl. **110** (2001), no. 1, 119–129.

[43] R. Weidmann *The Nielsen method for groups acting on trees*, Proc. London Math. Soc. **85** (2002), no. 1, 93–118

Dept. of Mathematics, University of Illinois at Urbana-Champaign, 1409 West Green Street, Urbana, IL 61801, USA
*E-mail address*: `kapovich@math.uiuc.edu`

Fakultät für Mathematik, Ruhr-Universität Bochum, Germany
*E-mail address*: `Richard.Weidmann@ruhr-uni-bochum.de`